\documentclass[11pt]{article}
\usepackage[left=1in, right=1in, top=1in, bottom=1in]{geometry}
\usepackage[table]{xcolor} 
\usepackage{setspace}
\usepackage{rotating}
\usepackage{multirow}
\usepackage[utf8]{inputenc}
\usepackage[T1]{fontenc}
\usepackage{caption}
\usepackage{subcaption}
\usepackage{enumitem}
\usepackage{amsthm}
\usepackage{booktabs}
\usepackage{amsmath}
\usepackage{bbm}
\usepackage{amssymb}
\usepackage{amsfonts}
\usepackage{amsmath}
\usepackage{color}
\usepackage{graphicx}
\usepackage{float}
\usepackage{soul}
\usepackage{algorithm}
\usepackage{bm}
\usepackage{url}
\usepackage{enumerate}
\usepackage{natbib}
\usepackage{comment}
\usepackage{mathtools}
\usepackage{xr}                  
\usepackage{xr-hyper}   
\usepackage{tikz}
\usepackage{multirow}
\usepackage{threeparttable}
\usetikzlibrary{positioning}
\usetikzlibrary{positioning, shapes.geometric}
\numberwithin{equation}{section}
\newtheorem{example}{Example}

\newtheorem{Thm}{Theorem}[section]
\newtheorem{Lem}{Lemma}[section]

\newtheorem{remark}{Remark}[section]

 \def\esssup{\mathop {\rm ess\,sup}}
  
  \def\Z{\bm Z}
  \def\bmtheta{\bm \theta}
  \def\bmvartheta{\bm \vartheta}

\begin{document}

\global\long\def\a{\alpha}%
 
\global\long\def\b{\beta}%
 
\global\long\def\g{\gamma}%
 
\global\long\def\d{\delta}%
 
\global\long\def\e{\epsilon}%
 
\global\long\def\l{\lambda}%
 
\global\long\def\t{\theta}%
 
\global\long\def\o{\omega}%
 
\global\long\def\s{\sigma}%

\global\long\def\G{\Gamma}%
 
\global\long\def\D{\Delta}%
 
\global\long\def\L{\Lambda}%
 
\global\long\def\T{\Theta}%
 
\global\long\def\O{\Omega}%
 
\global\long\def\R{\mathbb{R}}%
 
\global\long\def\N{\mathbb{N}}%
 
\global\long\def\Q{\mathbb{Q}}%
 
\global\long\def\I{\mathbb{I}}%
 
\global\long\def\P{\mathbb{P}}%
 
\global\long\def\E{\mathbb{E}}%
\global\long\def\B{\mathbb{\mathbb{B}}}%
\global\long\def\S{\mathbb{\mathbb{S}}}%

\global\long\def\X{{\bf X}}%
\global\long\def\cX{\mathscr{X}}%
 
\global\long\def\cY{\mathscr{Y}}%
 
\global\long\def\cA{\mathscr{A}}%
 
\global\long\def\cB{\mathscr{B}}%
 
\global\long\def\cM{\mathscr{M}}%
\global\long\def\cN{\mathcal{N}}%
\global\long\def\cG{\mathcal{G}}%
\global\long\def\cC{\mathcal{C}}%
\global\long\def\sp{\,}%

\global\long\def\es{\emptyset}%
 
\global\long\def\mc#1{\mathscr{#1}}%
 
\global\long\def\ind{\mathbf{\mathbbm1}}%
\global\long\def\indep{\perp}%

\global\long\def\any{\forall}%
 
\global\long\def\ex{\exists}%
 
\global\long\def\p{\partial}%
 
\global\long\def\cd{\cdot}%
 
\global\long\def\Dif{\nabla}%
 
\global\long\def\imp{\Rightarrow}%
 
\global\long\def\iff{\Leftrightarrow}%

\global\long\def\up{\uparrow}%
 
\global\long\def\down{\downarrow}%
 
\global\long\def\arrow{\rightarrow}%
 
\global\long\def\rlarrow{\leftrightarrow}%
 
\global\long\def\lrarrow{\leftrightarrow}%
\global\long\def\gto{\rightarrow}%

\global\long\def\abs#1{\left|#1\right|}%
 
\global\long\def\norm#1{\left\Vert #1\right\Vert }%
 
\global\long\def\rest#1{\left.#1\right|}%

\global\long\def\bracket#1#2{\left\langle #1\middle\vert#2\right\rangle }%
 
\global\long\def\sandvich#1#2#3{\left\langle #1\middle\vert#2\middle\vert#3\right\rangle }%
 
\global\long\def\turd#1{\frac{#1}{3}}%
 
\global\long\def\ellipsis{\textellipsis}%
 
\global\long\def\sand#1{\left\lceil #1\right\vert }%
 
\global\long\def\wich#1{\left\vert #1\right\rfloor }%
 
\global\long\def\sandwich#1#2#3{\left\lceil #1\middle\vert#2\middle\vert#3\right\rfloor }%

\global\long\def\abs#1{\left|#1\right|}%
 
\global\long\def\norm#1{\left\Vert #1\right\Vert }%
 
\global\long\def\rest#1{\left.#1\right|}%
 
\global\long\def\inprod#1{\left\langle #1\right\rangle }%
 
\global\long\def\ol#1{\overline{#1}}%
 
\global\long\def\ul#1{\underline{#1}}%
 
\global\long\def\td#1{\tilde{#1}}%

\global\long\def\upto{\nearrow}%
 
\global\long\def\downto{\searrow}%
 
\global\long\def\pto{\overset{p}{\longrightarrow}}%
 
\global\long\def\dto{\overset{d}{\longrightarrow}}%
 
\global\long\def\asto{\overset{a.s.}{\longrightarrow}}%

\global\long\def\tJ{\tilde{{\cal J}}}%

\title{Optimization via Strategic Law of Large Numbers}







\author{  Xiaohong Chen\\Cowles Foundation for Research in Economics, Yale University \and Zengjing Chen\\Zhongtai Securities Institute for Financial Studies, Shandong University\and Wayne Yuan Gao   \\Department of Economics, University of Pennsylvania \and Xiaodong Yan\\School of Mathematics and Statistics, Xi'an Jiaotong University \and Guodong Zhang \\School of Statistics and Mathematics, Shandong University of Finance and Economics }
\maketitle

\spacing{1.3}

\begin{abstract}
 This paper proposes a novel framework for the global optimization of a possibly multi-modal continuous function in a bounded rectangular domain. We first show that global optimization is equivalent to an optimal (sampling) strategy formation in a two-armed decision model with known distributions, based on the Strategic Law of Large Numbers we establish. There are many optimal strategies in general. We show that a concrete strategy using the sign of the partial gradient of the unique solution to a parabolic PDE is asymptotically optimal. Motivated by these results, we propose a class of {\bf S}trategic {\bf M}onte {\bf C}arlo {\bf O}ptimization (SMCO) algorithms, which uses a simple strategy that makes coordinate-wise two-armed decisions based on the signs of the partial gradient (or practically the first difference) of the objective function, without the need of solving PDEs. Under some sufficient conditions, we establish that our SMCO algorithm converges to a local optimizer from a single starting point, and to a global optimizer under a growing set of starting points. Extensive numerical studies demonstrate the suitability of our SMCO algorithms for global optimization well beyond the theoretical guarantees established herein. For a wide range of deterministic and random test functions with challenging landscapes (multi-modal, non-differentiable, discontinuous), our SMCO algorithms perform robustly well, even in high-dimensional ($d=200\sim1000$) settings. In fact, our algorithms outperform many state-of-the-art global optimizers, as well as local algorithms (with the same set of starting points as ours).\\
 \\
\textbf{Keywords:} Black box global optimization, multi-modal criterion, high dimension, two-armed machine, strategic sampling, strategic limit theorem, HJB equation, stochastic control
\end{abstract}
\spacing{1.5}





\section{Introduction}
It has been a notoriously difficult task to develop fast algorithms for global optimization of 
high-dimensional and potentially multi-modal continuous functions, which is known as the “curse of optimality” phenomenon. In this paper we propose a novel approach to solving the following optimization problem
\begin{equation}\label{opt0}
\max_{x\in\Gamma}f(x),
\end{equation}
where $\Gamma :=[\underline{\mu}_1,\overline{\mu}_1]\times\cdots\times[\underline{\mu}_d,\overline{\mu}_d] \subset \mathbb{R}^d$, and $f:\Gamma \to \mathbb{R}$ is a continuous function.

Problem \eqref{opt0} is well-defined in the sense that there exists a global maximizer of any continuous deterministic function $f$ over a compact region $\Gamma \subset \mathbb{R}^d$. Under mild regularity conditions, it is known that a continuous random function also has a global maximizer on a compact region $\Gamma$ almost surely without assuming $f$ being a global concave function (see, for example, Cox \cite{Cox2020}). However, when $f$ is a multi-modal and high-dimensional continuous function, such as Rastrigin function with $d=200,~400,~500,~1000$, it is very difficult to develop a fast algorithm to solve Problem \eqref{opt0} and to mathematically prove its convergence to a global maximizer. Various Stochastic Gradient Descent (SGD) algorithms rely on accurate estimation of gradients (some also rely on Hessians) for iterative updates, and are the key behind recent success of deep learning (neural networks). The gradient based algorithms are only shown to converge to a global maximum when the objective function $f$ is strictly concave, however. There are also many popular meta-heuristic global optimization algorithms for Problem \eqref{opt0}, such as Simulated Annealing (SA) of Kirkpatrick et al. \cite{Kirkpatrick1983} and the Particle Swarm Optimization (PSO) of Eberhart and Kennedy \cite{EK1995}. SA has been shown to converge weakly to global optimal points (Geman and Hwang \cite{GemanHwang1986}) but its convergence rate is known to be slow and it performs poorly in practice. Although popular in applied work, the PSO has not been mathematically proved to converge to a global maximum for a general multi-modal continuous function yet. 

We propose a novel framework for the global optimization problem \eqref{opt0}. Our core theoretical results, Theorems \ref{thm-lln-01} and \ref{thm-convergence2-sup}, establish that solving problem \eqref{opt0} is equivalent to the search for an optimal (sampling) strategy in a two-armed decision model with known distributions, based on the Strategic Law of Large Numbers (StLLN). Precisely, the StLLN we establish provides a new theoretical framework that translates the global optimization problem \eqref{opt0} into an equivalent problem of finding an approximately optimal strategy in a two-armed bandit model with known distributions. Both StLLN and the classical LLN concern the asymptotic behavior of the sample mean, but StLLN considers strategic sampling from two distributions, while the classical LLN focuses on (nonstrategic) sampling from a single distribution. Hence, the classical LLN can be viewed as a special case of StLLN in which the two distributions coincide. 

There are in general many approximately optimal strategies to solve problem \eqref{opt0}. Our Theorem \ref{thm-pde} further establishes that one approximately optimal strategy is to use the signs of partial gradient of the unique solution to a parabolic PDE. For low-dimensional (e.g., $d=1,~2$) problems, one could apply finite-difference method to numerically solve the PDE, and then design a two-armed sampling algorithm 
according to the sign of the partial gradient of the numerical solution of the PDE. Unfortunately, solving a high-dimensional PDE numerically is in general  very difficult. Interestingly, under symmetry and some smoothness condition on a one-dimensional $f$, the sign of gradient of the unique solution to the PDE coincides with the sign of the gradient of the original $f$ itself (see Remark \ref{remark-pde}).

Motivated by the above results, we propose a simple optimization algorithm based on a tractable strategy for the two-armed decision model, using only the sign information of the partial derivative of the original function $f$ (instead of the PDE solution). The running sample mean under iterative application of this strategy is then returned as the candidate optimizer by our algorithm. We call this the ``\textbf{S}trategic \textbf{M}onte \textbf{C}arlo \textbf{O}ptimization'' (SMCO) algorithm. Although this simple strategy may not coincide with the asymptotic optimal strategy using the sign of the gradient of the PDE solution in multi-dimensional problems, our Theorem \ref{thm-lln} shows that the SMCO algorithm is nevertheless guaranteed to converge to \textit{local} optimizers. 
We show in Theorem \ref{thm-ms} that the SMCO algorithm, when augmented with a growing set of starting points, converges to a global optimizer. In addition, when $f$ has a unique local maximizer $x^*$ on $\Gamma$, Theorem \ref{thm-convergence2} establishes a convergence rate of our simple strategic sample mean to $x^*$. In practice, we implement the SMCO algorithms with multiple randomly generated starting points, and use step-adaptive finite differences in place of the true gradients. We also provide two slightly modified versions of the basic SMCO algorithm, SMCO-R and SMCO-BR, with bookkeeping of running maximum and simple additional programming enhancements. The source code of SMCO/-R/-BR is publicly available on GitHub.

We conduct a wide range of numerical experiments to analyze the performance of our SMCO algorithms in solving complex global optimization problems in comparison with existing popular optimizers in R. We consider various multi-modal deterministic test functions widely used as benchmarks in the global optimization literature, such as Rastrigin, Ackely, Michalewicz, Griewank, and more, of dimensions ranging from 1 to 1000. In addition, we also consider a variety of random test functions with challenging optimization landscapes (some non-differentiable and even discontinuous), such as sample criterion function for maximum score estimation of Manski \cite{manski1975maximum}, empirical welfare function for optimal treatment assignment, conditional moment inequalities in dynamic binary choice models, in-sample approximation and out-of-sample prediction based on ReLU neural networks in regression and classification settings, and Cauchy likelihood functions. We compare our SMCO algorithms with a range of widely used state-of-the-art optimization algorithms/packages implemented in R, including six local algorithms (GD, SignGD, ADAM, SPSA, L-BFGS, and BOBYQA) and six global ones (GenSA, SA, DEoptim, STOGO, GA, and PSO). 

Our SMCO algorithms, designed using the simple sign-based strategy, perform remarkably well in complex global optimization problems, to an extent that goes well beyond the theoretical guarantees we establish.
Specifically, our SMCO algorithms are consistently able to solve almost all test optimization problems to a reasonable degree of accuracy, using only a very limited number of randomly generated starting points for dimensions up to $d=1000$. In contrast, all other algorithms considered here, with the sole exception of GenSA (an improved version of SA), suffer from drastic errors in at least a number of test cases. In particular, the six local optimizers perform significantly worse than ours in many test cases, even though they use the same set of starting points as ours. This provides strong evidence that the performance of our SMCO algorithms in global optimization cannot be explained by the use of multiple starting points alone. In addition, our SMCO algorithms compare favorably to most of the global optimizers, especially so for higher dimensional test functions. This is again remarkable given that our SMCO algorithms only use a small number of starting points even in higher dimensions. Although GenSA sometimes perform slightly better than ours, it is extremely slow as dimension increases, which explains why GenSA is rarely used in modern high dimensional optimization problems. The performance of our SMCO algorithms thus appears unexpectedly strong, especially given that all of these competing optimizers have been widely tested and improved over many years, and that there are likely to be many aspects of the SMCO algorithms and our simple R code that can be significantly improved (such as coding in C++). 

The numerical experiments demonstrate the promise of our SMCO algorithms, and more generally the optimization framework we proposed based on StLLN with two-armed decision strategy. We note that the SMCO algorithms are based on a simple sign-based strategy only, although robust to multi-modal landscape and noisy gradient, they may not be as fast as some other algorithms when $f$ is strictly concave and differentiable. This suggests an array of future research questions under our two-armed StLLN framework: to analyze the performances of other strategies under this framework; to search for and/or implement other (approximately) optimal strategies in a tractable way; or to provide additional theoretical results that articulate the observed performance of our SMCO algorithms in solving complex global optimization problems.  

\noindent{\bf Related Literature:}

(i) {\bf Gradient Descent and its variants}: Numerous optimization algorithms are based on gradient descent (GD) and stochastic gradient descent (SGD), including sign gradient descent (SignGD) \cite{Bern2018,moulay2019properties,safaryan2021stochastic} and adaptive gradient \cite{Kingma2014,Hinton2012,Zeiler2012,Reddi2019}. Unlike traditional GD, which adjusts based on gradient magnitude, SignGD focuses solely on the gradient's sign for updates. Li et al. \cite{li2023faster} showed that SignGD can improve stability and convergence in some cases, but its performance, like traditional GD, remains highly sensitive to the choices of step sizes and initial points. Guo and Fu \cite{GuoFu2022} proposed to apply multi-arm bandit to select initial starting points for SGD and demonstrated its convergence advantages numerically. Adaptive GD methods automatically adjust step sizes using past stochastic gradients, which potentially lead to faster convergence and mitigate gradient explosion. However, these methods still face challenges, particularly in high-dimensional settings, where step size determination is complex. Additionally, they are prone to gradient explosion in early training phases \cite{Bengio1994,Pascanu2013}, causing slow convergence or divergence. Our SMCO may look like SignGD in that both use sign information, but is quite different in that our sample mean formulation based on StLLN automatically adapts step sizes with respect to both the iteration counter and the current optimizer position, with randomness incorporated in the mean-zero distribution around the upper/lower bounds. We demonstrate in simulations how our algorithm compares favorably with GD, SignGD, and its variants in terms of optimization performances, especially in optimization problems with challenging landscapes.

(ii) {\bf Decision model-driven optimization:} Recently, bandit process, a simplified decision model of reinforcement learning, has shown significant success in optimization by employing random strategy with exploration, because this exploration expands the search space, enhances flexibility, and supports learning through trial and error. See \cite{Gittins2011,li2018hyperband,Lattimore2020,GuoFu2022} and the references therein for various applications of bandit algorithms with unknown arm distributions in optimization. Nevertheless, to the best of our knowledge, our work is the first to formally establish the equivalence of global optimization problem \eqref{opt0} to the optimal policy/strategy formation in a two-armed decision model with known arm distributions. We further propose a simple yet novel SMCO algorithm with multi-start, and rigorously establish its convergence to the global maximizer(s). The StLLN for the constructed decision model (Theorem \ref{thm-lln-01}) is our theoretical foundation for recasting the global optimization of a general continuous function as a decision rule problem, which is inspired by the LLN under the nonlinear expectation of Peng \cite{peng2019}.

(iii) \textbf{Two-armed bandit problems}: Our optimization framework is reminiscent of, but in fact distinct from, the classical two-armed bandit problems. While traditional bandit problems focus on maximizing cumulative rewards under uncertainty about arm distributions \cite{robbins,Lai-Robbins1985,Agrawal1995,Bubeck2012}, we search for global maximizer(s) of a possibly non-concave continuous function $f$ over $\Gamma$ according to a two-armed decision model with known arm distributions. Precisely, SMCO selects between two fixed arms (or $2\times d$ for $d$-dimensional problems) based on the sign changes of the objective function $f$'s coordinate-wise partial gradient, ensuring quicker convergence to maximum points while reverting to standard Monte Carlo sampling when the arms coincide. While the standard Gradient Descent method is sensitive to the choice of starting points, the updating step-size, and the accuracy of gradient estimates for iterative updates, our SMCO relies on two-armed decision rule-generated samples, in which the selection of arms is determined by the signs of coordinate-wise partial gradient, while the update step size is automatically adapted through the sample mean process. 

\noindent{\bf Outline of the rest of the paper:}

Section \ref{sec-theory} presents the Strategic Law of Large Numbers, the corresponding optimization framework, our SMCO algorithm design, and theoretical guarantees on the convergence of the SMCO algorithm. 
Section \ref{sec:Numerical} conducts a wide range of numerical comparisons with the existing popular optimization algorithms. We conclude with Section \ref{sec:end}. The SI Appendix contains supplemental materials, consisting of all the proofs, technical lemmas,  the details on the practical implementations of the SMCO algorithms, and additional numerical experiments.

\section{Theory and Algorithms}\label{sec-theory}
\subsection{New Framework: Optimization via Strategic Law of Large Numbers }
We use a simple decision model to illustrate our new approach to the optimization problem \eqref{opt0}. Suppose that there are $d$ two-armed decision models: these $d$ decision models are independent of each other and differ only in the mean value corresponding to their arms. For $j=1,\cdots,d$, machine $j$ has two arms with populations $X_j$ and $Y_j$ and different known means $\overline{\mu}_{j}$ and $\underline{\mu}_{j}$ (without loss of generality, we assume $\overline{\mu}_{j}> \underline{\mu}_{j}$). 
One can draw one $d$-dimensional sample from these $d$ two-armed decision models at each round; that is, the $d$ components of the sample come from the $d$ two-armed decision models, and decision model $j$ can be drawn from either arm $X_j$ or arm $Y_j$, but not both. If $X_j$ and $Y_j$ are identical, the decision model $j$ has a single arm and it becomes the classical Monte Carlo draws. 

\smallskip

We will focus on the global optimization of a continuous function $f$ on a bounded
rectangular region $\Gamma$,
    $\Gamma:=[\underline{\mu}_{1},\overline{\mu}_1]\times[\underline{\mu}_{2},\overline{\mu}_2]\times\cdots\times[\underline{\mu}_{d},\overline{\mu}_d]\subset \mathbb{R}^d.$
For a small positive $\delta$, we assume that the function $f$ is well-defined on the $\delta$-extension region  $\Gamma_{\delta}$ of $\Gamma$ defined by
    $\Gamma_{\delta}:=[\underline{\mu}_{1}-\delta,\overline{\mu}_1+\delta]\times[\underline{\mu}_{2}-\delta,\overline{\mu}_2+\delta]\times\cdots\times[\underline{\mu}_{d}-\delta,\overline{\mu}_d+\delta].$
In this paper all random variables are defined on a given probability space $(\Omega,\mathcal{F},P)$. Consider $d$ two-armed decision models, for $1\leq j \le d$, let $(X_j, Y_j)$ denote the random reward of one of these  $d$ two-armed decision models. For simplicity, we assume that $\{X_{j}: 1\le j\le d\}$ and $\{Y_{j}: 1\le j\le d\}$ satisfy
 \begin{equation}\label{denote-xy} 
     X_j=\overline{\mu}_j +\xi_j, \quad Y_j=\underline{\mu}_{j}+\xi'_j,
 \end{equation}
 where $\xi_j$ is a bounded random variable with mean $0,$ variance $\sigma^2$ and bounded by $\delta>0$ (that is, $|\xi_j|\leq \delta$), $\xi'_j$ is an independent copy of $\xi_j$. For all $1\leq j \le d$, let $\{X_{i,j}:i\ge1\}$ and $\{Y_{i,j}:i\ge1\}$ respectively denote the outcomes from the arm $X_{j}$ and arm $Y_{j}$ at the $i{\rm th}$ round. We define a sampling strategy $\bmtheta$ as a sequence of (binary) random vectors $\bmtheta:=\{\bmvartheta_1,\cdots, \bmvartheta_n,\cdots \}$,
where $\bmvartheta_{i}:=(\vartheta_{i,1},\cdots,\vartheta_{i,d})\in \{0,1\}^d$ and $\vartheta_{i,j}=1$ (respectively, $\vartheta_{i,j}=0$) means in the $i{\rm th}$ round we pull the arm $X_{j}$ (respectively, $Y_{j}$) of the $j{\rm th}$ slot machine.
Let $\{\Z_{i}^{\bmtheta}:=(Z_{i,1}^{\bmtheta},\cdots,Z_{i,d}^{\bmtheta}):i\ge1\}$ be a sequence of $d$-dimensional random vectors denoting the $i{\rm th}$  reward  obtained from these $d$ slot machines  related to the strategy $\bmtheta$.
 That is, for $j=1,\cdots,d$,
\begin{equation}\label{eq-1}
Z_{i,j}^{\bmtheta}:=\vartheta_{i,j} X_{i,j}+(1-\vartheta_{i,j})Y_{i,j}=\left\{\begin{array}{ll}
X_{i,j}, &\text{ if }\ \vartheta_{i,j}=1,\\
Y_{i,j}, &\text{ if }\ \vartheta_{i,j}=0.
\end{array}\right.
\end{equation} 
We call a sampling strategy $\bmtheta=\{\bmvartheta_1,\cdots, \bmvartheta_n,\cdots \}$ {\it admissible} if $\bmvartheta_n$ is $\mathcal{H}_{n-1}^{\bmtheta}$-measurable for all $n \geq 1$, where
\begin{equation*}
\mathcal{H}_{n-1}^{\theta }=\sigma \{\eta,\Z_{1}^{\bmtheta },\cdots ,\Z_{n-1}^{\bmtheta
},\bmvartheta _{1},...,\bmvartheta _{n-1}\}\text{ for }n>1\text{, and }\mathcal{H}%
_{0}^{\bmtheta }=\sigma\{\eta \}\text{,}
\end{equation*}  
where $\eta:=(\eta_1,\cdots,\eta_d)\in \Gamma$ denotes the initial point, which could be non-random or randomly drawn from a known distribution (such as uniform distribution on the domain or a deterministic grid distribution on the hyperdiagonal of $\Gamma$).
Let $\Theta$ denote the collection of all admissible sampling strategies. 
For any strategy $\bmtheta \in \Theta$, we let
\begin{equation}\label{S_n}
S_{0}^{\bmtheta}=\eta~,~~\text{ and for all}~ n\geq 1,~ S_{n}^{\bmtheta}=S_{n-1}^{\bmtheta} +\Z_{n}^{\bmtheta}.
\end{equation}
 
\begin{Thm}[\textbf{Strategic Law of Large Numbers and a Two-sided Framework for Optimization }]\label{thm-lln-01}
 Let the $d$ two-armed decision models $\{(X_j,Y_j):1\le j\le d\}$ satisfy \eqref{denote-xy}. Let $\Z_{n}^{\bmtheta}$ and $S_{n}^{\bmtheta}$ be defined in \eqref{eq-1} and \eqref{S_n}. Let $f$ be any continuous function defined on $\Gamma_{\delta}$. Then: the global optimization problem \eqref{opt0} is equivalent to the problem of asymptotically maximizing the ``reward'' with $f$
 \begin{equation}\label{lln-1}
\max\limits_{x\in\Gamma}f(x)=\lim\limits_{n\rightarrow\infty}\sup_{\bmtheta\in\Theta}E\left[f\left(\frac{S_n^{\bmtheta}}{n}\right) \right].
 \end{equation}
Furthermore, we have,
   \begin{align}\label{lln-non-prob}
\sup_{\bmtheta\in\Theta}P\left(\limsup_{n\to\infty}f\left(\frac{S_{n}^{\bmtheta}}{n}\right)=\max_{x\in\Gamma}f(x)\right)=1
\end{align}  
\end{Thm}

According to Theorem \ref{thm-lln-01}, the global optimization \eqref{opt0} can be transformed into two-armed decision model problem, which can be equivalently expressed as
   \begin{equation}\label{lln-1-regret}
\lim\limits_{n\rightarrow\infty}\inf_{\bmtheta\in\Theta}R_{f}^{\bmtheta}(n)=0,
\end{equation} 
where $R_{f}^{\bmtheta}(n):=\max_{x\in\Gamma}f(x)-E\left[f\left(\frac{S_n^{\bmtheta}}{n}\right) \right]$ denote the ``averaged regret'' with $f$ under $\bmtheta$. Therefore, to solve the global optimization problem  \eqref{opt0}, one just needs to find an asymptotically \textit{optimal (sampling) strategy} $\bmtheta^{o} \in \Theta$ such that 
$\lim\limits_{n\rightarrow\infty}R_{f}^{\bmtheta^{o}}(n)=0.$
In addition, if the objective function $f$ is Lipschitz continuous, we can obtain the following speed of convergence result.

\begin{Thm}[\textbf{Convergence Rate}]\label{thm-convergence2-sup} Under conditions of 
Theorem \ref{thm-lln-01}, if $f$ is Lipschitz continuous on $\Gamma_{\delta}$ (with a Lipschitz constant $L_f$), 
then:
   \begin{align}
\left|\sup_{\bmtheta\in\Theta}E\left[f\left(\frac{S_n^{\bmtheta}}{n}\right) \right]-\max_{x\in\Gamma}f(x)\right|\leq L_f\left(\frac{\sigma\sqrt{d}+\sum_{j=1}^d\frac{|\overline{\mu}_j-\underline{\mu}_j|}{2}}{\sqrt{n}}+\frac{\sqrt{E[\|\eta\|^2]}}{n}\right)\label{0convergence-rate-1-sup}.
   \end{align}
 \end{Thm}

\begin{remark}
    In our decision model, when these two arms (populations) coincide, our strategic laws of large numbers (StLLN) (Theorems \ref{thm-lln-01} and \ref{thm-convergence2-sup}) degenerates into the classical LLN, which serves as the theoretical foundation for the classical Monte Carlo simulation method. Consequently, we have the following viewpoints: repeatedly operating a single-armed decision model corresponds to the conclusions drawn from the classical laws of large numbers, and can thus be utilized for simulating the integration of a function. However, strategically operating a two-armed decision model aligns with our proposed strategic laws of large numbers, enabling its application in simulating the extremum of functions.

\end{remark}
Theorems \ref{thm-lln-01} and \ref{thm-convergence2-sup} establish that any optimization problem \eqref{opt0} can be solved by finding some optimal (sampling) strategy of the two-armed decision model.\footnote{For other recent studies on nonlinear and strategic limit theory for decision model problems, the reader is referred to \cite{chenepsteinzhang,chenepsteinzhang2} and references therein.}  Both results are significant and broadly applicable, and motivate us to consider a class of optimization algorithms that compute sample means over draws from one of two pre-specified distributions around the two boundary points (in each dimension), along with a potentially optimal strategy to choose from the two distributions.

\subsection{A Concrete Asymptotically Optimal Strategy}\label{sec-pde-strategy}

In the next theorem, we will suggest an asymptotically optimal (sampling) strategy using the signs of the partial gradients of the solutions to some parabolic PDEs associated with the Lipschitz continuous function $f$.
This PDE approach enables us to find the global optima of function $f$ on region $\Gamma$ via a two-armed decision model.

For any small enough $\epsilon>0$, consider the following parabolic PDE
 \begin{equation}
\left\{
\begin{array}{l}
\displaystyle\partial _{t}u_{\epsilon}(t,x)+ \sup_{p\in\Gamma}\left[D_{x}u_{\epsilon}(t,x)\cdot p\right] +\frac{\epsilon^2}{2} \textrm{tr}\left[D_{x}^2u_{\epsilon}(t,x) \right]=0,\quad (t,x)\in \lbrack 0,1+\epsilon)\times \mathbb{R}^d, \\
u_{\epsilon}(1+\epsilon,x)=\hat{f}(x),%
\end{array}%
\right.   \label{pde-1}
\end{equation}
where $\hat{f}$ is an extension of $f$ defined on $\mathbb{R}^d$, $D_x u_{\epsilon}$ and $D_{x}^2 u_{\epsilon}$ denote the gradient vector and the Hessian matrix of the function $x\to u_{\epsilon}(t,x)$ respectively. Under the bounded Lipschitz continuity assumption on $\hat{f}$, this PDE \eqref{pde-1} is actually the HJB equation, and has a unique solution $u_{\epsilon}$ (see Yong and Zhou \cite[Theorem 5.2, Ch. 4]{yongzhou1999}). Moreover, the solution  $u_{\epsilon}(t,x)$ is continuously differentiable in $x$ (see Peng \cite[Theorem C.4.5]{peng2019}, also see Appendix \ref{append-A} for details). We can establish the following refinement of Theorem  \ref{thm-convergence2-sup}.
 
\begin{Thm}\label{thm-pde} Let $f$ be a Lipschitz continuous function on $\Gamma_{\delta}$ with Lipschitz constant $L_{f}$, and $\hat{f}$ be a bounded  Lipschitz continuous (with Lipschitz constant $L_{f}$) extension of $f$ defined on $\mathbb{R}^d$.
    For any small $\epsilon>0$ let $u_{\epsilon}$ be the unique solution of PDE \eqref{pde-1}. For any $n\in\mathbb{N}$, define the strategy $\hat\bmtheta^{n,\epsilon}=(\hat\bmvartheta_1^{n,\epsilon},\cdots,\hat\bmvartheta_i^{n,\epsilon},\cdots)$ with $\hat\bmvartheta_i^{n,\epsilon}=(\hat\vartheta_{i,1}^{n,\epsilon},\cdots,\hat\vartheta_{i,d}^{n,\epsilon})$ as follows, for $j=1,...,d,$
\begin{equation}\label{strategy-e-optimal}
\hat\vartheta_{i,j}^{n,\epsilon}=\left\{\begin{array}{cc}
    1, & \displaystyle\textup{ if }\partial_{x_j}u_{\epsilon}\left(\frac{i}{n},\frac{S_{i-1}^{\hat\bmtheta^{n,\epsilon}}}{n}\right)>0, \\
   0,  & \textup{otherwise}.
\end{array}\right., \text{ for }1\le i\leq n, \text{ and }\hat\vartheta_{i,j}^{n,\epsilon}=1 \text{ for }i> n.
\end{equation}
Then: $\{\hat{\bmtheta}^{n,\epsilon}\}$ is \textbf{asymptotically $\epsilon$-optimal} in the sense that 
  \begin{align}
\left|E\left[f\left(\frac{S_n^{\hat{\bmtheta}^{n,\epsilon}}}{n}\right)\right]-\max_{x\in\Gamma}f(x)\right|\leq& L_f\epsilon(2C_{\Gamma}+\sqrt{\epsilon}+\sqrt{1+\epsilon})+\tfrac{C_1E[\|\eta\|]}{n}\nonumber\\
&+C\left(\tfrac{1}{n}\sum_{i=1}^d \left(\sigma^2+\tfrac{(\overline{\mu}_i-\underline{\mu}_i)^2}{4}\right)+\tfrac{\sum_{i=1}^d|\overline{\mu}_i-\underline{\mu}_i|}{2n}+2\epsilon\sqrt{\tfrac{d}{n}}+\epsilon^2\right),
    \end{align}
    where $C_{\Gamma}>0$ is a finite constant depending on the boundary of $\Gamma$, and $C,C_1>0$ are finite constants depending on $f$.
\end{Thm}

According to Theorem \ref{thm-pde}, we can use the sign of the gradient vector $D_x u_{\epsilon}$ to construct a sequence of strategies $\hat\bmtheta^{n,\epsilon}$ such that $f(S_n^{\hat\bmtheta^{n,\epsilon}}/n)$ approximates $\max_{x\in\Gamma}f(x)$. 

\begin{remark}\label{remark-pde}
To construct a sequence of asymptotically $\epsilon$-optimal strategies, we only need to know the sign of each component of the gradient vector $D_x u_{\epsilon}$. 
When $f$ and its extension $\hat{f}$ satisfy some assumptions, the sign of the gradient vector $D_x u_{\epsilon}$ can be denoted by the sign of the gradient $f$. For example, 
when $\hat{f}\in C^3(\mathbb{R})$, assume that for any $t\in[0,1]$, the nodal set of $\partial_x u_{\epsilon}(t,x)$ is denoted by $D=\{d_k,k=0,\pm1,\pm2,\cdots\}$. Moreover, if $\hat{f}$ is strictly monotone on $[d_k,d_{k+1}]$, and the
derivatives $\hat{f}^{(k)}\ (k=0,1,2,3)$ have at most polynomial growth, then (see more details in \cite{chenfengliu,LiuSun2025}): 
$$
\left\{\partial_{x}u_{\epsilon}\left(\frac{i}{n},\frac{S_{i-1}^{\hat\bmtheta^{n,\epsilon}}}{n}\right)>0\right\}\text{ is equivalent to }\left\{\partial_{x}\hat{f}\left(\frac{S_{i-1}^{\hat\bmtheta^{n,\epsilon}}}{n}\right)>0\right\}.
$$
Hence the optimal strategy defined in \eqref{strategy-e-optimal} can be replaced by the sign of the gradient vector $D_x \hat{f}$. A sufficient condition for the above assumption to hold is that $\hat{f}\in C^3(\mathbb{R})$ is a  symmetric function with symmetric points $D=\{d_k,k=0,\pm1,\pm2,\cdots\}$ and strictly monotone on $[d_k,d_{k+1}]$.
\end{remark}

It is well-known that obtaining an analytical solution for PDE \eqref{pde-1} is generally intractable, often necessitating numerical methods to approximate the sign of its solution's gradient. Consequently, the operational feasibility of locating the global optimum of function $f$ via the aforementioned approach faces limitations. However, Remark \ref{remark-pde} shows that for certain specific functions $f$, the sign of the gradient of the PDE solution coincides with that of the gradient of $f$ itself. This motivates us to investigate whether we can directly leverage the sign of function $f$'s gradient to construct a class of asymptotically optimal strategies. Such strategies would locate the optimum of $f$ within the two-armed decision-making framework. To address this, we present a simple but novel algorithm, which we call ``Strategic Monte Carlo Optimization'' in the next subsection.

\subsection{Strategic Monte Carlo Optimization (SMCO) Algorithm}\label{sec-smco}
 
Motivated by the theoretical results in previous subsections, we now introduce a simple strategy  $\bmtheta^*$ 
that just uses the sign information of the coordinate-wise partial derivative of $f$.

At round $m\ge 1$, we construct $\bmvartheta_m^*=(\vartheta_{m,1}^*,\cdots,\vartheta_{m,d}^*)$ as follows:\footnote{When $m=1$, $S_0^{\bmtheta^*}/0$ can be defined as $\eta$, which is uniform distributed on $\Gamma$.}
\begin{equation}\label{strategy-op}
\vartheta_{m,j}^*=\left\{\begin{array}{rl}
1,\quad &\text{if } \partial_{x_j}f\left(\frac{S_{m-1}^{\bmtheta^*}}{m-1}\right)\ge0;
\\
0,\quad & \text{otherwise},
\end{array}
\right.\quad j=1,\cdots,d.
\end{equation}

Intuitively, suppose the current sample mean is in a region where the function $f$ is increasing, then we aim to maximize $f$ by selecting a sample from the arm with a higher population mean in the next round; if the sample mean is in a region where the function $f$ is decreasing, then we select the next sample from the arm with a lower (population) mean. Figure \ref{cmc_smc}(c) shows the different iterative mechanisms for attaining the maximum value by the GD method (blue arrow line) and our SMCO (red arrow line) respectively. It also demonstrates that under the simple sign-based strategy \(\bmtheta^*\) given in \eqref{strategy-op}, the sample average \(S_n^{\bmtheta^*}/n\) converges to the one of the local maxima of \(f\). This result constitutes the Strategic Law of Large Numbers for the decision model problem (Theorem \ref{thm-lln} below).
In the one-dimensional setting, our proposed two-armed decision model process devises strategies that combine the two arms with known distributions to identify the maximum value of a non-concave function \(f\). This approach differs from the classical two-armed decision model, which explores optimal arms to achieve the maximum average reward for each arm under a monotonic function \(f\) with unknown distributions.

\begin{algorithm}
\caption{\label{alg:smco_algo}Strategic Monte Carlo Optimization (SMCO)}
{\bf Input: } \text{Function } $f$,   \text{ tolerance level } $\varepsilon$, initial point $\hat{x}_0$, \text{ with} $S_0^{\bmtheta^*} := \hat{x}_0$.

\hspace{0.2cm}\textbf{for }$n = 0$ \textbf{ to } $\infty$ \textbf{ do} 

\hspace{0.8cm}\textbf{for }$j = 1$ \textbf{ to } $d$ \textbf{ do} 

  \hspace{1.4cm}\textbf{if }$\partial_{x_j}f\left(\hat{x}_n\right) \ge 0$ \textbf{ then }
    
     \hspace{1.8cm}$\vartheta^*_{n+1,j}=1$, draw $Z_{n+1,j}^{\bmtheta^*}=X_{n+1,j}$ from arm $X_j$ (larger mean in [\ref{denote-xy}]);
   
   \hspace{1.4cm}\textbf{else } 
   
      \hspace{1.8cm} $\vartheta^*_{n+1,j}=0$, draw $Z_{n+1,j}^{\bmtheta^*}=Y_{n+1,j}$ from arm $Y_j$ (smaller mean in [\ref{denote-xy}]).
    
    \hspace{1.4cm}\textbf{end }    

\hspace{0.8cm}\textbf{end }  

 \hspace{0.8cm}Update  $S_{n+1}^{\bmtheta^*} :=S_n^{\bmtheta^*}+\Z_{n+1}^{\bmtheta^*}$, with $\Z_{n+1}^{\bmtheta^*}=(Z_{n+1,1}^{\bmtheta^*},\cdots,Z_{n+1,d}^{\bmtheta^*})$. 
    
 \hspace{0.8cm}Update $\hat{x}_{n+1} = S_{n+1}^{\theta^*} / (n+1).$

  \hspace{0.8cm}\textbf{if }$\left|f\left(\hat{x}_{n+1}\right)-f\left(\hat{x}_n\right)\right|\le \varepsilon$
\textbf{ then }

\hspace{1.2cm}\textbf{break }

 \hspace{0.8cm}\textbf{end }

\hspace{0.2cm}\textbf{end }\smallskip
    
{\bf Output: } Final point $\hat{x}^* := \hat{x}_{n+1}$ and final function value $\hat{f}^* := f\left(\hat{x}^*\right)$.\smallskip
\end{algorithm}




    

 The following theorem establishes the StLLN according to the sign-based strategy $\bmtheta^*$ described in Algorithm \ref{alg:smco_algo}, which ensures that under $\bmtheta^*$, the sample mean $S_n^{\bmtheta^*}/n$ converges to one of the local maxima of $f$ on $\Gamma$. 
\color{black}

\begin{Thm}[Local Optimization under the Sign-Based Strategy]\label{thm-lln} 

Let $f$ be a continuously differentiable function on $\Gamma_{\delta}$, and satisfy the following condition: For any $k=1,\cdots,d$ and a fixed integer $K_k\in\mathbb{N}$,  there exist two sets of points $\underline{\mu}_k\leq c_k^{1}<\cdots<c_{k}^{K_k}\leq\overline{\mu}_k$ and $\underline{\mu}_k-\delta=b_{k}^{0}<b_{k}^{1}<\cdots<b_{k}^{K_k-1}<b_{k}^{K_k}=\overline\mu_k+\delta$, 
such that for any fixed $x=(x_1,\cdots,x_d)\in\Gamma_{\delta}$, 
the function $f(x_1,..,x_{k-1},\cdot,x_{k+1},..,x_d)$  is  strictly increasing on $(b_{k}^{j-1},c_k^{j})$ and strictly decreasing on $[c_k^{j},b_{k}^{j})$ for $j=1,\cdots,K_{k}$.
 \smallskip 
 
 Let $\mathcal{M}$ denote all local maximum points of $f$ on $\Gamma$, and for any $\varepsilon>0$,  $\mathcal{M}_{\varepsilon}=\bigcup_{c\in\mathcal{M}}B(c,\epsilon)$, where $B(c,\varepsilon)=\{x\in \Gamma_{\delta}:\|x-c\|_{\infty}<\varepsilon\}$. Let   $\bmtheta^*=(\bmvartheta^*_1,\cdots,\bmvartheta_m^*,\cdots)$, with $\bmvartheta_m^*=(\vartheta^*_{m,1},\cdots,\vartheta_{m,d}^*)$, be the strategy defined in \eqref{strategy-op}. Then:
\begin{equation}\label{lln-inde}
\lim_{n\to\infty}P\left(\frac{S_n^{\bmtheta^*}}{n}\in\mathcal{M}_{\varepsilon}\right)=1.
\end{equation}
\end{Thm}
 \begin{remark}
     In Theorem \ref{thm-lln}, we make a sufficient condition that $f$  is continuously differentiable. In fact, the above conclusion can be obtained similarly for continuous function $f$ that is not differentiable. For a non-differentiable function, we can approximate it by a smooth function, for example, $f_{\varepsilon}(x)=\int f(x+\varepsilon y)\phi(y)dy\to f(x),$ as $\varepsilon\to0$, where $\phi$ is the probability density function of a standard normal distribution. Then, one can simulate the extremes of $f$ by applying the method to $f_{\varepsilon}$ for a small enough $\varepsilon$. 
 \end{remark}
 
\begin{remark}
The coordinate separability condition essentially ensures that (local) coordinate-wise updates are sufficient to increase the value of the function in each iteration. The condition plays a role that is akin to the coordinate-wise smoothness condition (Lipcchitz property of gradients) imposed for sign gradient descent \cite{Bern2018}  and coordinate gradient descent \cite{Richtarik2016}. Note that we do not impose additional smoothness conditions beyond the stated monotonicity requirements.
\end{remark}

Next, we show that SMCO, when augmented with a growing set of multiple starting points, attains the global maximum for functions with multiple local maxima. Specifically, we show theoretically that, given a starting point in a local neighborhood of the global maximizer, the SMCO algorithm modified with a ``boosted'' starting iteration index is guaranteed to converge to the global maximizer. 
\begin{Thm}\label{thm-ms}
    Let $f$ be a continuously differentiable function on $\Gamma_{\delta}$. For any global maximum $x^*$ of $f$ on $\Gamma$, assume that there exists $\delta_0>0$, such that for any $y\in B(x^*,\delta_0)$, $sgn(\partial_{x_i}f(y))=sgn(x_i^*-y_i),i=1,\cdots,d.$ Then there exists an integer $N_0>0$ such that for any $m\geq N_0$, whenever $S_m^{\bmtheta}/m\in B(x^*, \delta_0)$, it holds that for any $\varepsilon>0$,
    \begin{equation}
\lim_{n\to\infty}P\left(\left\|\frac{S_n^{\bmtheta'}}{n}-x^*\right\|_{\infty}<\varepsilon\right)=1, 
    \end{equation}
    where $B(x^*, \delta_0)=\{x\in \Gamma_{\delta}:\|x-x^*\|_{\infty}<\delta_0\}$ and $\bmtheta'=(\bmvartheta_1,\cdots,\bmvartheta_m,\bmvartheta_{m+1}^*, \bmvartheta_{m+2}^* \ldots)$. 
    
    Consequently, given a growing (deterministically, randomly, or quasi-randomly generated) set of starting points ${\cal X}_m$ such that $\#({\cal X}_m) = m$ and $P(B(x^*,\delta_0) \cap {\cal X}_m \neq \emptyset) \to 1 $ as $m\to\infty$, the SMCO algorithm run from all points in ${\cal X}_m$ with a ``boosted'' starting iteration index\footnote{Formally, running SMCO with a boosted iteration index $n_0$ means the following: given a starting point $x_0$, set the initial sum $S_0 := n_0 x_0 $, and in the $n$-th iteration, update the running sum as in Algorithm \ref{alg:smco_algo} but compute the sample mean $x_{(n+1)} := S_{(n+1)}/(n+n_0)$. Effectively, this is as if the algorithm starts with the iteration counter $n = n_0$ instead of $n=0$.} $n_0 = N_0$ is guaranteed to attain a global maximizer $x^*$ asymptotically as $n\to\infty$ and $m\to\infty$.
\end{Thm}

Finally, we present a rate of convergence of our SMCO algorithm 1 when $f$ has a unique local maximum point $x^*\in \Gamma$, which could be on the boundary of $\Gamma$. 
\begin{Thm}[\textbf{Convergence Rate II}]\label{thm-convergence2}
   Let $f$ be a continuous twice-differential function on $\Gamma_{\delta}$, and assume that $f$ has an unique local maximum point $x^*=(x_1^{*},\cdots,x_{d}^{*})$ on $\Gamma$, and further the $i{\rm th}$ component of $f$ is strictly increasing on $(\underline\mu_i-\delta,x_i^{*})$ and strictly decreasing on $[x_i^{*},\overline{\mu}_i+\delta)$. Then we have, 
   \begin{align}
\left|E\left[f\left(\frac{S_n^{\bmtheta^*}}{n}\right) \right]-f(x^*)\right|
\leq  \max\left\{L_f\left(\frac{\sigma\sqrt{d}}{\sqrt{n}}+\frac{\sqrt{E[\|\eta\|^2]}}{n}\right),\frac{K_f\left(d\sigma^2+\sum_{j=1}^d(\overline{\mu}_j-\underline{\mu}_j)^2\right)+L_fE[\|\eta\|]}{n}\right\},\label{0convergence-rate-1}
   \end{align}
 where $L_{f}$ is the Lipschitz constant of $f$, $K_f=\max_{i,j}\sup_{x\in\Gamma_{\delta}}|\partial_{x_ix_j}^2f(x)|$ and $\sigma^2$ is the uniform variances of $\{X_j\}$ and $\{Y_j\}$ as defined in \eqref{denote-xy}. Moreover,
 \begin{align}
    E\left[\left\|\frac{S_n^{\bmtheta^*}}{n}-x^*\right\|^2\right]\leq\frac{d\sigma^2+\sum_{j=1}^d(\overline{\mu}_j-\underline{\mu}_j)^2}{n}+\frac{E[\|\eta\|^2]}{n^2}.
\end{align}
\end{Thm}

\section{\label{sec:Numerical}Numerical Experiments}

\subsection{Three Implemented Versions of SMCO algorithms}
We implement the following three versions of our SMCO algorithms, which we use in the numerical experiments in Section \ref{sec:Numerical}. See Appendix \ref{append-B} for some additional details .

\medskip

\noindent\textbf{(1) SMCO:} 
Our baseline algorithm that implements Algorithm \ref{alg:smco_algo} with the adaptive-size finite difference as described above. We set the random variable $\xi_j$ in Equation [\ref{denote-xy}] to be $U[-\delta_j,\delta_j]$ with $\delta_j := \delta \left(\overline{\mu}_j - \underline{\mu}_j \right).$ We set the bound expansion parameter
$\delta=0.05$, the convergence tolerance at 1e-6, and the
maximum number of iterations to be 500. 

\medskip

\noindent\textbf{(2) SMCO-R:} A refined version of SMCO with two programming/meta-heuristic
enhancements. First, SMCO-R records and exploits the running maximum of
all test function evaluations: This is purely bookkeeping and does not add any additional evaluations of the test functions beyond those for finite-difference calculations and convergence checks. From a programming perspective, the running maximum generated in past iterations is ``free information'' that requires essentially no additional computation and memory storage.  Second, SMCO-R adopts a two-stage
procedure: it starts by running SMCO for at most 50\% iterations of
the designated maximum number of iterations, and then runs the rest 50\%
iterations with a more local search (and smaller step sizes). Specifically,
the end point of the initial stage is used as the start point of the
second stage, namely $x_{0}^{\left(2\right)}$, and the second stage
SMCO with the starting iteration counter set to $n=1000$. Specifically,
the partial sum is initialized at $S_{0}^{\left(2\right)}:=1000x_{0}^{\left(2\right)}$,
the next point is updated according to $x_{1}^{\left(2\right)}:=\left(S_{0}^{\left(2\right)}+\Z_{1}^{\left(2\right)}\right)/1001$, and
etc. Other hyperparameters are set to be the same as in SMCO.

\medskip
\noindent\textbf{(3) SMCO-BR}: A ``boosted'' version of SMCO-R that runs two passes
of SMCO-R with: (i) the first pass is exactly the same as SMCO-R,
(ii) the second pass runs a more local version of SMCO-R with the
starting iteration counter set to $n=100$, as motivated by Section \ref{sec-smco} and Theorem \ref{thm-ms}. Specifically, given
any starting point $x_{0}$, the second pass initializes the running sum as $S_{0}:=100x_{0}$, computes the first update as $x_{1}=\left(S_{0}+\Z_{1}\right)/101$,
and so on and so forth. Since SMCO-BR effectively runs two separate
passes of SMCO-R, to make the comparison fair, we set the maximum
number of iterations to be half of that for SMCO/SMCO-R, so that
the running times of the three versions are similar. Other hyperparameters are set to be the same as in
SMCO.

\begin{remark}
As we mentioned above, the implementation of our SMCO algorithms does not use signs of the true gradient of the objective function $f$, but just signs of coordinate-wise finite difference of $f$ for updating. This is why the SMCO algorithms perform robustly well even for discontinuous $f$ in the extensive numerical studies. In addition, we believe the SMCO algorithms can be extended to allow for non-rectangular or constrained domains via penalty methods. In particular, we consider the following general constrained problem:
\begin{align*}
   \max_{x \in \Gamma} f(x)
    \quad \text{subject to } 
\begin{cases}
g_i(x) = 0 & \text{for } i = 1, 2, \dots, m, \\
h_j(x) \leq 0 & \text{for } j = 1, 2, \dots, k,
\end{cases}
\end{align*}
where $f: \Gamma \to \mathbb{R}$, $g_i: \Gamma \to \mathbb{R}$ and $h_j: \Gamma \to \mathbb{R}$ are all continuous functions. The penalty method converts this to an unconstrained problem by embedding constraints into a modified objective function:
\[
F(x) = f(x) - \lambda_1 \sum_{i=1}^m \left[ g_i(x) \right]^2 - \lambda_2 \sum_{j=1}^k \left[ \max\left\{ 0, h_j(x) \right\} \right]^2~, 
\]
where $\lambda_1, \lambda_2 > 0$ are penalty parameters. Our SMCO then solve the transformed global optimization problem $\max_{\Gamma} F(x)$ on the original rectangular domain $\Gamma$, with penalties ensuring convergence to feasible points. We leave it to future work to study the theoretical properties of the SMCO algorithms applied to the modified optimization problem $\max_{x \in \Gamma} F(x)$.
\end{remark}

\subsection{Alternative Optimization Algorithms in R for Comparison}
Below we list two groups of state-of-the-art optimization algorithms (implemented in R) to which we
compare our SMCO/-R/-BR algorithms in terms of performance accuracy (in root mean squared error) and computation speed (in second). See Appendix \ref{append-C} for more information on these algorithms and their hyperparameter configurations.

\textbf{Algorithm Group I: Six Local Optimizers} This set consists of six gradient-based or derivative-free optimization algorithms that are designed to converge to a local optimum starting from a given starting point: (1) GD, (2) L-BFGS, (3) SignGD, (4) SPSA, (5) ADAM, (6) BOBYQA.  To make the comparison as fair as possible, we pass the exact sets
of starting points used for our SMCO algorithms to the following
algorithms as well: for single-start comparisons, our SMCO algorithms and the Group I algorithms are all run from a common randomly drawn starting point; for multi-start comparisons, they are all run from the same set of starting points.

\textbf{Algorithm Group II: Six Global Optimizers} This set consists of six R packages that implement meta-heuristic global optimization algorithms: (1) GenSA (2) SA (3) DEoptim (4) STOGO (5) GA (6) PSO. We use the default hyperparameters set in these packages. Since these algorithms are explicitly designed for global optimization (and many already implement multiple starting points from within), we do not pass any generated starting points to these algorithms whenever applicable.\footnote{A few of these global algorithms nevertheless require a user-designated starting point, in which case we pass a start point randomly chosen from the set of multiple starting points generated for our SMCO and the Group I algorithms.}


\subsection{Optimization of Deterministic Test Functions}

\begin{figure}
\caption{\label{fig:testf_plot}Visualization of Test Functions}

\begin{centering}
\subfloat[Rastrigin Function]{\includegraphics[scale=0.22]{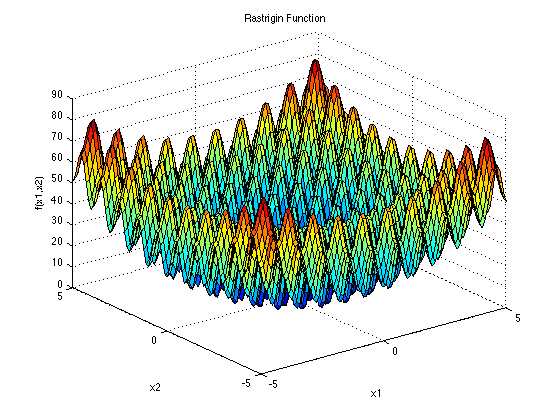}}
\subfloat[Ackley Function]{\includegraphics[scale=0.22]{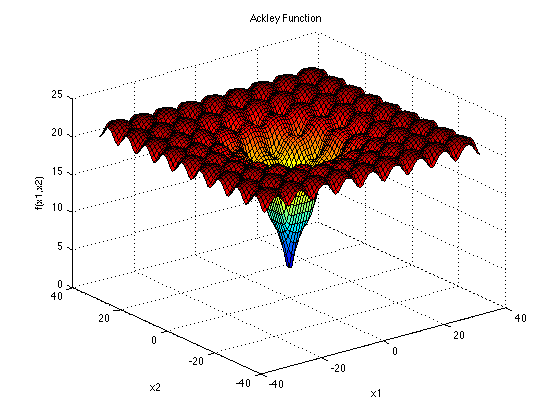}}
\subfloat[Michalewicz Function]{
\includegraphics[scale=0.22]{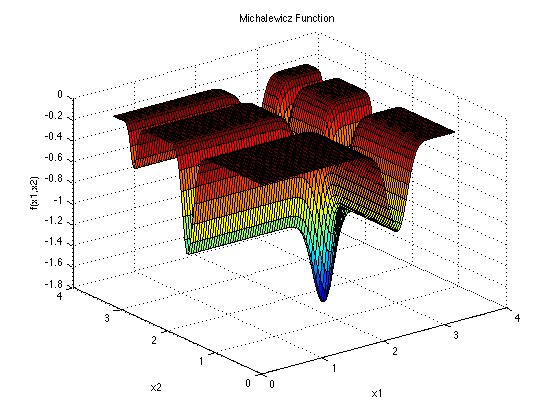}}
\subfloat[Griewank Function]{\includegraphics[scale=0.34]{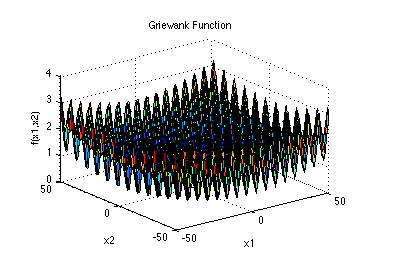}}
\end{centering}
\par\begin{centering}
Retrieved from: https://www.sfu.ca/\textasciitilde ssurjano/optimization.html
\par\end{centering}
\end{figure}

We first run the optimization algorithms on the following set of deterministic
test functions, which have been widely used to benchmark the performances
of global optimization algorithms. The function definitions (and hyperparameter
values, if applicable) follow those given in the optimization section
of the Virtual Library of Simulated Experiments: Test Functions and
Datasets\footnote{Available online at https://www.sfu.ca/\textasciitilde ssurjano/optimization.html}
\cite{Surja2013}. We have tested a wide array of test functions
from different categories with various dimensionalities. Due to page
limit, we only report results in the main text on the following four test
functions:
(1) \textbf{Rastrigin} function on $\left[-5.12,5.12\right]^{d}$. This function
has a unique global minimum of value $0$ at the origin and multiple
global maxima at $x\in\left\{ \pm4.52299\right\} ^{d}$. 
(2) \textbf{Ackley} function on $\left[-32.768,32.768\right]^{d}$. This function
has a unique global minimum of value $0$ at the origin and multiple
global maxima of unknown values.
(3) \textbf{Griewank} function on $\left[-600,600\right]^{d}$. This function has
a unique global minimum of value $0$ at the origin and multiple
global maxima of unknown values.
(4) \textbf{Michalewicz} function on $\left[0,\pi\right]^{d}$. This function has
a unique global minimum (of values reported up to $d=75$) and
multiple global maxima of the value $0$ at the vertices $\left\{ 0,\pi\right\} ^{d}$. These functions are constructed to have challenging optimization landscapes
with many local optima, boundary optima, (near) flat regions/directions,
steep drops/ridges, or a combination of these. Figure \ref{fig:testf_plot}
contains visualizations (retrieved from the Virtual Library of Simulated
Experiments) of these test functions in 2-dimensional forms.

One notable feature of several of the test functions above lies in their symmetric forms, i.e., they are constructed with symmetric domains and its global minimum attained at the origin (or, the center point of the domain). Such a kind of symmetry is unlikely to be found in realistic optimization problems. This would make the optimization problem drastically easier for any algorithm that explicitly or implicitly searches the center point of the domain, making the comparison between algorithms less clear and fair. Hence, we randomly rotate, shift, and asymmetrize the default domains of the test functions above, which we formally describe in the following. Specifically, for a given
test function $f$ on its default domain $\prod_{j=1}^{d}[\underline{\mu}_{j},\overline{\mu}_{j}]$,
write $\Delta\mu_{j}:=\overline{\mu}_{j}-\underline{\mu}_{j}$ be
the length of the $j$-th coordinate domain. Then, we shift and asymmetrize the
domain in each coordinate independently as follows: for each $j=1,...,d$,
we draw $r_{j}\sim Bernoulli\left(0.5\right),\ v_{j}\sim N(0,1)$, and $\nu_{j}\sim U[0,1]$,
and set \begin{align*}
    \underline{\mu}_{j}^{mod}:= &\ \underline{\mu}_{j}+\left[v_{j}+r_{j}\left(0.2+0.1\nu_{j}\right)-\left(1-r_{j}\right)\left(0.4+0.2\nu_{j}\right)\right]\Delta\mu_{j}, \\ 
\overline{\mu}_{j}^{mod}:=& \ \overline{\mu}_{j}+\left[v_{j}+r_{j}\left(0.4+0.2\nu_{j}\right)-\left(1-r_{j}\right)\left(0.2+0.1\nu_{j}\right)\right]\Delta\mu_{j}.
\end{align*}
Essentially, the above functions constitute the following operations in each
coordinate. (i) Randomly shift the center point of the original domain by $v_{j}\cdot\Delta\mu_{j}$; (ii) Randomly decide whether to asymmetrize the domain rightwards $\left(r_{j}=1\right)$
or leftwards $\left(r_{j}=0\right)$. 
If asymmetrizing rightwards, i.e., $r_{j}=1$, then push the
lower bound $\underline{\mu}_{j}$ rightwards by $25\%\pm5\%$ of
$\Delta\mu_{j}$, and push the upper bound $\overline{\mu}_{j}$ rightwards
by $50\%\pm10\%$ of $\Delta\mu_{j}$.  
If asymmetrizing leftwards, i.e., $r_{j}=0$, then push the lower
bound $\underline{\mu}_{j}$ leftwards by $50\%\pm10\%$ of $\Delta\mu_{j}$,
and push the upper bound $\overline{\mu}_{j}$ leftwards by $25\%\pm5\%$
of $\Delta\mu_{j}$. Note that the domain of each coordinate is perturbed independently,
so the ``asymmetrization directions'' generally differ across coordinates. After shifting and asymmetrizing the domain into 
${\cal X}^{mod}:=\prod_{j=1}^{d}\left[\underline{\mu}_{j}^{mod},\overline{\mu}_{j}^{mod}\right]$ as described above, we further randomly rotate the domain by randomly
drawing a rotation (orthonormal) matrix $Q$ of dimension $d$, and define the final test function by $$
f^{mod}\left(x\right)\ :=\ f\left(Q\left(x-v\odot\Delta\mu\right)\right),$$
where $\odot$ represents the coordinate-wise (Hadamard) product.

For each of the test functions (with rotation and domain modification), we test the optimization algorithms
both for minimization and maximization. Due to the randomness in the random starting point generation
as well as the stochastic nature of many optimization algorithms used, we run multiple replications of all optimization algorithms and report summary statistics for their performances across the replications. Specifically, we report the following performance measures: ``Time'', ``RMSE'', and ``AExx'', whose definitions we now explain. ``Time'' reports the average time (in second) it takes for a particular algorithm to return the final function value, averaged across all the replications. ``RMSE'' and ``AExx'' are relative accuracy measures defined as follows. Since the true maximum and minimum values may be different on the rotated, enlarged, and asymmetrized domain from those on the default domains, we first define ``Best Value'' as the best value found by all algorithms (ours and those in Groups I \& II) over all replications. We then define ``RMSE'' as the square root of the mean squared error between the value returned by a specific algorithm and the ``Best Value'' found, and define ``AExx'' as the xx-th percentile of the absolute error between the value returned by a specific algorithm and the ``Best Value''.

\textbf{Table \ref{tab:Comp_SS_d2}} summarizes how our SMCO algorithms  \textbf{Group I (local)} algorithms in maximizing and minimizing all four test functions of \textbf{dimension} $\bm{d=2}$, when run from a \textbf{single} uniformly generated starting point, based on $500$ replications. 
First,  our SMCO algorithms tend to produce significantly more accurate results than all Group I algorithms in most of the test configurations reported in Table \ref{tab:Comp_SS_d2}. The RMSE and AE percentiles of our SMCO algorithms are often orders of magnitude smaller than those of other algorithms in most configurations. Second, only our SMCO algorithms display robust performance across all eight test configurations in comparison with all other algorithm. Finally, our SMCO algorithms display remarkably stable performance in terms of the error percentiles: the differences between the 50th, 95th, and 99th error percentiles of our SMCO algorithms tend to be substantially smaller than those of all other Group I algorithms, demonstrating that the performance of our SMCO algorithms is much \emph{less sensitive} to the random starting point relative to others.

\begin{sidewaystable}
\centering{}\caption{\label{tab:Comp_SS_d2}Deterministic Functions, Dimension $d=2$: Comparison with Group I Algorithms}
\begin{threeparttable}
\begin{tabular}{|c|cccccc|ccccc|}
\hline 
\multicolumn{2}{|c|}{Performance} & RMSE & AE50 & AE95 & AE99 & Time & RMSE & AE50 & AE95 & AE99 & Time\tabularnewline
\hline 
\multirow{20}{*}{Min} & \textbf{Function: Best Value} & \multicolumn{5}{c|}{\textbf{Rastrigin: 0}} & \multicolumn{5}{c|}{\textbf{Griewank: 0}}\tabularnewline
\cline{3-12} \cline{4-12} \cline{5-12} \cline{6-12} \cline{7-12} \cline{8-12} \cline{9-12} \cline{10-12} \cline{11-12} \cline{12-12} 
 & \textbf{SMCO} & \textbf{7.8} & \textbf{8.0} & \textbf{13.0} & \textbf{13.0} & \textbf{0.016} & \textbf{9.9} & \textbf{9.8} & \textbf{11.1} & \textbf{11.6} & \textbf{0.017}\tabularnewline
 & \textbf{SMCO-R} & \textbf{2.8} & \textbf{2.0} & \textbf{5.0} & \textbf{8.0} & \textbf{0.016} & \textbf{0.04} & \textbf{0.03} & \textbf{0.08} & \textbf{0.12} & \textbf{0.017}\tabularnewline
 & \textbf{SMCO-BR} & \textbf{2.7} & \textbf{2.0} & \textbf{5.0} & \textbf{8.0} & \textbf{0.016} & \textbf{0.06} & \textbf{0.03} & \textbf{0.11} & \textbf{0.16} & \textbf{0.018}\tabularnewline
 \cline{3-12} \cline{4-12} \cline{5-12} \cline{6-12} \cline{7-12} \cline{8-12} \cline{9-12} \cline{10-12} \cline{11-12} \cline{12-12} 
 & GD & 36.7 & 17.0 & 80.9 & 104.7 & 0.010 & 196.2 & 141.4 & 345.6 & 437.4 & 0.032\tabularnewline
 & SignGD & 62.3 & 42.2 & 115.4 & 156.2 & 0.035 & 208.4 & 141.4 & 403.3 & 492.5 & 0.026\tabularnewline
 & SPSA & 39.6 & 16.9 & 84.6 & 135.4 & 0.002 & 137.9 & 72.6 & 277.6 & 330.7 & 0.003\tabularnewline
 & ADAM & 15.4 & 8.3 & 30.0 & 48.5 & 0.008 & 211.5 & 142.1 & 418.3 & 506.2 & 0.014\tabularnewline
 & L-BFGS & 43.8 & 15.9 & 97.5 & 115.5 & 0.0004 & 141.9 & 78.8 & 289.6 & 333.8 & 0.0005\tabularnewline
 & BOBYQA & 42.9 & 9.9 & 88.5 & 135.4 & 0.002 & 18.9 & 0.2 & 0.8 & 15.9 & 0.002\tabularnewline
\cline{2-12} \cline{4-12} \cline{5-12} \cline{6-12} \cline{7-12} \cline{8-12} \cline{9-12} \cline{10-12} \cline{11-12} \cline{12-12} 

 & \textbf{Function: Best Value} & \multicolumn{5}{c|}{\textbf{Ackley: 4.44e-16}} & \multicolumn{5}{c|}{\textbf{Michalewicz: -1.801303}}\tabularnewline
\cline{3-12} \cline{4-12} \cline{5-12} \cline{6-12} \cline{7-12} \cline{8-12} \cline{9-12} \cline{10-12} \cline{11-12} \cline{12-12} 
 & \textbf{SMCO} & \textbf{3.1} & \textbf{0.4} & \textbf{0.9} & \textbf{19.8} & \textbf{0.018} & \textbf{1.25} & \textbf{1.15} & \textbf{1.64} & \textbf{1.64} & \textbf{0.018}\tabularnewline
 & \textbf{SMCO-R} & \textbf{2.3} & \textbf{0.01} & \textbf{0.01} & \textbf{19.5} & \textbf{0.019} & \textbf{1.11} & \textbf{1.14} & \textbf{1.37} & \textbf{1.43} & \textbf{0.018}\tabularnewline
 & \textbf{SMCO-BR} & \textbf{2.6} & \textbf{0.01} & \textbf{0.01} & \textbf{19.4} & \textbf{0.018} & \textbf{1.09} & \textbf{1.14} & \textbf{1.15} & \textbf{1.37} & \textbf{0.018}\tabularnewline
 \cline{3-12} \cline{4-12} \cline{5-12} \cline{6-12} \cline{7-12} \cline{8-12} \cline{9-12} \cline{10-12} \cline{11-12} \cline{12-12}
 & GD & 19.0 & 19.9 & 20.0 & 20.3 & 0.021 & 1.74 & 1.84 & 1.99 & 1.99 & 0.012\tabularnewline
 & SignGD & 19.0 & 19.9 & 20.0 & 20.0 & 0.034 & 1.55 & 1.63 & 1.99 & 1.99 & 0.023\tabularnewline
 & SPSA & 19.0 & 19.9 & 20.0 & 20.0 & 0.002 & 1.45 & 1.34 & 1.99 & 1.99 & 0.004\tabularnewline
 & ADAM & 19.1 & 19.9 & 20.0 & 20.7 & 0.001 & 1.84 & 1.99 & 2.00 & 2.00 & 0.006\tabularnewline
 & L-BFGS & 18.8 & 19.9 & 20.0 & 20.4 & 0.0005 & 1.65 & 1.63 & 1.99 & 1.99 & 0.001\tabularnewline
 & BOBYQA & 15.1 & 18.0 & 20.0 & 20.2 & 0.002 & 1.47 & 1.37 & 1.99 & 1.99 & 0.003\tabularnewline

\hline 
\multirow{18}{*}{Max} & \multicolumn{1}{c}{\textbf{Function: Best Value}} & \multicolumn{5}{c|}{\textbf{Rastrigin: 213.5824}} & \multicolumn{5}{c|}{\textbf{Griewank: 622.3434}}\tabularnewline
\cline{3-12} \cline{4-12} \cline{5-12} \cline{6-12} \cline{7-12} \cline{8-12} \cline{9-12} \cline{10-12} \cline{11-12} \cline{12-12} 
 & \textbf{SMCO} & \textbf{17.4} & \textbf{10.3} & \textbf{10.4} & \textbf{88.7} & \textbf{0.015} & \textbf{2.5} & \textbf{1.0} & \textbf{5.5} & \textbf{6.7} & \textbf{0.016}\tabularnewline
 & \textbf{SMCO-R} & \textbf{13.4} & \textbf{10.2} & \textbf{10.2} & \textbf{58.5} & \textbf{0.015} & \textbf{1.1} & \textbf{3e-13} & \textbf{2.1} & \textbf{3.6} & \textbf{0.015}\tabularnewline
 & \textbf{SMCO-BR} & \textbf{15.6} & \textbf{10.2} & \textbf{10.2} & \textbf{80.6} & \textbf{0.015} & \textbf{2.2} & \textbf{0.1} & \textbf{4.8} & \textbf{8.5} & \textbf{0.016}\tabularnewline
\cline{3-12} \cline{3-12} \cline{4-12} \cline{5-12} \cline{6-12} \cline{7-12} \cline{8-12} \cline{9-12} \cline{10-12} \cline{11-12} \cline{12-12} 
 & GD & 85.3 & 82.7 & 173.3 & 192.9 & 0.004 & 469.6 & 477.0 & 609.1 & 618.5 & 0.033\tabularnewline
 & SignGD & 130.2 & 130.9 & 171.1 & 173.1 & 0.035 & 471.1 & 477.0 & 609.1 & 618.5 & 0.027\tabularnewline
 & SPSA & 85.1 & 80.6 & 167.2 & 176.0 & 0.004 & 472.0 & 477.6 & 609.1 & 618.9 & 0.013\tabularnewline
 & L-BFGS & 82.8 & 81.9 & 165.0 & 171.1 & 0.0004 & 447.6 & 459.9 & 608.7 & 618.5 & 0.0005\tabularnewline
 & BOBYQA & 73.6 & 20.3 & 167.0 & 171.1 & 0.002 & 204.4 & 0.0 & 293.4 & 575.1 & 0.002\tabularnewline
\cline{2-12} \cline{4-12} \cline{5-12} \cline{6-12} \cline{7-12} \cline{8-12} \cline{9-12} \cline{10-12} \cline{11-12} \cline{12-12} 

 & \textbf{Function: Best Value} & \multicolumn{5}{c|}{\textbf{Ackley: 22.35029}} & \multicolumn{5}{c|}{\textbf{Michalewicz: 1.21406}}\tabularnewline
\cline{3-12} \cline{4-12} \cline{5-12} \cline{6-12} \cline{7-12} \cline{8-12} \cline{9-12} \cline{10-12} \cline{11-12} \cline{12-12} 
 & \textbf{SMCO} & \textbf{0.20} & \textbf{0.18} & \textbf{0.29} & \textbf{0.37} & \textbf{0.017} & \textbf{0.25} & \textbf{0.21} & \textbf{0.41} & \textbf{0.41} & \textbf{0.017}\tabularnewline
 & \textbf{SMCO-R} & \textbf{0.03} & \textbf{0.02} & \textbf{0.06} & \textbf{0.08} & \textbf{0.018} & \textbf{0.17} & \textbf{6.5e-7} & \textbf{0.41} & \textbf{0.41} & \textbf{0.016}\tabularnewline
 & \textbf{SMCO-BR} & \textbf{0.02} & \textbf{0.01} & \textbf{0.05} & \textbf{0.08} & \textbf{0.018} & \textbf{0.15} & \textbf{7.6e-7} & \textbf{0.41} & \textbf{0.41} & \textbf{0.016}\tabularnewline
 \cline{3-12} \cline{4-12} \cline{5-12} \cline{6-12} \cline{7-12} \cline{8-12} \cline{9-12} \cline{10-12} \cline{11-12} \cline{12-12} 
 & GD & 1.90 & 0.08 & 4.44 & 8.83 & 0.017 & 0.85 & 0.72 & 1.21 & 1.21 & 0.011\tabularnewline
 & SignGD & 2.20 & 0.08 & 4.75 & 9.55 & 0.032 & 0.65 & 0.41 & 1.21 & 1.21 & 0.024\tabularnewline
 & SPSA & 2.19 & 0.08 & 4.75 & 9.55 & 0.002 & 0.94 & 1.21 & 1.21 & 1.21 & 0.002\tabularnewline
 & L-BFGS & 1.66 & 0.08 & 3.64 & 7.71 & 0.000 & 0.77 & 0.41 & 1.21 & 1.21 & 0.001\tabularnewline
 & BOBYQA & 0.25 & 0.02 & 0.61 & 1.19 & 0.002 & 0.58 & 0.41 & 1.21 & 1.21 & 0.002\tabularnewline
\hline 
\end{tabular} 
\begin{tablenotes}
\item Based on 500 reps. In each rep, all algorithms run from a \textbf{single uniformly drawn starting point}. Computed using Posit Cloud. SMCO/-R/-BR run with 1e-6 tolerance and 500 maximum iterations. 
\end{tablenotes}
\end{threeparttable}
\end{sidewaystable}

In Appendix \ref{append-D}, \textbf{Tables \ref{tab:Comp2d}-\ref{tab:Comp50d}} report comparisons of our SMCO algorithms with both Group I and II algorithms for \textbf{dimensions} $\bm{d = 2}$, $\bm{10}$,  $\bm{20}$,   $\bm{50}$, in which our SMCO algorithms and all Group I algorithms are run with $10\sqrt{d} \approx 14$, $32$, $45$, and $71$ \textbf{uniformly drawn starting points}. The general pattern there shows that: First, our SMCO algorithms perform robustly well across most settings, both in terms of RMSE and AE99. Second, Group I (local) algorithms perform significantly worse than SMCO and Group II (global) algorithms, especially when the dimension increases. Consequently, due to the space limit, we only compare our SMCO algorithms with Group II algorithms in higher-dimensional settings in the main text.

We now report results on how our SMCO algorithms compare with Group II (global) algorithms in a \textbf{high-dimensional settings} ($d=200$). Here our SMCO  algorithms are run from $\sqrt{d} \approx 14$  uniformly drawn starting points in a parallelized manner. The Group II (global) algorithms are run without any specified starting point (when the specification of a starting point is not required by the algorithm) or one randomly drawn starting point (when the specification of a starting point is required by the algorithm), again at their default settings. Again, our (multi-start) SMCO algorithms demonstrate robust performance across almost all test configurations, with speed and accuracy levels that compare competitively, and in fact often favorably, to Group II (global) algorithms.  Among Group II (global) algorithms, GenSA appears to be the only one that shows overall better performance than our SMCO algorithms in Table \ref{tab:Comp_200d}. That said, it still yields noticeably larger errors (in terms of RMSE and AE99) than ours, say, in the minimization of Rastrigin and Ackley functions. 
Additional results in higher dimensions $d= 400, 500$ and $1000$, as well as various sensitivity checks with respect to key hyperparameters of SMCO/-R/-BR, are available in Appendix \ref{append-D}.


\begin{sidewaystable}
\caption{\label{tab:Comp_200d} Deterministic  Functions, Dimension $d = 200$: Comparison with Group II Algorithms}
\begin{threeparttable}
\begin{tabular}{|c|cccc|ccc|ccc|ccc|}
\hline 
\multicolumn{2}{|c|}{Test Function} & \multicolumn{3}{c|}{\textbf{Rastrigin}} & \multicolumn{3}{c|}{\textbf{Griewank}} & \multicolumn{3}{c|}{\textbf{Ackley}} & \multicolumn{3}{c|}{\textbf{Michalewicz}}\tabularnewline
\hline 
\multicolumn{2}{|c|}{Performance} & \multicolumn{1}{c|}{RMSE} & AE99 & Time & \multicolumn{1}{c|}{RMSE} & \multicolumn{1}{c|}{AE99} & Time & \multicolumn{1}{c|}{RMSE} & \multicolumn{1}{c|}{AE99} & Time & \multicolumn{1}{c|}{RMSE} & AE99 & Time\tabularnewline
\hline 
\multirow{10}{*}{Min} & \textbf{Best Value} & \multicolumn{3}{c|}{\textbf{1596.819 }} & \multicolumn{3}{c|}{\textbf{2.35e-13 }} & \multicolumn{3}{c|}{\textbf{19.85938 }} & \multicolumn{3}{c|}{\textbf{-62.77543 }}\tabularnewline
\cline{3-14} \cline{4-14} \cline{5-14} \cline{6-14} \cline{7-14} \cline{8-14} \cline{9-14} \cline{10-14} \cline{11-14} \cline{12-14} \cline{13-14} \cline{14-14} 
 & \textbf{SMCOms} & \textbf{285 } & \textbf{375 } & \textbf{34.9 } & \textbf{1.07 } & \textbf{1.07 } & \textbf{34.4 } & \textbf{0.101 } & \textbf{0.106 } & \textbf{31.9 } & \textbf{42.22 } & \textbf{44.08 } & \textbf{42.7 }\tabularnewline
 & \textbf{SMCOms-R} & \textbf{262 } & \textbf{308 } & \textbf{33.4 } & \textbf{0.95 } & \textbf{0.95 } & \textbf{37.2 } & \textbf{0.009 } & \textbf{0.014 } & \textbf{32.3 } & \textbf{28.42 } & \textbf{29.41 } & \textbf{42.5 }\tabularnewline
 & \textbf{SMCOms-BR} & \textbf{272 } & \textbf{324 } & \textbf{32.6 } & \textbf{0.96 } & \textbf{0.96 } & \textbf{36.5 } & \textbf{0.010 } & \textbf{0.015 } & \textbf{32.6 } & \textbf{30.33 } & \textbf{32.21 } & \textbf{43.0 }\tabularnewline
\cline{3-14} \cline{4-14} \cline{5-14} \cline{6-14} \cline{7-14} \cline{8-14} \cline{9-14} \cline{10-14} \cline{11-14} \cline{12-14} \cline{13-14} \cline{14-14} 
 & GenSA & 1718  & 2011  & 268.3  & 3.05e-11  & 6.10e-11  & 302.2  & 0.107  & 0.114  & 271.7  & 6.27  & 12.03  & 813.4 \tabularnewline
 & SA & 3947  & 4338  & 0.2  & 11626  & 14002  & 0.1  & 1.651  & 1.679  & 0.3  & 52.76  & 54.53  & 0.36 \tabularnewline
 & DEoptim & 3813  & 3968  & 31.5  & 5856  & 6086  & 33.4  & 1.545  & 1.563  & 32.2  & 47.41  & 48.91  & 43.5 \tabularnewline
 & STOGO & 1293  & 1293  & 285.6  & 4.38e-09  & 4.38e-09  & 307.1  & 0  & 0 & 295.4  & 43.48  & 43.48  & 404.4 \tabularnewline
 & GA & 2496  & 2623  & 0.5  & 7372  & 7920  & 0.5  & 1.450  & 1.518  & 0.5  & 49.28  & 52.30  & 0.61 \tabularnewline
 & PSO & 590 & 812  & 4.9  & 81  & 143 & 5.1  & 1.546  & 1.590  & 5.1 & 43.13  & 45.36  & 6.19 \tabularnewline
\hline 
\multirow{10}{*}{Max} & \textbf{Best Value} & \multicolumn{3}{c|}{\textbf{18239.57 }} & \multicolumn{3}{c|}{\textbf{72331.4}} & \multicolumn{3}{c|}{\textbf{22.33939 }} & \multicolumn{3}{c|}{\textbf{69.53472 }}\tabularnewline
\cline{3-14} \cline{4-14} \cline{5-14} \cline{6-14} \cline{7-14} \cline{8-14} \cline{9-14} \cline{10-14} \cline{11-14} \cline{12-14} \cline{13-14} \cline{14-14} 
 & \textbf{SMCOms} & \textbf{7328 } & \textbf{7648 } & \textbf{34.7 } & \textbf{192 } & \textbf{199 } & \textbf{38.0 } & \textbf{0.119 } & \textbf{0.125 } & \textbf{31.6 } & \textbf{36.67 } & \textbf{38.25 } & \textbf{42.48 }\tabularnewline
 & \textbf{SMCOms-R} & \textbf{4419 } & \textbf{4939 } & \textbf{33.3 } & \textbf{278 } & \textbf{287 } & \textbf{36.2 } & \textbf{0.105 } & \textbf{0.110 } & \textbf{30.8 } & \textbf{20.91 } & \textbf{22.16} & \textbf{42.47 }\tabularnewline
 & \textbf{SMCOms-BR} & \textbf{4512 } & \textbf{5043 } & \textbf{34.7 } & \textbf{601 } & \textbf{627 } & \textbf{37.4 } & \textbf{0.030 } & \textbf{0.034 } & \textbf{31.7 } & \textbf{25.43 } & \textbf{26.87 } & \textbf{42.49 }\tabularnewline
\cline{3-14} \cline{4-14} \cline{5-14} \cline{6-14} \cline{7-14} \cline{8-14} \cline{9-14} \cline{10-14} \cline{11-14} \cline{12-14} \cline{13-14} \cline{14-14} 
 & GenSA & 2605  & 2975  & 275.9  & 0  & 0  & 312.7  & 0.004  & 0.007  & 273.5  & 4.65  & 8.98  & 814.1 \tabularnewline
 & SA & 5653  & 7180  & 0.3  & 32221  & 36635  & 0.3  & 0.488  & 0.505  & 0.3  & 55.48  & 57.14  & 0.37 \tabularnewline
 & DEoptim & 7978  & 8191  & 31.8  & 36012  & 36465  & 33.7  & 0.421  & 0.432 & 32.6  & 49.83  & 50.98  & 43.8 \tabularnewline
 & STOGO & 10987  & 10987  & 290.0  & 9908  & 9908  & 310.0  & 0.120  & 0.120  & 299.9  & 45.26  & 45.26  & 407.8 \tabularnewline
 & GA & 8977  & 9310  & 0.5  & 45978  & 46953  & 0.5  & 0.450  & 0.489  & 0.5  & 49.86  & 51.20  & 0.61 \tabularnewline
 & PSO & 998 & 1719  & 5.0  & 20792  & 23940  & 5.2  & 0.429  & 0.448  & 5.2  & 46.50 & 48.88  & 6.23 \tabularnewline
\hline 
\end{tabular} 
\begin{tablenotes}
\item Based on 10 reps. SMCO/-R/-BR run from \textbf{$\sqrt{d}\approx14$ uniformly drawn starting points} in parallel. Run on high performance cluster using 15 cores on 1 node.
\end{tablenotes}
\end{threeparttable}
\end{sidewaystable}

\subsection{Optimization of Random  Functions}

In this section, we investigate two randomly generated test functions that arise from two important economic contexts. The first test function is the sample criterion for maximum score estimation, as introduced by Manski \cite{manski1975maximum,manski1985semiparametric}. The maximum score estimator can be used to estimate model parameters in discrete choice models under weak assumptions on the error terms, and since \cite{manski1975maximum,manski1985semiparametric} there has been much subsequent work in economics that either extends the applicability of the maximum score estimation approach or analyzes its intriguing properties.  Specifically, the sample criterion function for maximum score estimator is defined as
 $f_{MS}\left(\b\right):=\frac{1}{N}\sum_{i=1}^{N}Y_{i}\ind\left\{ X_{i1}+X_{i,-1}^{'}\b\geq0\right\} $
where $\left(X_{i},Y_{i}\right)_{i=1}^{N}$ is a random sample with
$X_{i}=\left(X_{i,1},X_{i,-1}\right)\in\R^{d+1}$ being a vector of
covariates and $Y_{i}$ being a binary outcome variable generated
by the model 
$Y_{i}=\ind\left\{ X_{i1}+X_{i,-1}^{'}\b_{0}+\e_{i}\geq0\right\}$ 
under the conditional median restriction $\text{med}\left(\rest{\e_{i}}X_{i}\right)=0$
at the true parameter $\b_{0}$. The maximum score estimator is then defined as any solution to the maximization problem
$\max_{\b}f_{MS}\left(\b\right).$
Note that the parameter to be optimized over, $\b$, enters through the discontinuous indicator function, and furthermore the outcome variable $Y_i$ is also a binary variable itself. Due to the discreteness issues, the maximization of the maximum score criterion function is well-known to be computationally challenging, especially when the dimension of parameters becomes (moderately) high. We investigate the performance of our SMCO algorithms in optimizing the maximum score criterion function. Specifically, we build the test function $f_{MS}$ using simulated data generated
as follows. We simulate $X_{i,1}\sim\mathcal{N}\left(0,d\right)$ and $X_{i,2}\sim...\sim X_{i,d+1}\sim\mathcal{N}\left(0,1\right)$.
We set $\b_{0}\in\R^{d}$ through random draws of $\b_{0,j}\sim\mathcal{N}\left(0,1\right)$.
We simulate $\e_{i}\sim\mathcal{N}\left(0,d\right)$ and set $N=500$. The variances of $X_{i,1}$ and $\e_{i}$ are set equal to $d$ so that
the effects from $X_{i,1}$ and $\e_{i}$ will not be dominated by
those from $X_{i,-1}\in\R^{d}$ as the dimension $d$ increases. We consider $d=5, 10, 20 $ and $50$, and constrain the parameter search in the rectangle $[-20,20]^d.$

The second randomly generated test function we consider arises from empirical welfare maximization as introduced in \cite{kitagawa2018should}, which has inspired a very active line of related research in economics on optimal treatment assignment and policy evaluation. Specifically, the empirical welfare function, say, under linear treatment assignment, is defined in \cite{kitagawa2018should} as
$
f_{EW}\left(\b\right):=\frac{1}{N}\sum_{i=1}^{N}\left(\frac{D_{i}}{p\left(X_{i}\right)}-\frac{1-D_{i}}{1-p\left(X_{i}\right)}\right)Y_{i}\ind\left\{ \b_{1}+X_{i}^{'}\b_{-1}\geq0\right\} $
where $\left(D_{i},X_{i},Y_{i}\right)_{i=1}^{N}$ is a random sample
with $D_{i}$ being a binary treatment status variable with (known/chosen) propensity scores $p\left(X_{i}\right):=\E\left[\rest{D_{i}}X_{i}\right]$,
$X_{i}\in\R^{d-1}$ being a vector of covariates, and $Y_{i}$ is
a scalar-valued outcome variable. The optimal
linear treatment assignment rule is then obtained by solving the optimization problem $\max_{\b}f_{EW}\left(\b\right)$. 
Similarly to the  maximum score criterion function, $f_{EW}$ also features the discontinuous indicator function that makes the optimization problem challenging. 
We also investigated the performance of the SMCO algorithms in solving empirical welfare maximization problems. Specifically, we build the test function $f_{MS}$ based on a combination of real and simulated data. The real data are constructed based on the replication package of \cite{kitagawa2018should}, which is constructed from the Job Training Partnership Act (JTPA) data on job training programs. The specific data set we use contains $N=9223$ observations on the binary treatment status variable $D_{i}$,
the outcome variable $Y_{i}$ (post-program earning), and two covariates:
pre-program earning ($X_{i1}$) and education ($X_{i2}$). The propensity
score is known to be constant at $p\left(X_{i}\right) = 2/3$.
Following \cite{kitagawa2018should}, we set $X_{i3}=X_{i2}^{2}$ and $X_{i4}=X_{i2}^{3}$
to be quadratic and cubic transformation of education. We then standardize
$X_{i1},...,X_{i4}$ by demeaning it and rescaling its standard deviation
to 1. We also augment the data by simulating additional covariates
$X_{i5},...,X_{i,d-1}\sim\cN\left(0,1\right)$ for higher dimensions
$d\geq6$.  Again, we consider $d=5, 10, 20 $ and $50$, and constrain the parameter space to $[-20,20]^d.$

\begin{sidewaystable}
\centering
\caption{Maximization of Random Functions, Dimension $d=5,10,20$ \& $50$\label{tab:MSEW}}
\begin{threeparttable}
\begin{tabular}{|ccccc|ccc|ccc|ccc|}
\hline 
\multicolumn{2}{|c}{Dimension} & \multicolumn{3}{c|}{\textbf{$d=5$}} & \multicolumn{3}{c|}{\textbf{$d=10$}} & \multicolumn{3}{c|}{\textbf{$d=20$}} & \multicolumn{3}{c|}{\textbf{$d=50$}}\tabularnewline
\hline 
\multicolumn{2}{|c}{Performance} & \multicolumn{1}{c|}{RMSE} & \multicolumn{1}{c|}{AE99} & Time & \multicolumn{1}{c|}{RMSE} & \multicolumn{1}{c|}{AE99} & Time & \multicolumn{1}{c|}{RMSE} & \multicolumn{1}{c|}{AE99} & Time & \multicolumn{1}{c|}{RMSE} & \multicolumn{1}{c|}{AE99} & Time\tabularnewline
\hline 
\multicolumn{14}{|c|}{\textbf{Maximum Score Estimation: $\max_{\b}f_{MS}\left(\b\right)\phantom{\frac{\frac{1}{1}}{\frac{1}{2}}}$}}\tabularnewline
\hline 
 & \textbf{Best Value} & \multicolumn{3}{c|}{\textbf{0.42 }} & \multicolumn{3}{c|}{\textbf{0.448 }} & \multicolumn{3}{c|}{\textbf{0.462 }} & \multicolumn{3}{c|}{\textbf{0.424 }}\tabularnewline
\hline 
\multirow{6}{*}{(I)} & GD & 0.135  & 0.196  & 0.06  & 0.146  & 0.190  & 0.03  & 0.167  & 0.198  & 0.08 & 0.157  & 0.174  & 0.13 \tabularnewline
 & SignGD & 0.135  & 0.196  & 0.08  & 0.146  & 0.190  & 0.03  & 0.167  & 0.198  & 0.09  & 0.157  & 0.174  & 0.11 \tabularnewline
 & ADAM & 0.135  & 0.196  & 0.07  & 0.146  & 0.190  & 0.06  & 0.167  & 0.198  & 0.08  & 0.157  & 0.174  & 0.11 \tabularnewline
 & SPSA & 0.135  & 0.196  & 0.11  & 0.146  & 0.190  & 0.08 & 0.167  & 0.198  & 0.11  & 0.157  & 0.174  & 0.12 \tabularnewline
 & L-BFGS & 0.135  & 0.196  & 0.07  & 0.146  & 0.190  & 0.07  & 0.167  & 0.198  & 0.09  & 0.157  & 0.174  & 0.13 \tabularnewline
 & BOBYQA & 0.084  & 0.116  & 0.07  & 0.100  & 0.141  & 0.08  & 0.130  & 0.153  & 0.09  & 0.134  & 0.155  & 0.22 \tabularnewline
\hline 
\multirow{3}{*}{} & \textbf{SMCO} & \textbf{0.008 } & \textbf{0.012 } & \textbf{0.20 } & \textbf{0.009 } & \textbf{0.012 } & \textbf{0.33 } & \textbf{0.013 } & \textbf{0.019 } & \textbf{0.76 } & \textbf{0.056 } & \textbf{0.068 } & \textbf{5.09 }\tabularnewline
 & \textbf{SMCO-R} & \textbf{0.002 } & \textbf{0.004 } & \textbf{0.21 } & \textbf{0.005 } & \textbf{0.006 } & \textbf{0.33 } & \textbf{0.006 } & \textbf{0.008 } & \textbf{0.77 } & \textbf{0.010 } & \textbf{0.016 } & \textbf{5.19 }\tabularnewline
 & \textbf{SMCO-BR} & \textbf{0.003 } & \textbf{0.008 } & \textbf{0.21 } & \textbf{0.007 } & \textbf{0.010 } & \textbf{0.35 } & \textbf{0.010 } & \textbf{0.015 } & \textbf{0.81 } & \textbf{0.010 } & \textbf{0.017 } & \textbf{5.41 }\tabularnewline
\hline 
\multirow{6}{*}{(II)} & GenSA & 0.020  & 0.054  & 1.71  & 0.031  & 0.034  & 4.35  & 0.032  & 0.038  & 15.69  & 0.010  & 0.018  & 89.93 \tabularnewline
 & SA & 0.069  & 0.076  & 0.06  & 0.073  & 0.085  & 0.07  & 0.108  & 0.124  & 0.14  & 0.114  & 0.124  & 0.34 \tabularnewline
 & DEoptim & 0.013  & 0.054  & 0.30  & 0.036  & 0.038  & 0.76  & 0.045  & 0.050  & 2.78  & 0.056  & 0.064  & 13.18 \tabularnewline
 & STOGO & 0.052  & 0.052  & 3.49  & 0.068  & 0.068  & 8.18  & 0.080  & 0.080  & 28.43  & 0.082  & 0.082  & 182.28 \tabularnewline
 & GA & 0.050  & 0.058  & 0.22  & 0.036  & 0.049  & 0.25  & 0.054  & 0.074  & 0.38  & 0.059  & 0.070  & 0.65 \tabularnewline
 & PSO & 0.040  & 0.056  & 0.65  & 0.034  & 0.038  & 0.90  & 0.038  & 0.046  & 1.64  & 0.017  & 0.026  & 4.32 \tabularnewline
\hline 
\multicolumn{14}{|c|}{\textbf{Empirical Welfare Maximization: $\max _\b f_{EW}\left(\b\right)\phantom{\frac{\frac{1}{1}}{\frac{1}{2}}}$}}\tabularnewline
\hline 
 & \textbf{Best Value} & \multicolumn{3}{c|}{\textbf{1541.032 }} & \multicolumn{3}{c|}{\textbf{2061.409 }} & \multicolumn{3}{c|}{\textbf{2649.911 }} & \multicolumn{3}{c|}{\textbf{3825.49 }}\tabularnewline
\hline 
\multirow{6}{*}{(I)} & GD & 752  & 1448  & 0.12  & 1288  & 1666  & 0.16  & 1827  & 2058  & 0.24  & 2934  & 3067  & 0.60 \tabularnewline
 & SignGD & 752  & 1448  & 0.13  & 1288  & 1666  & 0.17  & 1827  & 2058  & 0.25  & 2934  & 3067  & 0.60 \tabularnewline
 & ADAM & 752  & 1448  & 0.13  & 1288  & 1666  & 0.16  & 1827  & 2058  & 0.25  & 2934  & 3067  & 0.62 \tabularnewline
 & SPSA & 752  & 1448  & 0.11  & 1305  & 1711  & 0.14  & 1848  & 2114  & 0.21  & 2920  & 3063  & 0.69 \tabularnewline
 & L-BFGS & 721  & 1448  & 0.12  & 1199  & 1621  & 0.18  & 1799  & 2050  & 0.28  & 2927  & 3067  & 0.99 \tabularnewline
 & BOBYQA & 398  & 955  & 0.13  & 846  & 1267  & 0.20  & 1365  & 1622  & 0.40  & 2338  & 2599  & 1.51 \tabularnewline
\hline 
\multirow{3}{*}{} & \textbf{SMCO} & \textbf{207 } & \textbf{301 } & \textbf{1.50 } & \textbf{195 } & \textbf{310 } & \textbf{3.41 } & \textbf{271 } & \textbf{379 } & \textbf{9.49 } & \textbf{428 } & \textbf{498 } & \textbf{78.37 }\tabularnewline
 & \textbf{SMCO-R} & \textbf{128 } & \textbf{230 } & \textbf{1.51 } & \textbf{69 } & \textbf{141 } & \textbf{3.41 } & \textbf{98 } & \textbf{171 } & \textbf{9.52 } & \textbf{158 } & \textbf{233 } & \textbf{78.59 }\tabularnewline
 & \textbf{SMCO-BR} & \textbf{140 } & \textbf{229 } & \textbf{1.55 } & \textbf{85 } & \textbf{198 } & \textbf{3.46 } & \textbf{144 } & \textbf{252 } & \textbf{9.57 } & \textbf{274 } & \textbf{361 } & \textbf{78.71 }\tabularnewline
\hline 
\multirow{6}{*}{(II)} & GenSA & 104  & 150  & 13.37  & 220  & 368  & 32.37  & 639  & 918  & 91.25  & 1586  & 1762  & 455.24 \tabularnewline
 & SA & 214  & 318  & 0.41  & 431  & 638  & 0.56  & 881  & 1268  & 0.85  & 1939  & 2398  & 1.57 \tabularnewline
 & DEoptim & 116  & 161  & 2.37  & 233  & 279  & 6.11  & 640  & 723  & 17.54  & 1722  & 1805  & 87.03 \tabularnewline
 & STOGO & 134  & 134  & 27.94  & 470  & 470  & 66.33  & 851  & 851  & 181.67  & 2096  & 2096  & 891.62 \tabularnewline
 & GA & 160  & 238  & 1.08  & 288  & 453  & 1.36  & 689  & 973  & 1.90  & 1685  & 1831  & 3.69 \tabularnewline
 & PSO & 108 & 192  & 3.60  & 138  & 272  & 5.26  & 297  & 433  & 8.26  & 848  & 1051  & 21.52 \tabularnewline
\hline 
\end{tabular}
\begin{tablenotes}
\item Based on $500/d = 100, 50, 25, 10$ reps. SMCO/-R/-BR and Group (I) algorithms run from \textbf{$\sqrt{d} \approx  2, 3, 4, 7$ uniformly drawn starting points} in parallel. Run on high performance cluster  with 8 cores on 1 node.
\end{tablenotes}
\end{threeparttable}
\end{sidewaystable}

Table \ref{tab:MSEW} reports the performance of our SMCO algorithms, along with the competing algorithms, in maximizing the two randomly generated test functions $f_{MS}$ and $F_{EW}$ as described above. First, all gradient-based local algorithms perform significantly worse than our SMCO algorithms and the global algorithms: this is expected given that the gradients of $f_{MS}$ and $F_{EW}$ are both degenerate almost everywhere given that the parameter enters through the indicator function. Second, our SMCO algorithms (especially SMCO-R/BR) perform robustly well across almost all the configruations in Table \ref{tab:MSEW}: in fact, SMCO-R/BR yields the best results among all algorithms except for one configuration ($f_{EW}$ with $d=5$). 

Appendix \ref{append-E} contains supplemental numerical results on maximum score estimation and empirical welfare maximization up to $d=100$, including replication of Table \ref{tab:MSEW} using 10 replications only, as well as supplemental robustness checks of maximum score estimation with Cauchy errors and different bound buffer parameters of SMCO/-R/-BR.

\subsection{Additional Numerical Experiments} We have conducted many more numerical experiments to check the wide applicability of our SMCO algorithms. Some of which are reported in Appendices \ref{append-F} and \ref{append-G}. Specifically, see: 
Appendix \ref{append-F} for results on two additional random test functions: one based on conditional moment inequalities in a dynamic binary choice model from economic contexts, the other based on mean squared in-class approximation error for ReLU neural networks; 
Appendix \ref{append-G} for results on additional random test functions in statistics, such as Cauchy log-likelihood, ReLU neural network out of sample predictions,  ReLU neural network out of sample misclassification problems. This section contains statistical performance measurements such as ``out-of-sample'' MSEs and the ``out-of-sample'' misclassification error rate. These are different from but complementary to the optimization performance measurements we focus on in the main text. We have also tested many other deterministic test functions, such as Bukin No.6, Damavandi, Dropwave, Easom, EggHolder, JennrichSampson, McCormick, Mishra No.6, Qing, Rosenbrock, Schwefel, Shubert, Six-Hump Camel, Trefethen, and Zakharov functions, with dimensions up to 200 (when applicable). We do not report them here due to the lack of space. Such results are available upon request, and can also be easily replicated using our SMCO source code on GitHub. The main findings from all these exercises are that, our SMCO algorithms, especially SMCO-R and SMCO-BR, perform very well across all these different configurations, essentially without any change in the hyper-parameters other than the dimension-adaptive number of starting points used. In fact, this observation suggests that the suitability of the SMCO algorithms for global optimization goes well beyond the theoretical guarantees we provide in Section \ref{sec-theory}.


\section{Conclusion}\label{sec:end}

Based on the Strategic Law of Large Numbers we establish, this paper proposes a unified two-armed decision model for black-box optimization. We transform the difficult problem of global optimization over a bounded rectangle region to optimal strategy formation of multiple two-armed decision models in infinite policy sets. 
We propose a novel class of Strategic Monte Carlo Optimization (SMCO) algorithms, for which we establish local and global convergence results. 
We also show via an extensive array of numerical experiments that our proposed SMCO algorithms 
perform consistently well across different complex multi-modal test functions and different dimensions (up to $d=1000$), often outperforming a portfolio of state-of-the-art existing algorithms, even though our SMCO algorithms are currently not coded in C++. Our work opens the door to a promising area of future research that explores and improves the theoretical properties as well as the practical implementation of optimization algorithms under our StLLN-based framework.

\subsection*{Data, Materials, and Software Availability} 
The Job Training Partnership Act (JTPA) data used in Section 3.4 can be constructed based on the replication package of \cite{kitagawa2018should}, which is publicly available at: \\ https://doi.org/10.3982/ECTA13288. 
The R code used to generate the numerical studies in Section 3 and the Appendix are freely available for academic use on Github at https://github.com/wayne-y-gao/SMCO and
https://github.com/yanxiaodong128/SMCO

\subsection*{Acknowledgments} 
We thank Michael Fu, Jiaqiao Hu, Faming Liang, Shige Peng, Andre Wibisono and Yinyu Ye for helpful discussions. This work was supported in part by the National Key R$\&$D Program of China (Grant No. 2023YFA1008701), the  Key Project of the National Natural Science Foundation of China (No. 12431017). the National Natural Science Foundation of China (Grant No. 12371292, 11901352, 12371148). Xiaohong Chen's research is partially supported by Cowles Foundation for Research in Economics, Yale University.


\bibliography{pnas-sample}
\bibliographystyle{ecta}

\newpage
\appendix

~
\begin{center}
    \textbf{\Large Appendix}
\end{center}

~

We first present a graphic illustration of our two-armed decision model as well as the difference between our novel SMCO algorithm against the standard gradient descent algorithm in one dimension:
\begin{figure}[H] 
\centering  
  \captionsetup[subfloat]{labelfont=scriptsize,textfont=tiny}
  \subfloat[Classical Monte Carlo]  
  {
      \label{Classical Monte Carlo (one arm)}\includegraphics[width=0.15\textwidth]{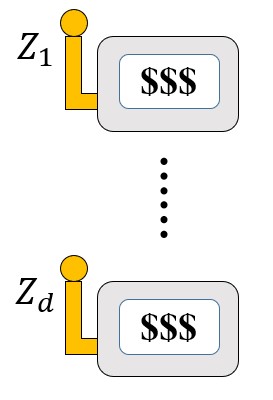}
  }
  \subfloat[Strategic Monte Carlo]  
  {
      \label{Strategic Monte Carlo(two arms)}\includegraphics[width=0.175\textwidth]{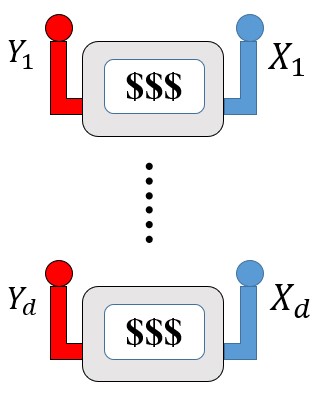}
  }
  \subfloat[Iterative Mechanism by GDA and SMCO]  
  {
      \label{Iterative Mechanism }\includegraphics[width=0.45\textwidth]{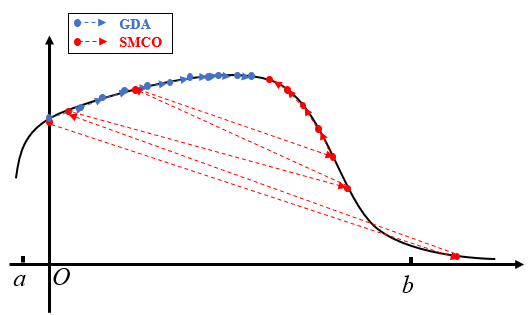}
  }
 \caption{ (a). Classic Monte Carlo sampling from a sequence of one-armed slot machines, with arm population $\{Z_1,\cdots, Z_d\}$; (b). Strategic Monte Carlo sampling
   from a sequence of two-armed slot machines, with arms populations $\{X_1,\cdots, X_d\}$ and $\{Y_1,\cdots, Y_d\}$; (c) The iterative mechanism for attaining the extreme value by GD (blue arrow line) and SMCO (red arrow line) under $d=1$, where red points denote the iterative value $S_{n}^{\bmtheta^*}/n$ for SMCO given in Algorithm \ref{alg:smco_algo}.}
  \label{cmc_smc}
\end{figure}

The Appendix contains several sections.  Appendix \ref{append-A} contains additional lemmas and   proofs of the mathematical results in Section \ref{sec-theory} of the main paper. Appendix \ref{append-B}  provides more details on the practical implementations of the SMCO algorithms. Appendix \ref{append-C} shows the details on comparsion-group algorithms. Appendix \ref{append-D} presents numerical experiments results that supplement those in Section 3.3 on deterministic test functions. Appendix \ref{append-E} provides results that supplement those in Section 3.4 on randomly generated test functions. Appendix \ref{append-F} presents results on two additional randomly generated functions constructed from ReLU neural networks. Appendix \ref{append-G}  presents comparisons using a variety of statistical random objective functions.

\section{Lemmas and Proofs}\label{append-A}
For convenience, we assume the rectangular domain $\Gamma=\prod_{i=1}^d[\underline{\mu}_{i},\overline{\mu}_i]$ is symmetric about the origin (i.e. $\underline{\mu}_i = -\overline{\mu}_i$ for all $i$). If not satisfied, apply the translation $T(x) = x - t$ where $t = \left(\frac{\underline{\mu}_1+\overline{\mu}_1}{2}, \dots, \frac{\underline{\mu}_d+\overline{\mu}_d}{2}\right)$. This transformation yields a shifted domain $\widetilde{\Gamma} = \prod_{i=1}^d[\widetilde{\underline{\mu}}_{i},\widetilde{\overline{\mu}}_{i}]$, where $\widetilde{\underline{\mu}}_{i}=\underline{\mu}_{i}-t_i,\widetilde{\overline{\mu}}_{i}=\overline{\mu}_i-t_i$, satisfying the symmetry condition, and a shifted function $\tilde{f}(x) = f(x + t)$. Consequently, for $\widetilde{\Gamma}$ we have $\widetilde{\underline{\mu}}_{i}+\widetilde{\overline{\mu}}_{i}=0$ and $\widetilde{\overline{\mu}}_{i}=-\widetilde{\underline{\mu}}_{i}=\frac{\widetilde{\overline{\mu}}_{i}-\widetilde{\underline{\mu}}_{i}}{2}$
for all $i$, ensuring origin symmetry.

We first present two lemmas that will be used in the rest of the proofs. 
Let $\psi\in C_{b,Lip}(\mathbb{R}^d)$ be a bounded Lipschitz continuous function on $\mathbb{R}^d$. 
For any $\epsilon>0$ and $t\in \lbrack 0,1+\epsilon]$, define
    \begin{equation}\label{function-H}
        H_t^{\epsilon}(x):=u_{\epsilon}(t,x):=\sup_{a\in\mathcal{A}(t,1+\epsilon)}E\left[\psi\left(x+\int_{t}^{1+\epsilon}a_sds+\int_t^{1+\epsilon}\epsilon\boldsymbol{I}_ddB_s\right)\right]
    \end{equation}
    where $\boldsymbol{I}_d$ is a $d$-dimensional identity matrix and $(B_t)$ is a $d$-dimensional Brownian motion  and $\mathcal{F}_t$ is the natural filtration generated by $(B_t)$. Let $\mathcal{A}(t,T)$ denote the set of all $\{\mathcal{F}_{s}\}$%
-adapted measurable processes valued in $\Gamma$.

\begin{Lem}
\label{lemma-ddp-0} (a) $\max_{x\in\Gamma}\psi(x)=\sup_{a\in\mathcal{A}(0,1)}E\left[\psi\left(\int_0^{1}a_sds\right)\right]$; for any $t\in[0,1],x\in \mathbb{R}^d$,
 $$u(t,x):=\sup_{a\in\mathcal{A}(t,1)}E\left[\psi\left(x+\int_{t}^{1}a_sds\right)\right]=\max_{y\in\Gamma}\psi(x+(1-t)y).$$

(b) There are finite constants $c_1,c_2>0$ such that 
$$\sup_{x\in \mathbb{R}^d} |u_{\epsilon}(t,x)-u(t,x) |\leq \epsilon \times [c_1+c_2 (d(\epsilon + 1-t))^{1/2}]~~~\text{uniformly over}~t\in[0,1].$$
\newline
(c) For any fixed $\epsilon>0$, the functions $\{H_t^{\epsilon}(x)\}_{t\in \lbrack 0,1]}$ satisfy the
following properties:

\begin{description}
\item[(1)] $H_t^{\epsilon}\in C_{b}^{2}(\mathbb{R}^d)$ and each element of the gradient vector and Hessian matrix of $H_t^{\epsilon}$ are uniformly bounded for all $t\in \lbrack 0,1]$.

\item[(2)] Dynamic programming principle: For any $\gamma \in \lbrack
0,1+\epsilon-t]$,
\begin{equation*}
H_t^{\epsilon}(x) =\sup_{a\in  \mathcal{A}(t,t+\gamma )}E
\left[ H_{t+\gamma}^{\epsilon}\left(x+\int_{t}^{t+\gamma
}a_{s}ds+\int_t^{t+\gamma}\epsilon\boldsymbol{I}_ddB_s\right) \right]
,\;x\in \mathbb{R}^d.
\end{equation*}

\item[(3)]There exists  a constant $C$, depends only on  the uniform bound of $D_{x}^2\, H_t^{\epsilon}$, such that for any $n\in\mathbb{N}$,
\begin{equation*}
\sum\limits_{m=1}^{n}\sup\limits_{x\in
\mathbb{R}^d}\left\vert H_{\frac{m-1}{n}}^{\epsilon}(x)-H_{\frac{m}{n}}^{\epsilon}\left(x\right) -\tfrac{1}{n}\sup_{p\in\Gamma}\left[ D_x H_{\frac{m}{n}}^{\epsilon}\left(
x\right)\cdot p\right]
\right\vert \le C\left(\frac{\sum_{i=1}^d|\overline{\mu}_i-\underline{\mu}_i|}{2n}+2\epsilon\sqrt{\frac{d}{n}}+\epsilon^2\right).
\end{equation*}

\end{description}
\end{Lem}

\noindent \textbf{Proof of Lemma \ref{lemma-ddp-0}:} 
For Result (a), we note that for any $x\in \Gamma$, let $a_s^{x}=x$ for $s\in[0,1]$, we have $a_s^{x}\in\mathcal{A}(0,1)$, then
$$
\sup_{a\in\mathcal{A}(0,1)}E\left[\psi\left(\int_0^{1}a_sds\right)\right]\geq
\sup_{x\in\Gamma}E\left[\psi\left(\int_0^{1}a_s^xds\right)\right]=\sup_{x\in\Gamma}E\left[\psi\left(x\right)\right] =\max_{x\in\Gamma}\psi(x).
$$
In addition, since any $a\in\mathcal{A}(0,1)$ takes values in $\Gamma$, we have
$$
\sup_{a\in\mathcal{A}(0,1)}\psi\left(\int_0^{1}a_s(\omega)ds\right)\leq \max_{x\in\Gamma}\psi(x), \quad \omega \in \Omega.
$$
Therefore,
$$
\sup_{a\in\mathcal{A}(0,1)}E\left[\psi\left(\int_0^{1}a_sds\right)\right]\leq
\sup_{a\in\mathcal{A}(0,1)}E\left[\sup_{a\in\mathcal{A}(0,1)}\psi\left(\int_0^{1}a_sds\right)\right]
\leq\max_{x\in\Gamma}\psi(x).
$$
Given the definition of $u(t,x)$, the relation $u(t,x)=\max_{y\in\Gamma}\psi(x+(1-t)y)$ can be established similarly.

For Result (b), follows from Result (a), we have uniformly over $t\in [0,1],~x\in \mathbb{R}^d$, 
\begin{align*}
   & |u_{\epsilon}(t,x)-u(t,x)|\\
    =&\left|\sup_{a\in\mathcal{A}(t,1+\epsilon)}E\left[\psi\left(x+\int_{t}^{1+\epsilon}a_sds+\int_t^{1+\epsilon}\epsilon\boldsymbol{I}_ddB_s\right)\right]-\sup_{a\in\mathcal{A}(t,1)}E\left[\psi\left(x+\int_{t}^{1}a_sds\right)\right]\right|\\
    \leq&\sup_{a\in\mathcal{A}(t,1+\epsilon)}E\left[\left|\psi\left(x+\int_{t}^{1+\epsilon}a_sds+\int_t^{1+\epsilon}\epsilon\boldsymbol{I}_ddB_s\right)-\psi\left(x+\int_{t}^{1}a_sds\right)\right|\right]\\
    \leq&L_{\psi}\sup_{a\in\mathcal{A}(t,1+\epsilon)}E\left[\left\|\int_{1}^{1+\epsilon}a_sds\right\|+\left\|\int_t^{1+\epsilon}\epsilon\boldsymbol{I}_ddB_s\right\|\right]\\
    \leq& L_{\psi}\epsilon\left( \sqrt{\sum_{i=1}^d\frac{(\overline{\mu}_i-\underline{\mu}_i)^2}{4}}
    +\sqrt{d(1+\epsilon-t)}\right).
\end{align*}
For Result (c) Part (1), for any $(t,x)\in \lbrack 0,1+\epsilon]\times \mathbb{R}^d $, let $u_{\epsilon}(t,x)=H_{t}^{\epsilon}(x)$,  it can be checked that $u_{\epsilon}$ is the solution of the PDE \eqref{pde-proof1} (HJB-equation, see Yong and Zhou \cite[Theorem 5.2, Ch. 4]{yongzhou1999}).
\begin{equation}
\left\{
\begin{array}{l}
\displaystyle\partial _{t}u_{\epsilon}(t,x)+ \sup_{p\in\Gamma}\left[D_{x}u_{\epsilon}(t,x)\cdot p\right] +\frac{\epsilon^2}{2} \textrm{tr}\left[D_{x}^2u_{\epsilon}(t,x) \right]=0,\quad (t,x)\in \lbrack 0,1+\epsilon)\times \mathbb{R}^d \\
u_{\epsilon}(1+\epsilon,x)=\psi(x).%
\end{array}%
\right.   \label{pde-proof1}
\end{equation}
By Theorem C.4.5 in Peng \cite{peng2019}, $\exists \beta \in (0,1)$ such that
\begin{equation*}
\Vert u_{\epsilon}\Vert _{C^{1+\beta /2,2+\beta }([0,1]\times \mathbb{R}^{d})}<\infty .
\end{equation*}%
Here $\Vert \cdot \Vert _{C^{1+\beta /2,2+\beta }([0,1]\times
\mathbb{R}^{d})}$ is the Krylov \cite{krylov} norm on $C^{1+\beta /2,2+\beta
}([0,1 ]\times \mathbb{R}^{d})$, the set of continuous and
suitably differentiable functions on $[0,1 ]\times \mathbb{R}^{d}$. This proves Part (1).

\medskip

\noindent Result (c) Part (2) follows directly from the dynamic programming
principle (see Yong and Zhou \cite[Theorem 3.3, Ch. 4]{yongzhou1999})\textbf{.}

\medskip

\noindent For Result (c) Part (3), by Results Parts (1) and (2), and Ito's formula, we have:
\small
\begin{align*}
& \sum\limits_{m=1}^{n}\sup\limits_{x\in
\mathbb{R}^d}\left\vert H_{\frac{m-1}{n}}^{\epsilon}(x)-H_{\frac{m}{n}}^{\epsilon}\left(x\right) -\frac{1}{n}\sup_{p\in\Gamma}\left[ D_x H_{\frac{m}{n}}^{\epsilon}\left(
x\right)\cdot p\right]
\right\vert \\
=& \sum_{m=1}^{n}\sup_{x\in \mathbb{R}^d}\left\vert \sup_{a\in\mathcal{A}(\frac{m-1}{n},\frac{m}{n})}\!\!E\left[ H_{\frac{m}{n}}^{\epsilon}\left( x+\int_{\frac{m-1}{n}}^{\frac{m}{n}}a_{s}ds+\int_{\frac{m-1%
}{n}}^{\frac{m}{n}}\epsilon\boldsymbol{I}_ddB_{s}\right) \right]  -H_{\frac{m}{n}}^{\epsilon}\left(x\right) -\frac{1}{n}\sup_{p\in\Gamma}\left[ D_x H_{\frac{m}{n}}^{\epsilon}\left(
x\right)\cdot p\right] \right\vert\\
=& \sum_{m=1}^{n}\sup_{x\in \mathbb{R}^d}\left\vert \sup_{a \in\mathcal{A}(\frac{m-1}{n},\frac{m}{n})}\!\!E\left[ \int_{\frac{m-1}{n}}^{\frac{m}{n}}D_x H_{\frac{m}{n}}^{\epsilon}\left(x+\int_{\frac{m-1}{n%
}}^{s}a_{s}ds+\int_{\frac{m-1}{n}}^{s}\epsilon\boldsymbol{I}_ddB_{s}\right)\cdot 
a_{s}ds\right. \right. \\
& \hspace{2cm}\left. +\frac{\epsilon^2}{2}\int_{\frac{m-1}{n}}^{\frac{m}{n}} \textup{tr}\left[D_x^2 H_{\frac{m}{n}}^{\epsilon}\left(  x+\int_{\frac{m-1}{n}}^{s}a_{s}ds+\int_{\frac{m-1}{n}}^{s}\epsilon\boldsymbol{I}_ddB_{s}\right)\right]ds%
\right]-\frac{1}{n}\sup_{p\in\Gamma}\left[ D_x H_{\frac{m}{n}}^{\epsilon}\left(
x\right)\cdot p\right]  \Bigg\vert\\
\leq & \frac{C}{n}\sum_{m=1}^{n}\left\vert
\sup_{a\in \mathcal{A}(\frac{m-1}{n},\frac{m}{n})}\!\!E%
\left[ \sup_{s\in \lbrack \frac{m-1}{n},\frac{m}{n}]}\left( \left\| \int_{%
\frac{m-1}{n}}^{s}a_{s}ds\right\| +\left\| \int_{\frac{m-1}{n}%
}^{s}\epsilon\boldsymbol{I}_ddB_{s}\right\| \right)  \right] \right\vert +C\epsilon^2\\
\leq& C\left(\frac{\sum_{i=1}^d|\overline{\mu}_i-\underline{\mu}_i|}{2n}+2\epsilon\sqrt{\frac{d}{n}}+\epsilon^2\right),
\end{align*}%
\normalsize
where $C$ is a constant that depends only on  the uniform bound of $D_x^2\,H_{t}^{\epsilon}$.
\hfill $\blacksquare $

\begin{Lem}\label{lemma-taylor-0} For any $1\le m\le n$, define a family of functions $\{L_{m,n}\}_{m=1}^n$ by
\begin{equation}\label{function-L}
   L_{m,n}\left(x\right):= H_{\frac{m}{n}}^{\epsilon}\left(x\right)+\frac{1}{n}\sup_{p\in\Gamma}\left[ D_x  H_{\frac{m}{n}}^{\epsilon}\left(
x\right)\cdot p\right]. 
\end{equation}
Let $S_{n}^{\bmtheta}=S_{n-1}^{\bmtheta} +\Z_{n}^{\bmtheta}$ with $\Z_{n}^{\bmtheta}$ defined in \eqref{eq-1}. Then we have:
\begin{align}
\lim_{n\to\infty}\sum_{m=1}^n\left\{\sup_{\bmtheta\in\Theta}E\left[ H_{\frac{m}{n}}^{\epsilon}\left(\frac{S_{m}^{\bmtheta}}{n}\right)\right]-\sup_{\bmtheta\in\Theta}E\left[L_{m,n}\left(\frac{S_{m-1}^{\bmtheta}}{n}\right)\right]\right\}=0.
\end{align}
\end{Lem}

\noindent \textbf{Proof of Lemma \ref{lemma-taylor-0}:} It is sufficient to prove the following two results:\small
\begin{align}
\lim_{n\to\infty}&\sum_{m=1}^n \sup_{\bmtheta\in\Theta}E\left[ \left|  H_{\frac{m}{n}}^{\epsilon}\left(\frac{S_{m}^{\bmtheta}}{n}\right)-H_{\frac{m}{n}}^{\epsilon}\left(\frac{S_{m-1}^{\bmtheta}}{n}\right)-D_x\,H_{\frac{m}{n}}^{\epsilon}\left(\frac{S_{m-1}^{\bmtheta}}{n}\right)\cdot\frac{\Z_m^{\bmtheta}}{n} \right|\right]=0,\label{prove-thm21-1}\\
& \sup_{\bmtheta\in\Theta}E\left[ H_{\frac{m}{n}}^{\epsilon}\left(\frac{S_{m-1}^{\bmtheta}}{n}\right)+D_x\,H_{\frac{m}{n}}^{\epsilon}\left(\frac{S_{m-1}^{\bmtheta}}{n}\right)\cdot\frac{\Z_m^{\bmtheta}}{n}\right] =\sup_{\bmtheta\in\Theta}E\left[L_{m,n}\left(\frac{S_{m-1}^{\bmtheta}}{n}\right)\right].\label{prove-thm21-2}
\end{align}\normalsize

For Result \eqref{prove-thm21-1}, by  Result (1) of Lemma~\ref{lemma-ddp-0},  there exists a constant $C>0$ such that
$$\sup\limits_{t\in[0,1]}\sup\limits_{x\in\mathbb{R}^d} \sum_{j,k=1}^d|\partial_{x_jx_k}^2H_{t}^{\epsilon}(x)|\leq C.$$
It follows from Taylor's expansion that  for any $x,y\in \mathbb{R}^d$, and $t\in[0,1]$,
\begin{equation}\label{le0-thm21-1}
\left|H_{t}^{\epsilon}(x+y)-H_{t}^{\epsilon}(x)- D_x\,H_{t}^{\epsilon}(x)\cdot y\right| \leq C\|y\|^2.
\end{equation}
For any   $1\le m \le n,$  taking $x=\frac{S_{m-1}^{\bmtheta}}{n},y=\frac{\Z_m^{\bmtheta}}{n}$ in  (\ref{le0-thm21-1}), we obtain
\begin{align*}
&\sum_{m=1}^n \sup_{\bmtheta\in\Theta}E\left[ \left| H_{\frac{m}{n}}^{\epsilon}\left(\frac{S_{m}^{\bmtheta}}{n}\right)-H_{\frac{m}{n}}^{\epsilon}\left(\frac{S_{m-1}^{\bmtheta}}{n}\right)-D_x\,H_{\frac{m}{n}}^{\epsilon}\left(\frac{S_{m-1}^{\bmtheta}}{n}\right)\cdot\frac{\Z_m^{\bmtheta}}{n}  \right|\right]\\
\leq&\frac{C}{n^2}\sum_{m=1}^n\sup_{\bmtheta\in\Theta}E[\|\Z_m^{\bmtheta}\|^2]\to0, \text{ as }n\to\infty.
\end{align*}
Next, for Result \eqref{prove-thm21-2}, we have for any $1\le m\le n$, 
\small
\begin{align*}
&\sup_{\bmtheta\in\Theta}E\left[ H_{\frac{m}{n}}^{\epsilon}\left(\frac{S_{m-1}^{\bmtheta}}{n}\right)+D_x\,H_{\frac{m}{n}}^{\epsilon}\left(\frac{S_{m-1}^{\bmtheta}}{n}\right)\cdot\frac{\Z_m^{\bmtheta}}{n}\right]\\
=&\sup_{\bmtheta\in\Theta}E\left[ H_{\frac{m}{n}}^{\epsilon}\left(\frac{S_{m-1}^{\bmtheta}}{n}\right)+\esssup_{\bmtheta\in\Theta}E\left[D_x\,H_{\frac{m}{n}}^{\epsilon}\left(\frac{S_{m-1}^{\bmtheta}}{n}\right)\cdot\frac{\Z_m^{\bmtheta}}{n}\bigg|\mathcal{H}_{m-1}^{\bmtheta}\right]\right]\\
=& \sup_{\bmtheta\in\Theta}E\left[  H_{\frac{m}{n}}^{\epsilon}\left(\frac{S_{m-1}^{\bmtheta}}{n}\right)+\esssup_{\bmtheta\in\Theta}\frac{1}{n}\sum_{j=1}^d\Bigg[\partial_{x_j}H_{\frac{m}{n}}^{\epsilon}\left(\frac{S_{m-1}^{\bmtheta}}{n}\right)E\left[\vartheta_{m,j}X_{m,j}+(1-\vartheta_{m,j})Y_{m,j}|\mathcal{H}_{m-1}^{\bmtheta}\right]\Bigg]\right]\\
=&\sup_{\bmtheta\in\Theta}E\left[  H_{\frac{m}{n}}^{\epsilon}\left(\frac{S_{m-1}^{\bmtheta}}{n}\right)+\sum_{j=1}^d\Bigg[\frac{\overline{\mu}_j}{n}\partial_{x_j}H_{\frac{m}{n}}^{\epsilon}\left(\frac{S_{m-1}^{\bmtheta}}{n}\right)I\left\{\partial_{x_j}H_{\frac{m}{n}}^{\epsilon}\left(\frac{S_{m-1}^{\bmtheta}}{n}\right)\ge0\right\}\right.\\
&\hspace{4.5cm}\left.+\frac{\underline{\mu}_j}{n}\partial_{x_j}H_{\frac{m}{n}}^{\epsilon}\left(\frac{S_{m-1}^{\bmtheta}}{n}\right)I\left\{\partial_{x_j}H_{\frac{m}{n}}^{\epsilon}\left(\frac{S_{m-1}^{\bmtheta}}{n}\right)<0\right\}\Bigg]\right]\\
=& \sup_{\bmtheta\in\Theta}E\left[  L_{m,n}\left(\frac{S_{m-1}^{\bmtheta}}{n}\right)\right].
\end{align*}\normalsize
Thus we obtain Result (\ref{prove-thm21-2}).  \hfill $\blacksquare $

\noindent \textbf{Proof of Theorem \ref{thm-lln-01}:} Let $\hat{f}$ be a bounded continuous extension of $f$ defined on $\mathbb{R}^d$ (we can assume that $\hat{f}$ is constant outside a compact set $\Gamma'\supset\Gamma_{\delta}$). For any $\varepsilon>0$, there exists $\psi\in C_{b,Lip}(\mathbb{R}^d)$ such that $\sup_{x\in\mathbb{R}^d}|\hat{f}(x)-\psi(x)|\leq \varepsilon$.

For any $\epsilon>0$, let $\{H_t^{\epsilon}(x)\}_{t\in[0,1]}$ be the functions defined by \eqref{function-H} via $\psi$. 
It can be checked that
\begin{align*}
& \sup_{\bmtheta\in\Theta}E\left[H_{1}^{\epsilon}\left(\frac{S_n^{\bmtheta}}{n}\right)\right]-E\left[H_{0}^{\epsilon}\left(\frac{\eta}{n}\right)\right]\\
=&\sum_{m=1}^n\left\{ \sup_{\bmtheta\in\Theta}E\left[H_{\frac{m}{n}}^{\epsilon}\left(\frac{S_{m}^{\bmtheta}}{n}\right)\right]- \sup_{\bmtheta\in\Theta}E\left[H_{\frac{m-1}{n}}^{\epsilon}\left(\frac{S_{m-1}^{\bmtheta}}{n}\right)\right]\right\}\\
=&\sum_{m=1}^n\left\{ \sup_{\bmtheta\in\Theta}E\left[H_{\frac{m}{n}}^{\epsilon}\left(\frac{S_{m}^{\bmtheta}}{n}\right)\right]- \sup_{\bmtheta\in\Theta}E\left[L_{m,n}\left(\frac{S_{m-1}^{\bmtheta}}{n}\right)\right]\right\}\\
&+\sum_{m=1}^n\left\{ \sup_{\bmtheta\in\Theta}E\left[L_{m,n}\left(\frac{S_{m-1}^{\bmtheta}}{n}\right)\right]- \sup_{\bmtheta\in\Theta}E\left[H_{\frac{m-1}{n}}^{\epsilon}\left(\frac{S_{m-1}^{\bmtheta}}{n}\right)\right]\right\}\\
=:&\Delta_{n,1}+\Delta_{n,2}.
\end{align*}
By Lemma \ref{lemma-taylor-0} we have $|\Delta_{n,1} |\to0$ as $n\to\infty$. Applying Lemma \ref{lemma-ddp-0}(3) we have $\limsup_{n\to\infty}|\Delta_{n,2} |\le C\epsilon^2$, where $C$ is a constant depends only on the uniform bound of $D_x^2\,H_{t}^{\epsilon}$. Therefore, we have
\begin{align}\label{sym-prove-thm21-1}
\limsup_{n\to\infty}\left|\sup_{\bmtheta\in\Theta}E\left[H_{1}^{\epsilon}\left(\frac{S_n^{\bmtheta}}{n}\right)\right]-E\left[H_{0}^{\epsilon}\left(\frac{\eta}{n}\right)\right]\right|\le C\epsilon^2.
\end{align}
By Result (a) we have:
\begin{align*}
&\left|\max_{x\in\Gamma}\psi(x)-
E\left[H_{0}^{\epsilon}\left(\frac{\eta}{n}\right)\right]\right|\\
\leq&\left|\max_{x\in\Gamma}\psi(x)-
H_{0}^{\epsilon}\left(0\right)\right|+\left|H_{0}^{\epsilon}\left(0\right)-
E\left[H_{0}^{\epsilon}\left(\frac{\eta}{n}\right)\right]\right|\\
=&\left|\sup_{a\in\mathcal{A}(0,1)}E\left[\psi\left(\int_0^{1}a_sds\right)\right]-
\sup_{a\in\mathcal{A}(0,1+\epsilon)}E\left[\psi\left(\int_0^{1+\epsilon}a_sds+\epsilon \int_{0}^{1+\epsilon}\boldsymbol{I}_ddB_{s}\right)\right]\right|+\left|H_{0}^{\epsilon}\left(0\right)-
E\left[H_{0}^{\epsilon}\left(\frac{\eta}{n}\right)\right]\right|\\
\le& L_{\psi}\epsilon(C_{\Gamma}+\sqrt{1+\epsilon})+\frac{C_1E[\|\eta\|]}{n},
\end{align*}
where $L_{\psi}$ is the Lipschitz constant of $\psi$, $C_1$ is a finite constant depending on the uniform bound of $D_x H_{t}^\epsilon$ and $C_{\Gamma}$ is a finite constant depending on the boundary of $\Gamma$. Combine with (\ref{sym-prove-thm21-1}), we have
\begin{align*}
&\limsup_{n\to\infty}\left|\sup_{\bmtheta\in\Theta}E\left[f\left(\frac{S_n^{\bmtheta}}{n}\right)\right]-\max_{x\in\Gamma}f(x)\right|\\
\le&\limsup_{n\to\infty}\left|\sup_{\bmtheta\in\Theta}E\left[\hat{f}\left(\frac{S_n^{\bmtheta}}{n}\right)\right]-\sup_{\bmtheta\in\Theta}E\left[\psi\left(\frac{S_n^{\bmtheta}}{n}\right)\right]\right|\\
&+\limsup_{n\to\infty}\left|\sup_{\bmtheta\in\Theta}E\left[\psi\left(\frac{S_n^{\bmtheta}}{n}\right)\right]-\sup_{\bmtheta\in\Theta}E\left[H_{1}^{\epsilon}\left(\frac{S_n^{\bmtheta}}{n}\right)\right]\right|\\
&+\limsup_{n\to\infty}\left|\sup_{\bmtheta\in\Theta}E\left[H_{1}^{\epsilon}\left(\frac{S_n^{\bmtheta}}{n}\right)\right]-E\left[H_{0}^{\epsilon}\left(\frac{\eta}{n}\right)\right]\right|\\
&+\limsup_{n\to\infty}\left|E\left[H_{0}^{\epsilon}\left(\frac{\eta}{n}\right)\right]-\max_{x\in\Gamma}\psi(x)\right|\\
&+\left|\max_{x\in\Gamma}\psi(x)-\max_{x\in\Gamma}\hat{f}(x)\right|+\left|\max_{x\in\Gamma}\hat{f}(x)-\max_{x\in\Gamma}f(x)\right|\\
\le &2\varepsilon+L_{\psi}\epsilon(2 C_{\Gamma}+\sqrt{\epsilon}+\sqrt{1+\epsilon})+C\epsilon^2.
\end{align*}
Due to the arbitrariness of $\epsilon$ and $\varepsilon$, we have \eqref{lln-1}.

Now we prove \eqref{lln-non-prob}. Recall that $S_{n}^{\bmtheta}=(S_{n,1}^{\bmtheta},\cdots,S_{n,d}^{\bmtheta})$, and $\Gamma=[\underline{\mu}_1,\overline{\mu}_1]\times\cdots\times[\underline{\mu}_d,\overline{\mu}_d]$. By Theorem 3.1 of Chen \cite{chen2016}, we have the following results (A) and (B):
    \begin{align*}
(A)~~~\inf_{\bmtheta\in\Theta}P\left(\underline{\mu}_i\leq \liminf_{n\to\infty}\frac{S_{n,i}^{\bmtheta}}{n}\leq \limsup_{n\to\infty}\frac{S_{n,i}^{\bmtheta}}{n}\leq\overline{\mu}_i,i=1,\cdots,d\right)=1
    \end{align*}
    
(B) Let $\mathcal{C}\left(\frac{S_{n}^{\bmtheta}}{n}\right)$ denote the cluster set of $\frac{S_{n}^{\bmtheta}}{n}$. For any $\mu \in \Gamma$, we have:
    \begin{align*}
    \sup_{\bmtheta\in\Theta}P\left(\mu \in \mathcal{C}\left(\frac{S_{n}^{\bmtheta}}{n}\right)\right)=1.
    \end{align*}

By result (A) and the continuity of $f$, we have
\begin{align}\label{2}
\inf_{\bmtheta\in\Theta}P\left(\limsup_{n\to\infty}f\left(\frac{S_{n}^{\bmtheta}}{n}\right)\leq\max_{x\in\Gamma}f(x)\right)=1 .
\end{align}

On the other hand, there exists $x^*\in\Gamma$ such that $f(x^*)=\max_{x\in\Gamma}f(x)$. Combined with result (B), we have:
\begin{align}
    \sup_{\bmtheta\in\Theta}P\left(f(x^*)\in \mathcal{C}\left(f\left(\frac{S_{n}^{\bmtheta}}{n}\right)\right)\right)=1
\end{align}
then 
\begin{align}\label{2b}
    \sup_{\bmtheta\in\Theta}P\left(f(x^*)\leq \limsup_{n\to\infty}f\left(\frac{S_{n}^{\bmtheta}}{n}\right)\right)=1
\end{align}
Relations \eqref{2} and \eqref{2b} together imply result \eqref{lln-non-prob}.
\hfill $\blacksquare $

\medskip

\noindent{\bf Proof of Theorem \ref{thm-convergence2-sup}:} 
For any $\bmtheta\in\Theta$, set $\Tilde{S}_{n}^{\bmtheta}=\sum_{i=1}^nE\left[\Z_{i}^{\bmtheta}|\mathcal{H}_{i-1}^{\bmtheta}\right]$, then we have $\frac{\Tilde{S}_{n}^{\bmtheta}}{n}\in \Gamma$, $P$-a.s., which implies that
$$
\sup_{\bmtheta\in\Theta}E\left[f\left(\frac{S_{n}^{\bmtheta}}{n}\right)\right]-\max_{x\in\Gamma}f(x)\leq \sup_{\bmtheta\in\Theta}E\left[f\left(\frac{S_{n}^{\bmtheta}}{n}\right)-f\left(\frac{\tilde{S}_{n}^{\bmtheta}}{n}\right)\right]\leq \frac{L_{f}}{n}\left(\sup_{\bmtheta\in\Theta}E\left[\|S_{n}^{\bmtheta}-\tilde{S}_{n}^{\bmtheta}\|^2\right]\right)^{\frac{1}{2}}.
$$
It can be checked that, for any $\bmtheta\in\Theta$
\begin{align*}
 E\left[\|S_{n}^{\bmtheta}-\tilde{S}_{n}^{\bmtheta}\|^2\right]
 \leq&E[\|\eta\|^2]+\sum_{j=1}^{d}\sum_{i=1}^n  E\left[(Z_{i,j}^{\bmtheta}-E[Z_{i,j}^{\bmtheta}|\mathcal{H}_{i-1}^{\bmtheta}])^2\right] \\
 \leq&E[\|\eta\|^2]+\sum_{j=1}^{d}\sum_{i=1}^n\max\{Var(X_j),Var({Y_j})\}\\
 = &nd\sigma^2+E[\|\eta\|^2].
\end{align*}

Therefore, we have
$$
\sup_{\bmtheta\in\Theta}E\left[f\left(\frac{S_{n}^{\bmtheta}}{n}\right)\right]-\max_{x\in\Gamma}f(x)\leq L_{f}\left(\frac{\sigma\sqrt{d}}{\sqrt{n}}+\frac{\sqrt{E[\|\eta\|^2]}}{n}\right).
$$

On the other hand, there exists $x^*=(x_1^*,\cdots,x_d^*)\in\Gamma$ such that $f(x^*)=\max_{x\in\Gamma}f(x)$. Then there exists a strategy $\hat{\bmtheta}$ such that
$$
P\left(\hat{\vartheta}_{ij}=1\right)=\frac{x^*_j-\underline{\mu}_j}{\overline{\mu}_j-\underline{\mu}_j},\quad P\left(\hat{\vartheta}_{ij}=0\right)=\frac{\overline{\mu}_j-x^*_j}{\overline{\mu}_j-\underline{\mu}_j}.
$$
Therefore, 
$$E[Z_{ij}^{\hat{\bmtheta}}]=E\left[E[Z_{ij}^{\hat{\bmtheta}}|\mathcal{H}_{i-1}^{\hat{\bmtheta}}]\right]=E\left[\overline{\mu}_j I\left\{\hat{\vartheta}_{i,j}=1\right\}+\underline{\mu}_j I\left\{\hat{\vartheta}_{i,j}=0\right\}\right]=x_j^*
$$
which leads to
\begin{align*} \sup_{\bmtheta\in\Theta}E\left[f\left(\frac{S_{n}^{\bmtheta}}{n}\right)\right]-\max_{x\in\Gamma}f(x)\geq E\left[f\left(\frac{S_{n}^{\hat{\bmtheta}}}{n}\right)\right]-f(x^*)\geq-\frac{L_f}{n}\left(E[\|S_n^{\hat{\bmtheta}}-nx^*\|^2]\right)^{\frac{1}{2}},
\end{align*}
and
\begin{align*}
 E\left[\|S_{n}^{\hat{\bmtheta}}-nx^*\|^2\right]
 \leq&E[\|\eta\|^2]+\sum_{j=1}^{d}\sum_{i=1}^n  E\left[(Z_{i,j}^{\hat{\bmtheta}}-E[Z_{i,j}^{\hat{\bmtheta}}])^2\right] \\
 \leq&E[\|\eta\|^2]+\sum_{j=1}^{d}\sum_{i=1}^n\left\{\max\{Var(X_j),Var({Y_j})\}+E\left[(E[Z_{i,j}^{\hat{\bmtheta}}|\mathcal{H}_{i-1}^{\hat{\bmtheta}}]-E[Z_{ij}^{\hat{\bmtheta}}])^2\right]\right\}\\
 = &nd\sigma^2+n\sum_{j=1}^d\frac{(\overline{\mu}_j-\underline{\mu}_j)^2}{4}+E[\|\eta\|^2].
\end{align*}

Finally, we have
\begin{align*} \left|\sup_{\bmtheta\in\Theta}E\left[f\left(\frac{S_{n}^{\bmtheta}}{n}\right)\right]-\max_{x\in\Gamma}f(x)\right|\leq L_f\left(\frac{\sigma\sqrt{d}+\sum_{j=1}^d\frac{|\overline{\mu}_j-\underline{\mu}_j|}{2}}{\sqrt{n}}+\frac{\sqrt{E[\|\eta\|^2]}}{n}\right).
\end{align*}

\hfill $\blacksquare $
\medskip

Next, we give the proof of Theorem \ref{thm-pde}. 
For any $\epsilon>0$ and $(t,x)\in[0,1+\epsilon]\times\mathbb{R}^d$, define $u_{\epsilon}(t,x)=H_{t}^{\epsilon}(x)$, where $H_{t}^{\epsilon}$ is defined by \eqref{function-H} with $\psi$ replaced by $\hat{f}$, which is a bounded Lipschitz continuous (with Lipschitz constant $L_{f}$) extension of $f$ from $\Gamma_{\delta}$ to $\mathbb{R}^d$. Then $u_{\epsilon}(t,x)=H_{t}^{\epsilon}(x)$ is the unique solution of the HJB-equation [\ref{pde-1}] and satisfies Lemma \ref{lemma-ddp-0}.

\begin{Lem}\label{lemma-taylor} For any $1\le m\le n$, let functions $\{L_{m,n}\}_{m=1}^n$  be defined by \eqref{function-L}, which can be equivalently denoted by 
$$
L_{m,n}\left(x\right)=u_{\epsilon}\left(\tfrac{m}{n},x\right)+\tfrac{1}{n}\sup_{p\in\Gamma}\left[ D_x u_{\epsilon}\left(\tfrac{m}{n},
x\right)\cdot p\right].
$$
Let the strategies $\{\hat{\bmtheta}^{n,\epsilon}\}_{n\geq1}$ be defined by \eqref{strategy-e-optimal}. Then we have,
\begin{align}
\sum_{m=1}^n\left\{E\left[u_{\epsilon}\left(\frac{m}{n},\frac{S_{m}^{\hat{\bmtheta}^{n,\epsilon}}}{n}\right)\right]-E\left[L_{m,n}\left(\frac{S_{m-1}^{\hat{\bmtheta}^{n,\epsilon}}}{n}\right)\right]\right\}\leq\frac{C}{n}\sum_{i=1}^d \left(\sigma^2+\frac{(\overline{\mu}_i-\underline{\mu}_i)^2}{4}\right), 
\end{align}
where $C$ is a constant only depending on  the uniform bound of $D_{x}^2\, u_{\epsilon}(t,\cdot)$.
\end{Lem}
\noindent \textbf{Proof:} It is sufficient to prove
\begin{align}
&\sum_{m=1}^n E\left[ \left| u_{\epsilon}\left(\frac{m}{n},\frac{S_{m}^{\hat{\bmtheta}^{n,\epsilon}}}{n}\right)-u_{\epsilon}\left(\frac{m}{n},\frac{S_{m-1}^{\hat{\bmtheta}^{n,\epsilon}}}{n}\right)-D_x\,u_{\epsilon}\left(\frac{m}{n},\frac{S_{m-1}^{\hat{\bmtheta}^{n,\epsilon}}}{n}\right)\cdot\frac{\Z_m^{\hat{\bmtheta}^{n,\epsilon}}}{n} \right|\right]\nonumber\\
\leq &\frac{C}{n}\sum_{i=1}^d \left(\sigma^2+\frac{(\overline{\mu}_i-\underline{\mu}_i)^2}{4}\right),\label{prove-1}\\
\text{ and }& \sum_{m=1}^n E\left[ u_{\epsilon}\left(\frac{m}{n},\frac{S_{m-1}^{\hat{\bmtheta}^{n,\epsilon}}}{n}\right)+D_x\,u_{\epsilon}\left(\frac{m}{n},\frac{S_{m-1}^{\hat{\bmtheta}^{n,\epsilon}}}{n}\right)\cdot\frac{\Z_m^{\hat{\bmtheta}^{n,\epsilon}}}{n} -L_{m,n}\left(\frac{S_{m-1}^{\hat{\bmtheta}^{n,\epsilon}}}{n}\right)\right]=0.\label{prove-2}
\end{align}

In fact, by Lemma \ref{lemma-ddp-0} Part (1) with $H_{t}^{\epsilon}(x)=u_{\epsilon}(t,x)$, there exists a constant $C>0$ such that
$$\sup\limits_{t\in[0,1]}\sup\limits_{x\in\mathbb{R}^d} \sum_{j,k=1}^d|\partial_{x_jx_k}^2u_{\epsilon}(t,x)|\leq C.$$
It follows from Taylor's expansion that  for any $x,y\in \mathbb{R}^d$, and $t\in[0,1]$,
\begin{equation}\label{le0-1}
\left|u_{\epsilon}(t,x+y)-u_{\epsilon}(t,x)- D_x\,u_{\epsilon}(t,x)\cdot y\right| \leq C\|y\|^2.
\end{equation}
For any   $1\le m \le n,$  taking $x=\frac{S_{m-1}^{\hat{\bmtheta}^{n,\epsilon}}}{n},y=\frac{\Z_m^{\hat{\bmtheta}^{n,\epsilon}}}{n}$ in  \eqref{le0-1}, we obtain \eqref{prove-1} as follows
\begin{align*}
&\sum_{m=1}^n E\left[ \left| u_{\epsilon}\left(\frac{m}{n},\frac{S_{m}^{\hat{\bmtheta}^{n,\epsilon}}}{n}\right)-u_{\epsilon}\left(\frac{m}{n},\frac{S_{m-1}^{\hat{\bmtheta}^{n,\epsilon}}}{n}\right)-D_x\,u_{\epsilon}\left(\frac{m}{n},\frac{S_{m-1}^{\hat{\bmtheta}^{n,\epsilon}}}{n}\right)\cdot\frac{\Z_m^{\hat{\bmtheta}^{n,\epsilon}}}{n}  \right|\right]\\
\leq&\frac{C}{n^2}\sum_{m=1}^nE[\|\Z_m^{\hat{\bmtheta}^{n,\epsilon}}\|^2]\leq \frac{C}{n^2}\sum_{m=1}^n\sum_{i=1}^dE\left[|\xi_i+\overline{\mu}_i|^2\right]\vee E\left[|\xi_i+\underline{\mu}_i|^2\right]\\
\leq&\frac{C}{n}\sum_{i=1}^d\left(\sigma^2+\frac{(\overline{\mu}_i-\underline{\mu}_i)^2}{4}\right).
\end{align*}

On the other hand,  we have
\begin{align*}
&\sum_{m=1}^n E\left[ u_{\epsilon}\left(\frac{m}{n},\frac{S_{m-1}^{\hat{\bmtheta}^{n,\epsilon}}}{n}\right)+D_x\,u_{\epsilon}\left(\frac{m}{n},\frac{S_{m-1}^{\hat{\bmtheta}^{n,\epsilon}}}{n}\right)\cdot\frac{\Z_m^{\hat{\bmtheta}^{n,\epsilon}}}{n}\right]\\
=&\sum_{m=1}^n E\left[  u_{\epsilon}\left(\frac{m}{n},\frac{S_{m-1}^{\hat{\bmtheta}^{n,\epsilon}}}{n}\right)+\sum_{j=1}^d\Bigg[\partial_{x_j}u_{\epsilon}\left(\frac{m}{n},\frac{S_{m-1}^{\hat{\bmtheta}^{n,\epsilon}}}{n}\right)E\left[\frac{X_{m,j}}{n}|\mathcal{H}_{m-1}^{\hat{\bmtheta}^{n,\epsilon}}\right]I\left\{\partial_{x_j}u_{\epsilon}\left(\frac{m}{n},\frac{S_{m-1}^{\hat{\bmtheta}^{n,\epsilon}}}{n}\right)\ge0\right\}\right.\\
&\hspace{4cm}\left.+\partial_{x_j}u_{\epsilon}\left(\frac{m}{n},\frac{S_{m-1}^{\hat{\bmtheta}^{n,\epsilon}}}{n}\right)E\left[\frac{Y_{m,j}}{n}|\mathcal{H}_{m-1}^{\hat{\bmtheta}^{n,\epsilon}}\right]I\left\{\partial_{x_j}u_{\epsilon}\left(\frac{m}{n},\frac{S_{m-1}^{\hat{\bmtheta}^{n,\epsilon}}}{n}\right)<0\right\}\Bigg]\right]\\
=&\sum_{m=1}^n E\left[  u_{\epsilon}\left(\frac{m}{n},\frac{S_{m-1}^{\hat{\bmtheta}^{n,\epsilon}}}{n}\right)+\sum_{j=1}^d\Bigg[\frac{\overline{\mu}_j}{n}\partial_{x_j}u_{\epsilon}\left(\frac{m}{n},\frac{S_{m-1}^{\hat{\bmtheta}^{n,\epsilon}}}{n}\right)I\left\{\partial_{x_j}u_{\epsilon}\left(\frac{m}{n},\frac{S_{m-1}^{\hat{\bmtheta}^{n,\epsilon}}}{n}\right)\ge0\right\}\right.\\
&\hspace{4.5cm}\left.+\frac{\underline{\mu}_j}{n}\partial_{x_j}u_{\epsilon}\left(\frac{m}{n},\frac{S_{m-1}^{\hat{\bmtheta}^{n,\epsilon}}}{n}\right)I\left\{\partial_{x_j}u_{\epsilon}\left(\frac{m}{n},\frac{S_{m-1}^{\hat{\bmtheta}^{n,\epsilon}}}{n}\right)<0\right\}\Bigg]\right]\\
=&\sum_{m=1}^n E\left[  L_{m,n}\left(\frac{S_{m-1}^{\hat{\bmtheta}^{n,\epsilon}}}{n}\right)\right].
\end{align*}
Then we have \eqref{prove-2}. The proof is completed.  \hfill $\blacksquare $

\noindent \textbf{Proof of Theorem \ref{thm-pde}:} 
It can be checked that
\begin{align*}
&E\left[u_{\epsilon}\left(1,\frac{S_n^{\hat{\bmtheta}^{n,\epsilon}}}{n}\right)\right]-E\left[u_{\epsilon}\left(0,\frac{\eta}{n}\right)\right]\\
=&\sum_{m=1}^n\left\{E\left[u_{\epsilon}\left(\frac{m}{n},\frac{S_{m}^{\hat{\bmtheta}^{n,\epsilon}}}{n}\right)\right]-E\left[u_{\epsilon}\left(\frac{m-1}{n},\frac{S_{m-1}^{\hat{\bmtheta}^{n,\epsilon}}}{n}\right)\right]\right\}\\
=&\sum_{m=1}^n\left\{E\left[u_{\epsilon}\left(\frac{m}{n},\frac{S_{m}^{\hat{\bmtheta}^{n,\epsilon}}}{n}\right)\right]-E\left[L_{m,n}\left(\frac{S_{m-1}^{\hat{\bmtheta}^{n,\epsilon}}}{n}\right)\right]\right\}+\sum_{m=1}^n\left\{E\left[L_{m,n}\left(\frac{S_{m-1}^{\hat{\bmtheta}^{n,\epsilon}}}{n}\right)\right]-E\left[u_{\epsilon}\left(\frac{m-1}{n},\frac{S_{m-1}^{\hat{\bmtheta}^{n,\epsilon}}}{n}\right)\right]\right\}\\
=:&\hat{\Delta}_{n}^1+\hat{\Delta}_n^2.
\end{align*}
By Lemma \ref{lemma-taylor} and Lemma \ref{lemma-ddp-0} Part (3) with $H_{t}^{\epsilon}(x)=u_{\epsilon}(t,x)$, we have
\begin{align}\label{sym-prove-1}
\left|E\left[u_{\epsilon}\left(1,\frac{S_n^{\hat{\bmtheta}^{n,\epsilon}}}{n}\right)\right]-E\left[u_{\epsilon}\left(0,\frac{\eta}{n}\right)\right]\right|\le \frac{C}{n}\sum_{i=1}^d \left(\sigma^2+\frac{(\overline{\mu}_i-\underline{\mu}_i)^2}{4}\right)+C\left(\frac{\sum_{i=1}^d|\overline{\mu}_i-\underline{\mu}_i|}{2n}+2\epsilon\sqrt{\frac{d}{n}}+\epsilon^2\right),
\end{align}
where $C$ is a finite constant depending on $f$.

Next, by Lemma \ref{lemma-ddp-0} Result (a) with $\psi$ replaced by $\hat{f}$ or $f$ on $\Gamma$ we have:
$$
\max_{x\in\Gamma}f(x)=\max_{x\in\Gamma}\hat{f}(x)=\sup_{a\in\mathcal{A}(0,1)}E\left[\hat{f}\left(\int_0^{1}a_sds\right)\right].
$$
Then we have,
\begin{align*}
&\left|\max_{x\in\Gamma}f(x)-
E\left[u_{\epsilon}\left(0,\frac{\eta}{n}\right)\right]\right|\\
\leq& \left|\max_{x\in\Gamma}f(x)-
u_{\epsilon}(0,0)\right|+\left|
E\left[u_{\epsilon}\left(0,\frac{\eta}{n}\right)\right]-u_{\epsilon}(0,0)\right|\\
=&\left|\sup_{a\in\mathcal{A}(0,1)}E\left[\hat{f}\left(\int_0^{1}a_sds\right)\right]-
\sup_{a\in\mathcal{A}(0,1+\epsilon)}E\left[\hat{f}\left(\int_0^{1+\epsilon}a_sds+\epsilon \int_{0}^{1+\epsilon}\boldsymbol{I}_ddB_{s}\right)\right]\right|+\left|
E\left[u_{\epsilon}\left(0,\frac{\eta}{n}\right)\right]-u_{\epsilon}(0,0)\right|\\
\le& L_f\epsilon(C_{\Gamma}+\sqrt{1+\epsilon})+\frac{C_1E[\|\eta\|]}{n},
\end{align*}
where $L_f$ is the Lipschitz constant of $\hat{f}$, $C_1$ is a finite constant depending on the uniform bound of $D_x u_{\epsilon}$ and $C_{\Gamma}$ is a finite constant depending on the boundary of $\Gamma$. Combine with \eqref{sym-prove-1}, we have
\begin{align*}
&\left|E\left[f\left(\frac{S_n^{\hat{\bmtheta}^{n,\epsilon}}}{n}\right)\right]-\max_{x\in\Gamma}f(x)\right|\\
\le&\left|E\left[\hat{f}\left(\frac{S_n^{\hat{\bmtheta}^{n,\epsilon}}}{n}\right)\right]-E\left[u_{\epsilon}\left(1,\frac{S_n^{\hat{\bmtheta}^{n,\epsilon}}}{n}\right)\right]\right|+\left|E\left[u_{\epsilon}\left(1,\frac{S_n^{\hat{\bmtheta}^{n,\epsilon}}}{n}\right)\right]-E\left[u_{\epsilon}\left(0,\frac{\eta}{n}\right)\right]\right|
+\left|E\left[u_{\epsilon}\left(0,\frac{\eta}{n}\right)\right]-\max_{x\in\Gamma}f(x)\right|\\
\le &L_f\epsilon(2C_{\Gamma}+\sqrt{\epsilon}+\sqrt{1+\epsilon})+\frac{C_1E[\|\eta\|]}{n} +C\left(\frac{1}{n}\sum_{i=1}^d \left(\sigma^2+\frac{(\overline{\mu}_i-\underline{\mu}_i)^2}{4}\right)+\frac{\sum_{i=1}^d|\overline{\mu}_i-\underline{\mu}_i|}{2n}+2\epsilon\sqrt{\frac{d}{n}}+\epsilon^2\right).
\end{align*}

The proof of Theorem \ref{thm-pde} is completed.
\hfill $\blacksquare $
\medskip

Finally, we present the proofs of Theorems \ref{thm-lln}, \ref{thm-convergence2} and \ref{thm-ms}. The following lemma can be regarded as a weak form of strategic law for large numbers.
\medskip

\begin{Lem}\label{lemma-strategic-wlln2}
Let $\varphi\in C_b^3(\mathbb{R}^d)$, the space of bounded and three-times continuously differentiable functions with bounded derivatives of all orders less than or equal to three, satisfy the assumptions  in Theorem \ref{thm-lln} with $f$ replaced by $\varphi$.
Furthermore, we assume that there exists a constant $L>0$ such that for any $j=1,\cdots,d$, $|X_j|<L$ and $|Y_j|<L$. Let $N>0$ be the smallest $n$ such that
$$
\frac{2L}{n+1}<\min_{1\le i\le d}\min_{1\le j\le c_i^{k_i}}\{b_{i}^{j}-c_i^{j},c_i^j-b_i^{j-1}\}
$$
 
 Consider the strategy $\bmtheta^c=(\bmvartheta^c_{1},\cdots,\bmvartheta_{m}^c,\cdots)$,
 where $\bmvartheta_m^c=(\vartheta^c_{m1},\cdots,\vartheta_{md}^c)$ is defined by
 \begin{equation}\label{proof-new-strategy-op}
\vartheta_{mi}^c=\left\{\begin{array}{rl}
1,\quad &\text{if } \partial_{x_i}\varphi\left(\frac{S_{m-1}^{\bmtheta^c}}{m-1}\right)\ge0;
\\
0,\quad & \text{otherwise},
\end{array}
\right.\quad i=1,\cdots,d.
\end{equation}
where $S_0^{\bmtheta^c}/0$ is set to be $\eta$. 
Then we have
\begin{equation}\label{eq-strategic-wlln-new}
\liminf_{n\to\infty}\left\{E\left[\varphi\left(\frac{S_n^{\bmtheta^c}}{n}\right)\right]-E\left[\sum_{j_1=1}^{K_1}\cdots\sum_{j_d=1}^{K_d}\varphi\left(c_{1}^{j_1},c_{2}^{j_2},\cdots,c_d^{j_d}\right)I\left\{\frac{S_N^{\bmtheta^c}}{N}\in\Gamma(j_1,\cdots,j_d)\right\}\right]\right\}\ge0,
\end{equation}
 where $\Gamma(j_1,\cdots,j_d)=(b_1^{j_1-1},b_1^{j_1}]\times\cdots\times(b_{d}^{j_d-1},b_d^{j_d}].$   
\end{Lem}

\noindent \textbf{Proof of Lemma \ref{lemma-strategic-wlln2}:}  For convenience, we only give the proof for the case $d=2$, it can be proved similarly for the general case.

For any large enough $n>N$, it can be checked that 
\begin{align*}
&E\left[\varphi\left(\frac{S_n^{\bmtheta^c}}{n}\right)\right]-E\left[\sum_{j_1=1}^{K_1}\sum_{j_2=1}^{K_2}\varphi\left(c_{1}^{j_1},c_{2}^{j_2}\right)I\left\{\frac{S_N^{\bmtheta^c}}{N}\in\Gamma(j_1,j_2)\right\}\right]\\
=&\sum_{m=1}^n\sum_{j_1=1}^{K_1}\sum_{j_2=1}^{K_2} E\left[\left[\varphi\left(\frac{S_m^{\bmtheta^c}}{n}+\frac{n-m}{n}\left(c_{1}^{j_1},c_{2}^{j_2}\right)\right)-\varphi\left(\frac{S_{m-1}^{\bmtheta^c}}{n}+\frac{n-m+1}{n}\left(c_{1}^{j_1},c_{2}^{j_2}\right)\right)\right]I\left\{\frac{S_N^{\bmtheta^c}}{N}\in\Gamma(j_1,j_2)\right\}\right] \\
=:&\sum_{m=1}^n\sum_{j_1=1}^{K_1}\sum_{j_2=1}^{K_2} E\left[D_x\,\varphi\left(\frac{S_{m-1}^{\bmtheta^c}}{n}+\frac{n-m+1}{n}\left(c_{1}^{j_1},c_{2}^{j_2}\right)\right)\cdot\frac{\Z_{m}^{\bmtheta^c}-\left(c_{1}^{j_1},c_{2}^{j_2}\right)}{n}I\left\{\frac{S_N^{\bmtheta^c}}{N}\in\Gamma(j_1,j_2)\right\}\right]+R_n,
\end{align*}
where $ D_x \varphi = (\partial_{x_1} \varphi,\partial_{x_2} \varphi) $ denotes the gradient vector of $\varphi$. It follows from Taylor's expansion that  for any $x,y\in \mathbb{R}^d$, 
\begin{equation}\label{le0}
\left|\varphi(x+y)-\varphi(x)- D_x\,\varphi(x)\cdot y\right| \leq C\|y\|^2,
\end{equation}
where $C$ is a constant depend only on the bound of the second derivatives of $\varphi$. 
For any   $1\le m \le n,$  taking $x=\frac{S_{m-1}^{\bmtheta^c}}{n}+\frac{n-m+1}{n}\left(c_{1}^{j_1},c_{2}^{j_2}\right),y=\frac{\Z_m^{\bmtheta^c}-\left(c_{1}^{j_1},c_{2}^{j_2}\right)}{n}$ in  [\ref{le0}], we obtain
\begin{align*}
|R_n|\leq&\frac{C_{\varphi}}{n}E[\|\eta\|]+\frac{C}{n^2}\sum_{m=1}^n\sum_{j_1=1}^{K_1}\sum_{j_2=1}^{K_2}E\left[\|\Z_m^{\bmtheta^c}\|^2+\|\left(c_{1}^{j_1},c_{2}^{j_2}\right)\|^2\right]\to0, \text{ as }n\to\infty,
\end{align*}
where $C_{\varphi}$ is a finite constant depending on the bound of $D_{x}\varphi$.

On the other hand, for any $N< m\le n$, if $\frac{S_{N}^{\bmtheta^c}}{N}\in\Gamma(j_1,j_2)$, then one have $\frac{S_{m}^{\bmtheta^c}}{m}\in\Gamma(j_1,j_2)$.   Set $\lambda_n^m=\frac{m-1}{n}$, then
\begin{align*}
&E\left[D_x\,\varphi\left(\tfrac{S_{m-1}^{\bmtheta^c}}{n}+\tfrac{n-m+1}{n}\left(c_{1}^{j_1},c_{2}^{j_2}\right)\right)\cdot\tfrac{\Z_{m}^{\bmtheta^c}-\left(c_{1}^{j_1},c_{2}^{j_2}\right)}{n}I\left\{\tfrac{S_N^{\bmtheta^c}}{N}\in\Gamma(j_1,j_2)\right\}\right]\\
=&\sum_{i=1}^2E\left[\partial_{x_i}\varphi\left(\tfrac{S_{m-1}^{\bmtheta^c}}{n}+\tfrac{n-m+1}{n}\left(c_{1}^{j_1},c_{2}^{j_2}\right)\right)\tfrac{Z_{mi}^{\bmtheta^c}-c_i^{j_i}}{n}I\left\{\tfrac{S_N^{\bmtheta^c}}{N}\in\Gamma(j_1,j_2)\right\}\right]\\
=&\sum_{i=1}^2E\left[\partial_{x_i}\varphi\left(\lambda_{n}^m\tfrac{S_{m-1}^{\bmtheta^c}}{m-1}+(1-\lambda_{n}^m)\left(c_{1}^{j_1},c_{2}^{j_2}\right)\right)\tfrac{X_{m,i}-c_i^{j_i}}{n}I\left\{\partial_{x_i}\varphi\left(\tfrac{S_{m-1}^{\bmtheta^c}}{m-1}\right)\ge0\right\}I\left\{\tfrac{S_N^{\bmtheta^c}}{N}\in\Gamma(j_1,j_2)\right\}\right]\\
&+\sum_{i=1}^2E\left[\partial_{x_i}\varphi\left(\lambda_{n}^m\tfrac{S_{m-1}^{\bmtheta^c}}{m-1}+(1-\lambda_{n}^m)\left(c_{1}^{j_1},c_{2}^{j_2}\right)\right)\tfrac{Y_{m,i}-c_i^{j_i}}{n}I\left\{\partial_{x_i}\varphi\left(\tfrac{S_{m-1}^{\bmtheta^c}}{m-1}\right)<0\right\}I\left\{\tfrac{S_N^{\bmtheta^c}}{N}\in\Gamma(j_1,j_2)\right\}\right]\\
=&\sum_{i=1}^2E\left[\partial_{x_i}\varphi\left(\lambda_{n}^m\tfrac{S_{m-1}^{\bmtheta^c}}{m-1}+(1-\lambda_{n}^m)\left(c_{1}^{j_1},c_{2}^{j_2}\right)\right)\tfrac{\overline{\mu}_i-c_i^{j_i}}{n}I\left\{\partial_{x_i}\varphi\left(\tfrac{S_{m-1}^{\bmtheta^c}}{m-1}\right)\ge0\right\}I\left\{\tfrac{S_N^{\bmtheta^c}}{N}\in\Gamma(j_1,j_2)\right\}\right]\\
&+\sum_{i=1}^2E\left[\partial_{x_i}\varphi\left(\lambda_{n}^m\tfrac{S_{m-1}^{\bmtheta^c}}{m-1}+(1-\lambda_{n}^m)\left(c_{1}^{j_1},c_{2}^{j_2}\right)\right)\tfrac{\underline{\mu}_i-c_i^{j_i}}{n}I\left\{\partial_{x_i}\varphi\left(\tfrac{S_{m-1}^{\bmtheta^c}}{m-1}\right)<0\right\}I\left\{\tfrac{S_N^{\bmtheta^c}}{N}\in\Gamma(j_1,j_2)\right\}\right].
\end{align*}
Since 
$$
\left\{\partial_{x_i}\varphi\left(x\right)\ge0,x\in\Gamma(j_1,j_2)\right\}=\left\{x=(x_1,x_2)\in\Gamma(j_1,j_2):x_i\le c_i^{j_i}\right\}=:M_i^{j_i}
$$
Then for any $y\in M_i^{j_i}$ and any $\lambda\in(0,1)$, we have $\lambda y+(1-\lambda)\left(c_{1}^{j_1},c_{2}^{j_2}\right)\in M_i^{j_i}$. Therefore we have, for any $j=1,2$,
$$
\partial_{x_i}\varphi\left(\lambda_{n}^m\frac{S_{m-1}^{\bmtheta^c}}{m-1}+(1-\lambda_{n}^m)\left(c_{1}^{j_1},c_{2}^{j_2}\right)\right)I\left\{\partial_{x_i}\varphi\left(\tfrac{S_{m-1}^{\bmtheta^c}}{m-1}\right)\ge0\right\}I\left\{\tfrac{S_N^{\bmtheta^c}}{N}\in\Gamma(j_1,j_2)\right\}\ge0
$$
and similarly,
$$
\partial_{x_i}\varphi\left(\lambda_{n}^m\frac{S_{m-1}^{\bmtheta^c}}{m-1}+(1-\lambda_{n}^m)\left(c_{1}^{j_1},c_{2}^{j_2}\right)\right)I\left\{\partial_{x_i}\varphi\left(\tfrac{S_{m-1}^{\bmtheta^c}}{m-1}\right)<0\right\}I\left\{\tfrac{S_N^{\bmtheta^c}}{N}\in\Gamma(j_1,j_2)\right\}\le0.
$$
And these will lead to, for $N<m\le n$,
$$
E\left[D_x\,\varphi\left(\frac{S_{m-1}^{\bmtheta^c}}{n}+\frac{n-m+1}{n}\left(c_{1}^{j_1},c_{2}^{j_2}\right)\right)\cdot\frac{\Z_{m}^{\bmtheta^c}-\left(c_{1}^{j_1},c_{2}^{j_2}\right)}{n}I\left\{\tfrac{S_N^{\bmtheta^c}}{N}\in\Gamma(j_1,j_2)\right\}\right]\ge0.
$$
That is,
\begin{align*}
&\liminf_{n\to\infty}\left\{E\left[\varphi\left(\tfrac{S_n^{\bmtheta^c}}{n}\right)\right]-E\left[\sum_{j_1=1}^{K_1}\sum_{j_2=1}^{K_2}\varphi\left(c_{1}^{j_1},c_{2}^{j_2}\right)I\left\{\tfrac{S_N^{\bmtheta^c}}{N}\in\Gamma(j_1,j_2)\right\}\right]\right\}\\
\ge&\liminf_{n\to\infty}\sum_{m=1}^n\sum_{j_1=1}^{K_1}\sum_{j_2=1}^{K_2} E\left[D_x\,\varphi\left(\tfrac{S_{m-1}^{\bmtheta^c}}{n}+\tfrac{n-m+1}{n}\left(c_{1}^{j_1},c_{2}^{j_2}\right)\right)\cdot\tfrac{\Z_{m}^{\bmtheta^c}-\left(c_{1}^{j_1},c_{2}^{j_2}\right)}{n}I\left\{\tfrac{S_N^{\bmtheta^c}}{N}\in\Gamma(j_1,j_2)\right\}\right]-\lim_{n\to\infty}|R_n|\\
\ge&\liminf_{n\to\infty}\sum_{m=1}^N \sum_{j_1=1}^{K_1}\sum_{j_2=1}^{K_2}E\left[D_x\,\varphi\left(\tfrac{S_{m-1}^{\bmtheta^c}}{n}+\tfrac{n-m+1}{n}\left(c_{1}^{j_1},c_{2}^{j_2}\right)\right)\cdot\tfrac{\Z_{m}^{\bmtheta^c}-\left(c_{1}^{j_1},c_{2}^{j_2}\right)}{n}I\left\{\tfrac{S_N^{\bmtheta^c}}{N}\in\Gamma(j_1,j_2)\right\}\right]\\
&+\liminf_{n\to\infty}\sum_{m=N+1}^n \sum_{j_1=1}^{K_1}\sum_{j_2=1}^{K_2}E\left[D_x\,\varphi\left(\tfrac{S_{m-1}^{\bmtheta^c}}{n}+\tfrac{n-m+1}{n}\left(c_{1}^{j_1},c_{2}^{j_2}\right)\right)\cdot\tfrac{\Z_{m}^{\bmtheta^c}-\left(c_{1}^{j_1},c_{2}^{j_2}\right)}{n}I\left\{\tfrac{S_N^{\bmtheta^c}}{N}\in\Gamma(j_1,j_2)\right\}\right]-\lim_{n\to\infty}|R_n|\\
\geq&0.
\end{align*}
 The proof is completed.  \hfill $\blacksquare $

\smallskip

\noindent \textbf{Proof of Theorem \ref{thm-lln}:} 
Let $\psi\in C_b^3(\mathbb{R}^d)$ satisfy the assumption in Theorem \ref{thm-lln}. Furthermore suppose that for $x\notin \mathcal{M}_{\varepsilon}$, $\psi(x)\le0$, for $x\in\mathcal{M}_{\varepsilon}$, $0\le\psi(x)\le1$ and $\psi(c_1^{j_1},\cdots,c_{d}^{j_d})=1$. Then the strategy $\bmtheta^*$ given in [\ref{strategy-op}] is same as the strategy $\bmtheta^c$, which is defined by [\ref{proof-new-strategy-op}] via $\psi$, and we have
\begin{align*}
&\lim_{n\to\infty}P\left( \frac{S_n^{\bmtheta^*}}{n}\in\mathcal{M}_{\varepsilon}\right)\\
=&\lim_{n\to\infty}E\left[I{\left\{\frac{S_n^{\bmtheta^*}}{n}\in\mathcal{M}_{\varepsilon}\right\}}\right]\\
\ge&\lim_{n\to\infty}E\left[\psi\left( \frac{S_n^{\bmtheta^*}}{n}\right)\right]-1+1\\
=&\liminf_{n\to\infty}\left\{E\left[\psi\left( \frac{S_n^{\bmtheta^*}}{n}\right)\right]-E\left[\sum_{j_1=1}^{K_1}\cdots\sum_{j_d=1}^{K_d}\psi\left(c_{1}^{j_1},c_{2}^{j_2},\cdots,c_d^{j_d}\right)I\left\{\frac{S_N^{\bmtheta^*}}{N}\in\Gamma(j_1,\cdots,j_d)\right\}\right]\right\}+1\\
\ge&1,
\end{align*}
where the last inequality is due to Lemma \ref{lemma-strategic-wlln2}. Then we have $\lim_{n\to\infty}P\left( \frac{S_n^{\bmtheta^*}}{n}\in\mathcal{M}_{\varepsilon}\right)=1$.
\hfill $\blacksquare $

\smallskip

\noindent \textbf{Proof of Theorem \ref{thm-ms}:}
Suppose $x^*$ is any global maximum point of $f(x)$ on $\Gamma$.  
By the definition of the arms in \eqref{denote-xy}, since $\xi_j$ and $\xi_j'$ are bounded random variables, assume that $\underset{1\le j \le d}{\max} \{|X_j|,|Y_j|\}< L$. 
For strategy $\bmtheta$, we have: 
\begin{align*}
\left\|\frac{S_{n+1}^{\bmtheta}}{n+1}-\frac{S_{n}^{\bmtheta}}{n}\right\|_{\infty}&=\frac{1}{n+1}\left\|S_{n+1}^{\bmtheta}-\frac{(n+1)S_n^{\bmtheta}}{n}\right\|_{\infty}\\
&\le \frac{1}{n+1}\left(\left\|S_{n+1}^{\bmtheta}-S_{n}^{\bmtheta}\right\|_{\infty}+\left\|\frac{S_{n}^{\bmtheta}}{n}\right\|_{\infty}\right)\le \frac{2L}{n+1}.
\end{align*}
This indicates that for each $n\ge m > N=\left\lceil\frac{4L}{\delta_0} \right\rceil$, it follows $\left\|\frac{S_{n+1}^{\bmtheta}}{n+1}-\frac{S_{n}^{\bmtheta}}{n}\right\|_{\infty}<\delta_0/2$, which means that when $S_m^{\bmtheta}/m\in B(x^*,\delta_0)$, $S_n^{\bmtheta'}/n$ will not jump out of $B(x^*,\delta_0)$. 

Since $f(x)$ has only one maximum point $x^*$ in $B(x^*, \delta_0)$, then according to Theorem \ref{thm-lln}, it follows for any $\varepsilon>0$: 
\begin{equation}
\lim_{n\to\infty}P\left(\left\|\frac{S_n^{\bmtheta'}}{n}-x^*\right\|_{\infty}< \varepsilon\right)=1. 
\end{equation}
\hfill $\blacksquare$
\smallskip

\noindent{\bf Proof of Theorem \ref{thm-convergence2}:} 
For any $\bmtheta\in\Theta$, set $\Tilde{S}_{n}^{\bmtheta}=\sum_{i=1}^nE\left[\Z_{i}^{\bmtheta}|\mathcal{H}_{i-1}^{\bmtheta}\right]$, then we have $\frac{\Tilde{S}_{n}^{\bmtheta}}{n}\in \Gamma$, $P$-a.s., which implies that
$$
E\left[f\left(\frac{S_{n}^{\bmtheta^*}}{n}\right)\right]-\max_{x\in\Gamma}f(x)\leq E\left[f\left(\frac{S_{n}^{\bmtheta^*}}{n}\right)-f\left(\frac{\tilde{S}_{n}^{\bmtheta^*}}{n}\right)\right]\leq \frac{L_{f}}{n}\left(E\left[\|S_{n}^{\bmtheta^*}-\tilde{S}_{n}^{\bmtheta^*}\|^2\right]\right)^{\frac{1}{2}}.
$$
It can be checked that,
\begin{align*}
 E\left[\|S_{n}^{\bmtheta^*}-\tilde{S}_{n}^{\bmtheta^*}\|^2\right]
 \leq&\sum_{j=1}^{d}\sum_{i=1}^n  E\left[(Z_{ij}^{\bmtheta^*}-E[Z_{ij}^{\bmtheta^*}|\mathcal{H}_{i-1}^{\bmtheta^*}])^2\right] +E[\|\eta\|^2]\\
 \leq&\sum_{j=1}^{d}\sum_{i=1}^n\max\{Var(X_j),Var({Y_j})\}+E[\|\eta\|^2]\\
 = &nd\sigma^2+E[\|\eta\|^2].
\end{align*}

Therefore, we have
$$
E\left[f\left(\frac{S_{n}^{\bmtheta^*}}{n}\right)\right]-f(x^*)\leq L_{f}\left(\frac{\sigma\sqrt{d}}{\sqrt{n}}+\frac{\sqrt{E[\|\eta\|^2]}}{n}\right).
$$

On the other hand,
\begin{align*}
    &E\left[f\left(\frac{S_n^{\bmtheta^*}}{n}\right)\right]-f(x^*)\\
    =&\sum_{m=1}^nE\left[f\left(\frac{S_m^{\bmtheta^*}}{n}+\frac{(n-m)x^*}{n}\right)-f\left(\frac{S_{m-1}^{\bmtheta^*}}{n}+\frac{(n-m+1)x^*}{n}\right)
    \right]+E\left[f\left(\frac{\eta}{n}+x^*\right)-f\left(x^*\right)
    \right]\\
    \geq &\sum_{m=1}^nE\left[\frac{\Z_{m}^{\bmtheta^*}-x^*}{n}\cdot D_xf\left(\frac{S_{m-1}^{\bmtheta^*}}{n}+\frac{(n-m+1)x^*}{n}\right)
    \right]-\frac{K_f}{n^2}\sum_{m=1}^nE\left[\left\|\Z_{m}^{\bmtheta^*}-x^*\right\|^2
    \right]-\frac{L_f}{n}E[\|\eta\|]\\
    \geq&-\frac{K_f\left[d\sigma^2+\sum_{j=1}^d(\overline{\mu}_j-\underline{\mu}_j)^2\right]}{n}-\frac{L_f}{n}E[\|\eta\|],
\end{align*}
where $K_f=\max_{i,j}\sup_{x\in\Gamma_{\delta}}|\partial_{x_ix_j}^2f(x)|$. The last inequality is due to, for any $1\le m\le n$, set $\lambda_n^m=\frac{m-1}{n}$, we have
\begin{align*}
&E\left[D_x\,f\left(\frac{S_{m-1}^{\bmtheta^*}}{n}+\frac{n-m+1}{n}x^*\right)\cdot\frac{\Z_{m}^{\bmtheta^*}-x^*}{n}\right]\\
=&\sum_{j=1}^dE\left[\partial_{x_j}f\left(\frac{S_{m-1}^{\bmtheta^*}}{n}+\frac{n-m+1}{n}x^*\right)\frac{Z_{mj}^{\bmtheta^*}-x^*_j}{n}\right]\\
=&\sum_{j=1}^dE\left[\partial_{x_j}f\left(\lambda_{n}^m\frac{S_{m-1}^{\bmtheta^*}}{m-1}+(1-\lambda_{n}^m)x^*\right)\frac{X_{m,j}-x^*_j}{n}I\left\{\partial_{x_j}f\left(\tfrac{S_{m-1}^{\bmtheta^*}}{m-1}\right)\ge0\right\}\right]\\
&+\sum_{j=1}^dE\left[\partial_{x_j}f\left(\lambda_{n}^m\frac{S_{m-1}^{\bmtheta^*}}{m-1}+(1-\lambda_{n}^m)x^*\right)\frac{Y_{m,j}-x^*_j}{n}I\left\{\partial_{x_j}f\left(\tfrac{S_{m-1}^{\bmtheta^*}}{m-1}\right)<0\right\}\right]\\
=&\sum_{j=1}^dE\left[\partial_{x_j}f\left(\lambda_{n}^m\frac{S_{m-1}^{\bmtheta^*}}{m-1}+(1-\lambda_{n}^m)x^*\right)\frac{\overline{\mu}_j-x^*_j}{n}I\left\{\partial_{x_j}f\left(\tfrac{S_{m-1}^{\bmtheta^*}}{m-1}\right)\ge0\right\}\right]\\
&+\sum_{j=1}^dE\left[\partial_{x_j}f\left(\lambda_{n}^m\frac{S_{m-1}^{\bmtheta^*}}{m-1}+(1-\lambda_{n}^m)x^*\right)\frac{\underline{\mu}_j-x^*_j}{n}I\left\{\partial_{x_j}f\left(\tfrac{S_{m-1}^{\bmtheta^*}}{m-1}\right)<0\right\}\right].
\end{align*}
Since 
$$
\left\{\partial_{x_j}f\left(x\right)\ge0\right\}=\left\{x=(x_1,\cdots,x_d)\in\Gamma_{\delta}:x_j\le x^*_j\right\}=:M_j
$$
Then for any $y\in M_j$ and any $\lambda\in(0,1)$, we have $\lambda y+(1-\lambda)x^*\in M_j$. Therefore we have, for any $j=1,\cdots,d$,
$$
\partial_{x_j}f\left(\lambda_{n}^m\frac{S_{m-1}^{\bmtheta^*}}{m-1}+(1-\lambda_{n}^m)x^*\right)I\left\{\partial_{x_j}f\left(\tfrac{S_{m-1}^{\bmtheta^*}}{m-1}\right)\ge0\right\}\ge0
$$
and similarly,
$$
\partial_{x_j}f\left(\lambda_{n}^m\frac{S_{m-1}^{\bmtheta^*}}{m-1}+(1-\lambda_{n}^m)x^*\right)I\left\{\partial_{x_j}f\left(\tfrac{S_{m-1}^{\bmtheta^*}}{m-1}\right)<0\right\}\le0.
$$
And these will lead to
$$
E\left[D_x\,f\left(\frac{S_{m-1}^{\bmtheta^*}}{n}+\frac{n-m+1}{n}x^*\right)\cdot\frac{\Z_{m}^{\bmtheta^*}-x^*}{n}\right]\ge0.
$$
Therefore, we have
$$
\left|E\left[f\left(\frac{S_n^{\bmtheta^*}}{n}\right)\right]-f(x^*)\right|\leq \max\left\{L_f\left(\frac{\sigma\sqrt{d}}{\sqrt{n}}+\frac{\sqrt{E[\|\eta\|^2]}}{n}\right),\frac{K_f\left(d\sigma^2+\sum_{j=1}^d(\overline{\mu}_j-\underline{\mu}_j)^2\right)+L_fE[\|\eta\|]}{n}\right\}
$$

For any $x\in\mathbb{R}^d$, define $\phi(x)=\|x-x^*\|^2$, therefore we have $\phi(x^*)=0$ and by Taylor's formula
\begin{align*}
    &E\left[\left\|\frac{S_n^{\bmtheta^*}}{n}-x^*\right\|^2\right]\\
    =&E\left[\phi\left(\frac{S_n^{\bmtheta^*}}{n}\right)\right]-\phi(x^*)\\
    =&\sum_{m=1}^nE\left[\phi\left(\frac{S_m^{\bmtheta^*}}{n}+\frac{n-m}{n}x^*\right)-\phi\left(\frac{S_{m-1}^{\bmtheta^*}}{n}+\frac{n-m+1}{n}x^*\right)\right]+E\left[\phi\left(\frac{\eta}{n}+x^*\right)-\phi\left(x^*\right)\right]\\
    = &\sum_{m=1}^nE\left[\frac{\Z_{m}^{\bmtheta^*}-x^*}{n}\cdot D_x\phi\left(\frac{S_{m-1}^{\bmtheta^*}}{n}+\frac{(n-m+1)x^*}{n}\right)
    \right]+\frac{1}{n^2}\sum_{m=1}^nE\left[\left\|\Z_{m}^{\bmtheta^*}-x^*\right\|^2\right]+\frac{1}{n^2}E[\|\eta\|^2]
\end{align*}
It can be checked that,
\begin{align*}
&E\left[D_x\,\phi\left(\frac{S_{m-1}^{\bmtheta^*}}{n}+\frac{n-m+1}{n}x^*\right)\cdot\frac{\Z_{m}^{\bmtheta^*}-x^*}{n}\right]\\
=&\sum_{j=1}^dE\left[\partial_{x_j}\phi\left(\frac{S_{m-1}^{\bmtheta^*}}{n}+\frac{n-m+1}{n}x^*\right)\frac{Z_{mj}^{\bmtheta^*}-x^*_j}{n}\right]\\
=&\sum_{j=1}^dE\left[2\left(\frac{\sum_{i=1}^{m-1}Z_{ij}^{\bmtheta^*}}{n}-\frac{m-1}{n}x^*_j\right)\frac{Z_{mj}^{\bmtheta^*}-x^*_j}{n}\right]\\
=&\sum_{j=1}^dE\left[2\left(\frac{\sum_{i=1}^{m-1}Z_{ij}^{\bmtheta^*}}{n}-\frac{m-1}{n}x^*_j\right)\frac{X_{m,j}-x^*_j}{n}I\left\{\partial_{x_j}f\left(\tfrac{S_{m-1}^{\bmtheta^*}}{m-1}\right)\ge0\right\}\right]\\
&+\sum_{j=1}^dE\left[2\left(\frac{\sum_{i=1}^{m-1}Z_{ij}^{\bmtheta^*}}{n}-\frac{m-1}{n}x^*_j\right)\frac{Y_{m,j}-x^*_j}{n}I\left\{\partial_{x_j}f\left(\tfrac{S_{m-1}^{\bmtheta^*}}{m-1}\right)<0\right\}\right]\\
=&\sum_{j=1}^dE\left[\frac{2(m-1)}{n}\left(\frac{\sum_{i=1}^{m-1}Z_{ij}^{\bmtheta^*}}{m-1}-x^*_j\right)\frac{\overline{\mu}_{j}-x^*_j}{n}I\left\{\partial_{x_j}f\left(\tfrac{S_{m-1}^{\bmtheta^*}}{m-1}\right)\ge0\right\}\right]\\
&+\sum_{j=1}^dE\left[\frac{2(m-1)}{n}\left(\frac{\sum_{i=1}^{m-1}Z_{ij}^{\bmtheta^*}}{m-1}-x^*_j\right)\frac{\underline{\mu}_{j}-x^*_j}{n}I\left\{\partial_{x_j}f\left(\tfrac{S_{m-1}^{\bmtheta^*}}{m-1}\right)<0\right\}\right]
\end{align*}
Since $\underline{\mu}_j\leq x_j^*\leq \overline{\mu}_j$ and
$$
\left\{\partial_{x_j}f\left(x\right)\ge0\right\}=\left\{x=(x_1,\cdots,x_d)\in\Gamma_{\delta}:x_j\le x^*_j\right\},
$$
we have
$$
E\left[D_x\,\phi\left(\frac{S_{m-1}^{\bmtheta^*}}{n}+\frac{n-m+1}{n}x^*\right)\cdot\frac{\Z_{m}^{\bmtheta^*}-x^*}{n}\right]\leq0,
$$
and then
\begin{align*}
    E\left[\left\|\frac{S_n^{\bmtheta^*}}{n}-x^*\right\|^2\right]= &\sum_{m=1}^nE\left[\frac{Z_{m}^{\bmtheta^*}-x^*}{n}\cdot D_x\phi\left(\frac{S_{m-1}^{\bmtheta^*}}{n}+\frac{(n-m+1)x^*}{n}\right)
    \right]+\frac{1}{n^2}\sum_{m=1}^nE\left[\left\|Z_{m}^{\bmtheta^*}-x^*\right\|^2\right]\\
    \leq &\frac{d\sigma^2+\sum_{j=1}^d(\overline{\mu}_j-\underline{\mu}_j)^2}{n}+\frac{E[\|\eta\|^2]}{n^2}.
\end{align*}
\hfill $\blacksquare $

\smallskip

\section{Practical Implementation of the SMCO Algorithm}\label{append-B}

In this subsection, we describe in more details how we practically implement Algorithm \ref{alg:smco_algo} in R, and introduce three versions of our algorithms, namely ``SMCO'', ``SMCO-R'', and ``SMCO-BR'', which we use for the numerical experiments in Section \ref{sec:Numerical}. 

For the practical implementation of Algorithm \ref{alg:smco_algo}, which requires the evaluation of the gradient of the function, we compute the finite differences place of the gradients. We adopt this approach for two reasons. First, analytical gradient information is usually hard, or even impossible, to obtain in complex and/or black-box models. Second, we observe that many standard optimization packages in R, even if gradient based, do not require user specification of the gradient, and we follow this practice to: (i) make the comparisons with other optimization algorithms in R relatively fair, and (ii) maximize user friendliness of our algorithm as other algorithms do.

Motivated by the Strategic Law of Large Numbers, we adopt the following step-adaptive approach to calculate the finite differences. Specifically, at iteration $n$, we define the step size in the $j$-th coordinate 
$h_{n,j} := (\overline{\mu}_j - \underline{\mu}_j)/(n + 1)$
and compute the following finite difference $\Delta_{x_j}f(\hat{x}_n) := f(\hat{x}_n + h_{n,j}\cdot e_j) - f(\hat{x}_n - h_{n,j} \cdot e_j),$
where $e_j$ denotes the $j$-th elementary vector. We then use the sign information $\Delta_{x_j}f(\hat{x}_n) \geq 0 $ in place of $\partial_{x_j}f(\hat{x}_n) \geq 0$ in Algorithm \ref{alg:smco_algo}. Observe that $\Delta_{x_j}f(x) \geq 0$ is equivalent to $\Delta_{x_j}f(x)/(2h_{n,j}) \geq 0$, and that $\Delta_{x_j}f(x)/(2h_{n,j}) \to \partial_{x_j}f(x) $ as $n\to\infty$. Hence, the finite difference calculation above is a valid approximation of the gradient, establishing the coherence of our practical implementation with the gradient-based theoretical algorithm.

We now provide a heuristic explanation for our particular choice of step size $h_n$, which is motivated by our Strategic Law of Large Numbers. Specifically, under the Strategic Law of Large Numbers, $\hat{x}_n$ takes the form of a running average, and the update from $\hat{x}_n$ to $\hat{x}_{n+1}$ can be rewritten as $\hat{x}_{n+1} = \frac{n}{n+1}\hat{x}_n + \frac{1}{n+1} \Z^{\theta^*}_{n+1}$. Recall that $\Z^{\theta^*}_{n+1}$ is drawn from two fixed distributions concentrated near the upper/lower bounds $\overline{\mu}_j$ or $\underline{\mu}_j$ in each coordinate, and thus the step size of the update from $\hat{x}_n$ to $\hat{x}_{n+1}$ is of the order $(n+1)^{-1}$. Since the updating of $\hat{x}_n$ is guided by the (sign) information in the finite difference calculation, it is natural to set the step size $h_n$ for the finite difference calculation to be of the same rate $(n+1)^{-1}$ as the step size for the updating of $\hat{x}_n$.

We now describe the default choice of hyperparameters for our SMCO algorithms. For the distributions of $X_j=\overline{\mu}_j +\xi_j$, $ Y_j=\underline{\mu}_{j}+\xi'_j$ as introduced in \eqref{denote-xy}, we set $\xi_j,\xi_j' \sim U[-\d^*,\d^*]\cdot (\ol{\mu}_j-\ul{\mu}_j)$ with default $\delta^*=0.05$, i.e., 5\% of the length of the domain in each coordinate. We set the default tolerance level of our algorithms to be $1e-6$, the maximum number of iterations to be 500. By default, we generate the starting point(s) for our SMCO algorithms $\hat{x}_0 \sim \eta$ with $\eta$ being the uniform distribution over the domain. See Appendix \ref{append-D} and \ref{append-E} for sensitivity checks with respect to these hyperparameter configurations.
 
\section{Details on Comparison-Group Algorithms}\label{append-C}

$\quad$

\noindent\textbf{Algorithm Group I: Six Local Optimizers}
\medskip
\newline
This set consists of six gradient-based
or derivative-free optimization algorithms that are designed
to converge to a local optimum starting from a given starting point.

\begin{itemize}
\item GD: Gradient Descent. We implement this algorithm in R with learning
rate 0.1, momentum 0.9, max number of iterations 1000 and tolerance
1e-6.
\item L-BFGS: As implemented in R command \texttt{optim} through method ``L-BFGS-B'' in R package
\texttt{stats} (version 4.4.2), with default hyperparameter values.
\item SignGD: Sign Gradient Descent by \cite{Bern2018}. We implement this
in R with learning rate 0.1, learning decay rate 0.995, max number
of iterations 1000 and tolerance 1e-6.
\item SPSA: Simultaneous Perturbation Stochastic Approximation by \cite{Spall92}.
The tuning parameters of SPSA  are set following \cite{Spall98} and \cite{GuoFu2022}:
$A=50,a=0.1,\alpha=0.602,\gamma=0.101,$$c=$1e-3, and tolerance 1e-7.
\item ADAM: Adaptive Moment Estimation method by \cite{Kingma2014}, as
implemented in the R package \texttt{optimg} (version 0.1.2), with
default hyperparameter values.
\item BOBYQA: Bound Optimization by Quadratic Approximation by \cite{Powell09}, as implemented
in the R package \texttt{nloptr} (version 2.1.1), with default hyperparameter
values.
\end{itemize}

\noindent\textbf{Algorithm Group II: Six Global Optimizers}
\medskip
\newline
This set consists of six meta-heuristic global optimization algorithms. Since
these algorithms are explicitly designed for global optimization (and many already implement multiple starting points from within), we do not pass any generated starting points to these algorithms whenever applicable.\footnote{A few of these global algorithms nevertheless require a user-designated starting point, in which case we pass a start point randomly chosen from the set of multiple starting points generated for our SMCO and the Group I algorithms.}
\begin{itemize}
\item GenSA: Generalized Simulated Annealing implemented in the R package
\texttt{GenSA} (version 1.1.14.1), with default hyperparameter values.\footnote{Another R implementation for simulated annealing is availble through the \texttt{optim} command via method ``\texttt{SANN}''. However, initial tests of GenSA versus optim-SANN reveals that \texttt{optim-SANN} performs much worse than our algorithms and GenSA, and is thus removed from our comparison.}
\item SA: Simulated Annealing implement in the R package \texttt{SA} (version 1.0-9), with default hyperparameter values.
\item DEoptim: Differential Evolution implemented in the R package \texttt{DEoptim}
(version 2.2-8), with default hyperparameter values.
\item STOGO: Stochastic Global Optimization as implemented in the R package
\texttt{nloptr} (version 2.1.1), with default hyperparameter values.
\item GA: Genetic Algorithms as implemented in the R package \texttt{GA}
(version 3.2.4), with default hyperparameter values.
\item PSO: Particle Swarm Optimizer as implemented in the R package \texttt{PSO}
(version 1.0.4, with default hyperparameter values.\footnote{Another R package for particle swarm optimization, namely \texttt{hydroPSO}, was removed from CRAN and no longer available on R 4.4.2.}
\end{itemize}

%

\begin{remark}[About C++]
 Our algorithms SMCO/SMCO-R/SMCO-BR are fully written in
R without (explicit) use of C++. Translating the key functions in
our algorithms into C++ and calling these functions in R through the
\texttt{Rcpp} package are likely to significantly increase the running speed
of our algorithms. In fact, the documentation of some optimization
algorithms, such as \texttt{nloptr} and \texttt{GenSA}, explicitly
state that their core functions/algorithms are coded in C++ for faster
execution. Hence, one should keep this difference in mind when comparing the running speed between our algorithms and
those algorithms having exploited C++.
\end{remark}

\section{Supplemental Numerical Results on Deterministic Test Functions}\label{append-D}

This section contains additional numerical results on the four deterministic test functions (Rastrigin, Griewank, Ackley, Michaelwicz) considered in the main text, including robustness checks across an array of test configurations and hyperparameter settings.  

\textbf{Tables \ref{tab:Comp_SS_d1}} replicates the results in Table \ref{tab:Comp_SS_d2} for dimension $\bm{d=1}$, further confirming that our SMCO algorithms perform drastically better than all Group I algorithms, when all are run from a single common randomly generated starting point.

\textbf{Tables \ref{tab:Comp2d}-\ref{tab:Comp50d}} report results for the four test functions for {dimensions} $\bm{d = 2}$,  $\bm{10}$, $\bm{20}$, \textbf{and} $\bm{50}$, in which our SMCO algorithms and all Group I algorithms are run with ($10\sqrt{d} \approx 14$, $32$,  $45$, and $71$) uniformly drawn starting points, respectively.
\footnote{These results were generated in a previous using 200 maximum iterations and 1e-8 tolerance levels for our SMCO/-R/-BR algorithms.}. Again, our (multi-start) SMCO algorithms demonstrate robust performance across almost all test configurations, with speed and accuracy levels that compare competitively, and in fact often favorably, to Group II (global) algorithms.  Not very surprisingly, all Group I (local) algorithms, even with the use of multiple starting points, tend to suffer from very large errors in most test configurations, given the challenging landscapes of the four test functions in the moderate dimension (say, $d=10$). This demonstrates that the use of multiple starting points alone cannot explain the good and robust performance of our SMCO algorithms. Among Group II (global) algorithms, GenSA appears to be the only one that shows similar level of robustness in its performance across the configurations. That said, it still yields noticeably larger errors (in terms of RMSE and AE99) than ours, say, in the minimization of Rastrigin and Ackley functions. It is worth noting again that our SMCO algorithms are particularly stable in terms of the 99th error percentile: certain global optimization algorithms may outperform ours in terms of RMSE under a particular specification, but its 99th error percentile becomes notably larger than ours, showing that our algorithm is less sensitive to the randomness inherent in the starting points and the optimization algorithms.


\textbf{Table \ref{tab:Comp_200d_diag}} replicates Table \ref{tab:Comp_200d} on optimizaiton in the high-dimensional setting $d=200$, but uses \textbf{diagonal starting points} rather than uniformly drawn ones. Specifically, under the multiple starting points are generated deterministically as an equally spaced grid of points on the hyper-diagonal of the rectangular parameter space. \textbf{Table \ref{tab:Comp_d500_3glob}} further compares our SMCO algorithms with Group II algorithms in \textbf{500-dimensional} settings.  Observe that our SMCO-R/-BR algorithms perform consistently well across most of the test cases, with the maximization of the Rastrigin function being the only exception where our algorithms yield significant errors. 

\textbf{Table \ref{tab:shift_sensitivity}} analyzes the sensitivity of our results with respect to the magnitude of the asymmetrization  transformation, using the Rastrigin function of $d=500$. In the main text, we randomly push out the support by 25\% ±5\% on one side, and by 50\% ±10\% on the other side. We consider two additional variants: Variant 1 pushes out both sides by 25\% ±5\%, resulting in a smaller and more symmetric domain; Variant 2 pushes out one side by 25\% ±5\% but 75\% ±10\% on the other, resulting in a larger and more asymmetric domain.  Note that larger domains make the optimization problem harder for all optimization algorithms.  
\textbf{Table \ref{tab:tolerance} } analyzes the sensitivity of our results with respect to the tolerance level, using the original Rastrigin function with $d=500$.  \textbf{Table \ref{tab:rastrigin400d_combined}} reports sensitivity to different design of starting values and different max iteration of SMCO algorithms using the original Rastrigin function with $d=400$. \textbf{Table \ref{tab:rastrigin400d_combined1}} reports sensitivity to different bound buffer $\delta^*$ of our SMCO algorithms using original Rastrigin function with $d=400$. Overall, the pattern of the performance comparison of our SMCO algorithms with other algorithms is quite stable. 

\textbf{Table \ref{tab:rastrigin1000_hd_combined}} reports results on the comparison of our SMCO algorithms with three best-performing global optimization algorithms in a very high-dimensional setting ($d=1000)$. In particular, the rows labeled GenSA* and PSO* contain results from GenSA and PSO with their hyperparameter adjusted so that the running time of GenSA* and PSO* are on the same order of magnitude as those of our SMCO algorithms (with parallelization). Since the running time of DEoptim under the default setting is about twice of SMCO algorithms, we view it as being on the same order of magnitude and did not tune DEoptim hyperparameters. The overall conclusion of this table is that our SMCO algorithms are able to deliver better results than GenSA*, PSO* and DEoptim when their running times are similar. 

\begin{figure}
\caption{Visualization of the Cross-Leg Function}
\label{fig:CrossLeg}
    \centering
    \includegraphics[scale = 0.4]{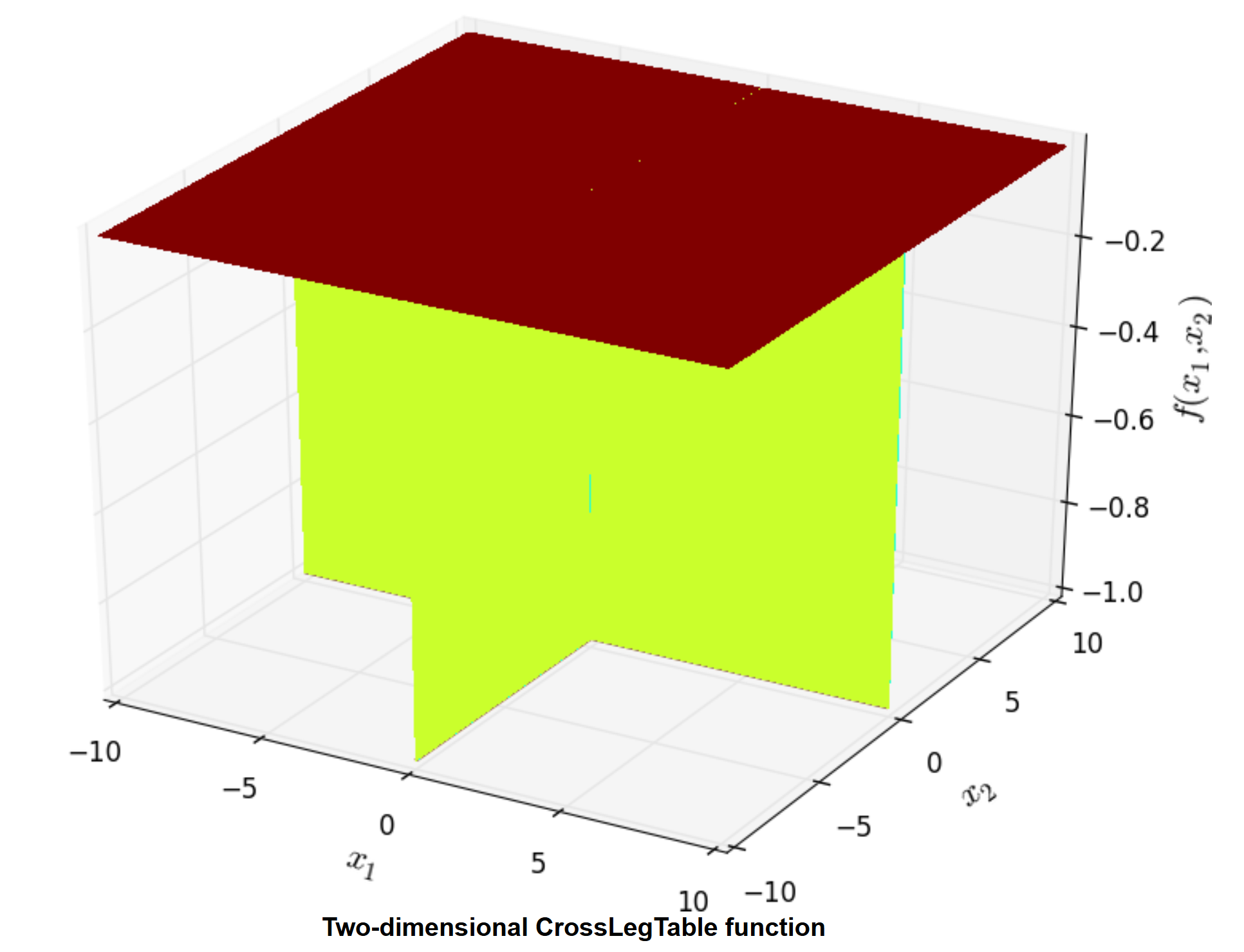}
    
\end{figure}

Finally, in \textbf{Table \ref{tab:CrossLeg}}, we investigate the performance of our SMCO algorithm in minimizing a very challenging test function, the so-called Cross Leg Table function, which is ranked among the top three in terms of global optimization difficulty; see the webpage by Gavana (2009):
``https://infinity77.net/global\_optimization/test\_functions.html\#test-functions-index''. The original Cross-Leg-Table function has domain $[-10,10]^2$ (with $d=2$), with global min value $-1$ attained at $(0,0)$. It features an almost flat landscape with very thin and deep drop near the minimum (visualized in Figure \ref{fig:CrossLeg}). Hence, it is similar to a Dirac delta function at $(0,0)$, making it naturally hard for our StLLN/average-based algorithm to hit the thin region where the minimum is located. We shift, rotate and asymmetrize the function in the same way as that in our paper for other deterministic test functions, and report the numerical results in Table \ref{tab:CrossLeg} below: the left panel is based on results with 3 starting points, while the right panel is based on results with 25 starting points (for SMCO and local algorithms only). We note that all the algorithms perform poorly in terms of AE99;  In terms of RMSE, GenSA performs the best and ADAM performs relatively well with 25 starting points, while all other algorithms perform poorly with similar RMSE. In summary, our new SMCO algorithms fail to find global min value in this example but so are most other state-of-the-art algorithms.

\begin{sidewaystable}
\centering{}\caption{\label{tab:Comp_SS_d1}Dimension $d=1$,
Uniform Single-Start}
\begin{tabular}{|c|ccccccccccc|}
\cline{1-2} \cline{2-8} \cline{8-12} \cline{9-12} \cline{10-12} \cline{11-12} \cline{12-12} 
\multicolumn{2}{|c|}{Test Function} & \multicolumn{5}{c|}{\textbf{Rastrigin}} & \multicolumn{5}{c|}{\textbf{Griewank}}\tabularnewline
\hline 
\multicolumn{2}{|c|}{Performance} & \multicolumn{1}{c|}{RMSE} & \multicolumn{1}{c|}{AE50} & \multicolumn{1}{c|}{AE95} & \multicolumn{1}{c|}{AE99} & \multicolumn{1}{c|}{Time} & \multicolumn{1}{c|}{RMSE} & \multicolumn{1}{c|}{AE50} & \multicolumn{1}{c|}{AE95} & \multicolumn{1}{c|}{AE99} & Time\tabularnewline
\hline 
\multirow{9}{*}{Max} & \textbf{Best Value} & \multicolumn{5}{c|}{\textbf{130.8169}} & \multicolumn{5}{c|}{\textbf{427.1695}}\tabularnewline
\cline{3-12} \cline{4-12} \cline{5-12} \cline{6-12} \cline{7-12} \cline{8-12} \cline{9-12} \cline{10-12} \cline{11-12} \cline{12-12} 
 & GD & 61.42 & 2.94 & 121.06 & 121.06 & 1.2 & 332.20 & 367.38 & 424.88 & 425.11 & 19\tabularnewline
 & SignGD & 83.58 & 90.46 & 110.57 & 110.57 & 18 & 332.20 & 367.38 & 424.88 & 425.11 & 15\tabularnewline
 & SPSA & 82.41 & 90.46 & 110.57 & 110.57 & 1.6 & 332.31 & 367.61 & 424.93 & 425.16 & 14\tabularnewline
 & L-BFGS & 58.60 & 2.94 & 121.06 & 121.06 & 0.19 & 330.66 & 365.86 & 424.88 & 425.11 & 0.26\tabularnewline
 & BOBYQA & 48.85 & 2.94 & 110.57 & 121.06 & 1.3 & 170.70 & 5.68E-14 & 409.13 & 409.77 & 1.4\tabularnewline
 & \textbf{SMCO} & 0.302 & 0.135 & 0.659 & 0.716 & 3.8 & 3.67 & 3.21 & 7.49 & 8.49 & 4.3\tabularnewline
 & \textbf{SMCO-R} & 0.0188 & 3.96e-7 & 1.66e-6 & 0.0714 & 4.1 & 0.185 & 0.012 & 0.085 & 0.150 & 4.6\tabularnewline
 & \textbf{SMCO-BR} & 0.438 & 7.12e-7 & 5.86e-6 & 2.94 & 4.0 & 0.489 & 0.024 & 0.117 & 0.178 & 4.5\tabularnewline
\hline 
\multirow{10}{*}{Min} & \textbf{Best Value} & \multicolumn{5}{c|}{\textbf{0}} & \multicolumn{5}{c|}{\textbf{0}}\tabularnewline
\cline{3-12} \cline{4-12} \cline{5-12} \cline{6-12} \cline{7-12} \cline{8-12} \cline{9-12} \cline{10-12} \cline{11-12} \cline{12-12} 
 & GD & 14.43 & 5.61 & 35.82 & 52.44 & 4.2 & 174.83 & 73.81 & 367.63 & 414.60 & 19\tabularnewline
 & SignGD & 52.15 & 15.92 & 100.54 & 120.38 & 19 & 174.83 & 73.81 & 367.63 & 414.60 & 15\tabularnewline
 & ADAM & 14.95 & 0.995 & 24.87 & 80.59 & 2.1 & 174.83 & 73.81 & 367.63 & 414.60 & 1.4\tabularnewline
 & SPSA & 42.72 & 15.92 & 80.59 & 120.38 & 2.0 & 174.91 & 73.83 & 367.92 & 414.60 & 16\tabularnewline
 & L-BFGS & 39.11 & 3.98 & 99.49 & 120.38 & 0.25 & 166.63 & 64.72 & 356.30 & 414.60 & 0.26\tabularnewline
 & BOBYQA & 33.65 & 3.98 & 99.49 & 120.38 & 1.3 & 26.76 & 0.0395 & 0.799 & 3.29 & 1.4\tabularnewline
\cline{3-12} \cline{4-12} \cline{5-12} \cline{6-12} \cline{7-12} \cline{8-12} \cline{9-12} \cline{10-12} \cline{11-12} \cline{12-12} 
 & \textbf{SMCO} & 1.73 & 1.08 & 4.06 & 4.27 & 4.1 & 8.11 & 7.92 & 10.44 & 11.55 & 4.4\tabularnewline
 & \textbf{SMCO-R} & 0.594 & 1.43e-6 & 0.995 & 0.995 & 4.3 & 0.0159 & 0.00986 & 0.0395 & 0.0400 & 4.9\tabularnewline
 & \textbf{SMCO-BR} & 0.580 & 2.89e-6 & 0.995 & 0.995 & 4.3 & 0.0203 & 0.00987 & 0.0395 & 0.0888 & 4.7\tabularnewline
\hline 
\end{tabular}

Based on 500 replications. ``Time'' is reported in $10^{-3}$ second. SMCO: max iter = 200, tol = 1e-8.
\end{sidewaystable}

\begin{sidewaystable}
\centering{}\caption{Dimension $d=2$, Uniform Multi-Start}\label{tab:Comp2d}
\begin{tabular}{|c|ccccccccccccc|}
\hline 
\multicolumn{2}{|c|}{Test Function} & \multicolumn{3}{c|}{\textbf{Rastrigin}} & \multicolumn{3}{c|}{\textbf{Griewank}} & \multicolumn{3}{c|}{\textbf{Ackley}} & \multicolumn{3}{c|}{\textbf{Michalewicz}}\tabularnewline
\hline 
\multicolumn{2}{|c|}{Performance} & \multicolumn{1}{c|}{RMSE} & \multicolumn{1}{c|}{AE99} & \multicolumn{1}{c|}{Time} & \multicolumn{1}{c|}{RMSE} & \multicolumn{1}{c|}{AE99} & \multicolumn{1}{c|}{Time} & \multicolumn{1}{c|}{RMSE} & \multicolumn{1}{c|}{AE99} & \multicolumn{1}{c|}{Time} & \multicolumn{1}{c|}{RMSE} & \multicolumn{1}{c|}{AE99} & Time\tabularnewline
\hline 
\multirow{16}{*}{Max} & \textbf{Best} & \multicolumn{3}{c|}{\textbf{213.5824}} & \multicolumn{3}{c|}{\textbf{622.3434}} & \multicolumn{3}{c|}{\textbf{22.35029 }} & \multicolumn{3}{c|}{\textbf{1.21406 }}\tabularnewline
\cline{3-14} \cline{4-14} \cline{5-14} \cline{6-14} \cline{7-14} \cline{8-14} \cline{9-14} \cline{10-14} \cline{11-14} \cline{12-14} \cline{13-14} \cline{14-14} 
 & GD & 1.99e-13  & 1.99e-13 & 0.048  & 192.52  & 370.71  & 0.55  & 2.56e-3  & 8.79e-3  & 0.11  & 0.220  & 0.413  & 0.18 \tabularnewline
 & SignGD & 58.70  & 100.71  & 0.44  & 215.45  & 370.71  & 0.44  & 2.46e-3  & 8.48e-3  & 0.12  & 0.180  & 0.413  & 0.39 \tabularnewline
 & SPSA & 1.44  & 1.99e-13  & 0.058  & 221.90  & 370.80  & 0.24  & 2.98e-3  & 0.0111  & 0.12  & 0.263  & 0.413  & 0.027 \tabularnewline
 & L-BFGS & 6.48  & 16.50  & 0.003  & 117.35  & 303.68  & 0.007  & 2.73e-3  & 0.0112  & 0.28  & 0.156  & 0.413  & 0.008 \tabularnewline
 & BOBYQA & 1.73e-13  & 1.99e-13  & 0.016  & 2.87e-13 & 3.41e-13  & 0.020  & 1.36e-3 & 5.01e-3  & 0.51  & 0.110  & 0.413  & 0.037 \tabularnewline
\cline{3-14} \cline{4-14} \cline{5-14} \cline{6-14} \cline{7-14} \cline{8-14} \cline{9-14} \cline{10-14} \cline{11-14} \cline{12-14} \cline{13-14} \cline{14-14} 
 & \textbf{SMCO} & 3.31  & 10.34  & 0.076 & 0.754  & 2.21  & 0.11  & 0.148  & 0.240  & 0.11  & 8.73e-4  & 2.85e-3  & 0.12 \tabularnewline
 & \textbf{SMCO-R} & 0.737  & 1.87  & 0.077 & 3.41e-13  & 3.41e-13  & 0.11  & 2.35e-3  & 6.70e-3  & 0.12  & 1.84e-7  & 4.82e-7  & 0.12 \tabularnewline
 & \textbf{SMCO-BR} & 0.867  & 4.24  & 0.079 & 4.00e-3  & 3.41e-13  & 0.12  & 1.73e-3  & 4.82e-3 & 0.12  & 3.24e-7  & 9.48e-7  & 0.12 \tabularnewline
\cline{3-14} \cline{4-14} \cline{5-14} \cline{6-14} \cline{7-14} \cline{8-14} \cline{9-14} \cline{10-14} \cline{11-14} \cline{12-14} \cline{13-14} \cline{14-14} 
 & GenSA & 1.99e-13  & 1.99e-13  & 0.24  & 3.41e-13  & 3.41e-13  & 0.26  & 1.66e-4  & 5.07e-4  & 0.35  & 2.21e-16  & 2.22e-16  & 0.50 \tabularnewline
 & SA & 2.76  & 10.29  & 0.010  & 0.0465  & 0.0684  & 0.015  & 0.00612  & 0.0141 & 0.01  & 0.0580  & 0.186  & 0.015 \tabularnewline
 & DEoptim & 8.19  & 20.30  & 0.024  & 0.0343  & 0.0529  & 0.034  & 0.00154  & 4.57e-3  & 0.04  & 0.0309  & 0.214  & 0.036 \tabularnewline
 & CMAES & 61.51  & 168.10  & 0.050  & 369.77  & 616.53  & 0.074  & 0.390  & 1.67 & 0.17  & 0.340  & 1.19  & 0.14 \tabularnewline
 & STOGO & 10.25  & 10.25  & 0.32  & 4.34  & 4.34  & 0.47  & 8.51e-5  & 8.51e-5  & 0.48  & 6.17e-10  & 6.17e-10  & 0.49 \tabularnewline
 & GA & 13.38  & 32.36  & 0.11  & 17.80  & 33.72  & 0.15  & 0.00607  & 0.0120  & 0.14  & 0.0754  & 0.413  & 0.14 \tabularnewline
 & PSO & 5.48 & 1.99e-13  & 0.25  & 3.41e-13 & 3.41e-13  & 0.35  & 7.10e-4 & 2.49e-3  & 0.34  & 0.119 & 0.315  & 0.34 \tabularnewline
\hline 
\multirow{17}{*}{Min} & \textbf{Best} & \multicolumn{3}{c|}{\textbf{0}} & \multicolumn{3}{c|}{\textbf{0}} & \multicolumn{3}{c|}{\textbf{4.440892e-16}} & \multicolumn{3}{c|}{\textbf{-1.995552 }}\tabularnewline
\cline{3-14} \cline{4-14} \cline{5-14} \cline{6-14} \cline{7-14} \cline{8-14} \cline{9-14} \cline{10-14} \cline{11-14} \cline{12-14} \cline{13-14} \cline{14-14} 
 & GD & 5.42  & 12.70  & 0.14  & 16.19  & 45.56  & 0.41  & 14.00  & 19.54  & 0.26  & 1.11  & 1.51  & 0.16 \tabularnewline
 & SignGD & 5.31  & 18.92  & 0.45  & 16.20  & 45.56  & 0.34  & 14.69  & 19.51  & 0.44  & 1.02  & 1.15  & 0.30 \tabularnewline
 & ADAM & 2.03  & 4.97  & 0.024  & 7.53  & 25.89  & 0.034  & 14.55  & 19.51  & 0.027  & 0.72  & 1.15  & 0.053 \tabularnewline
 & SPSA & 1.07  & 3.98  & 0.095  & 16.43  & 45.82  & 0.18  & 14.70  & 19.53  & 0.006  & 1.18  & 1.63  & 0.074 \tabularnewline
 & L-BFGS & 4.42  & 12.93  & 0.004  & 9.63  & 35.71  & 0.005  & 12.06  & 19.51  & 0.005  & 1.02  & 1.15  & 0.006 \tabularnewline
 & BOBYQA & 1.42  & 4.49  & 0.016  & 0.0151  & 0.0531  & 0.017 & 1.15  & 4.50  & 0.017  & 1.01  & 1.15  & 0.031\tabularnewline
\cline{3-14} \cline{4-14} \cline{5-14} \cline{6-14} \cline{7-14} \cline{8-14} \cline{9-14} \cline{10-14} \cline{11-14} \cline{12-14} \cline{13-14} \cline{14-14} 
 & \textbf{SMCO} & 2.92  & 5.03 & 0.080  & 0.560  & 0.909  & 0.090  & 0.322  & 0.872  & 0.092  & 1.14  & 1.16  & 0.10 \tabularnewline
 & \textbf{SMCO-R} & 0.563  & 1.00 & 0.084  & 0.00784  & 0.0219  & 0.093  & 0.00308  & 0.00670  & 0.097  & 1.00  & 1.00  & 0.10 \tabularnewline
 & \textbf{SMCO-BR} & 0.482  & 1.00 & 0.085  & 0.0108  & 0.0300  & 0.095 & 0.00436  & 0.00883  & 0.098  & 1.00  & 1.00  & 0.094 \tabularnewline
\cline{3-14} \cline{4-14} \cline{5-14} \cline{6-14} \cline{7-14} \cline{8-14} \cline{9-14} \cline{10-14} \cline{11-14} \cline{12-14} \cline{13-14} \cline{14-14} 
 & GenSA & 0.350  & 0.995  & 0.23  & 0.0109  & 0.0296  & 0.24  & 0 & 0 & 0.60  & 1.00  & 1.00  & 0.30 \tabularnewline
 & SA & 1.77  & 4.27  & 0.007  & 0.526  & 1.18  & 0.008  & 3.95  & 6.79 & 0.008  & 1.00  & 1.01  & 0.012 \tabularnewline
 & DEoptim & 0.192  & 0.995  & 0.022  & 0.00847  & 0.0219  & 0.025  & 0 & 0 & 0.026  & 1.00  & 1.00  & 0.027\tabularnewline
 & CMAES & 11.66  & 24.87  & 0.10  & 78.72  & 240.00  & 0.14  & 17.10  & 20.70  & 0.13  & 1.26  & 2.00  & 0.085 \tabularnewline
 & STOGO & 3.75e-5  & 3.75e-5  & 0.30  & 3.41e-10  & 3.41e-10  & 0.34  & 0.00376  & 0.00376  & 0.35  & 1.00  & 1.00  & 0.37 \tabularnewline
 & GA & 0.808  & 1.99  & 0.11  & 0.0724  & 0.241  & 0.11  & 0.311 & 1.86  & 0.11  & 1.00  & 1.01  & 0.11 \tabularnewline
 & PSO & 0.623 & 1.50  & 0.25  & 0.00377 & 0.00986 & 0.26  & 0 & 0 & 0.26  & 1.00  & 1.00  & 0.26 \tabularnewline
\hline 
\end{tabular}
Based on 250 replications. 14 starting points. SMCO: max iter = 200, tol = 1e-8.
\end{sidewaystable}

\begin{sidewaystable}
\centering{}\caption{\label{tab:Comp_10d}Dimension $d=10$,
Uniform Multi-Start}
\begin{tabular}{|c|ccccccccccccc|}
\cline{1-11} \cline{2-11} \cline{3-11} \cline{4-11} \cline{5-11} \cline{6-11} \cline{7-11} \cline{8-11} \cline{9-11} \cline{10-11} \cline{11-14} 
\multicolumn{2}{|c|}{Test Function} & \multicolumn{3}{c}{\textbf{Rastrigin}} & \multicolumn{3}{c|}{\textbf{Griewank}} & \multicolumn{3}{c|}{\textbf{Ackley}} & \multicolumn{3}{c|}{\textbf{Michalewicz}}\tabularnewline
\hline 
\multicolumn{2}{|c|}{Performance} & \multicolumn{1}{c|}{RMSE} & AE99 & \multicolumn{1}{c|}{Time} & \multicolumn{1}{c|}{RMSE} & \multicolumn{1}{c|}{AE99} & \multicolumn{1}{c|}{Time} & \multicolumn{1}{c|}{RMSE} & \multicolumn{1}{c|}{AE99} & \multicolumn{1}{c|}{Time} & \multicolumn{1}{c|}{RMSE} & AE99 & Time\tabularnewline
\hline 
\multirow{16}{*}{Max} & \textbf{Best Value} & \multicolumn{3}{c|}{\textbf{1232.878}} & \multicolumn{3}{c|}{\textbf{3771.88}} & \multicolumn{3}{c|}{\textbf{22.3502}} & \multicolumn{3}{c|}{\textbf{8.051224}}\tabularnewline
\cline{3-14} \cline{4-14} \cline{5-14} \cline{6-14} \cline{7-14} \cline{8-14} \cline{9-14} \cline{10-14} \cline{11-14} \cline{12-14} \cline{13-14} \cline{14-14} 
 & GD & 89.17 & 184.95 & 1.63 & 966.68 & 1,459.17 & 4.0 & 2.51e-3 & 5.09e-3 & 2.9 & 4.81 & 5.81 & 2.8\tabularnewline
 & SignGD & 504.36 & 596.61 & 3.49 & 1918.83 & 2,314.83 & 4.0 & 2.31e-3 & 5.15e-3 & 3.1 & 4.51 & 5.57 & 0.45\tabularnewline
 & SPSA & 108.40 & 159.60 & 0.52 & 1988.77 & 2,372.31 & 0.19 & 0.0257 & 0.0491 & 0.09 & 5.61 & 6.81 & 0.11\tabularnewline
 & L-BFGS & 144.81 & 273.17 & 0.057 & 72.10 & 311.97 & 0.14 & 2.22e-3 & 4.77e-3 & 0.04 & 3.14 & 4.20 & 0.28\tabularnewline
 & BOBYQA & 135.30 & 278.02 & 0.16 & 63.74 & 306.40 & 0.04 & 2.03e-3 & 4.73e-3 & 0.18 & 2.97 & 3.98 & 0.42\tabularnewline
\cline{3-14} \cline{4-14} \cline{5-14} \cline{6-14} \cline{7-14} \cline{8-14} \cline{9-14} \cline{10-14} \cline{11-14} \cline{12-14} \cline{13-14} \cline{14-14} 
 & \textbf{SMCO} & 84.15 & 167.17 & 0.63 & 12.93 & 18.28 & 0.72 & 0.150 & 0.179 & 0.76 & 4.38 & 4.88 & 0.81\tabularnewline
 & \textbf{SMCO-R} & 52.58 & 92.77 & 0.64 & 16.91 & 23.49 & 0.74 & 7.17e-3 & 0.0163 & 0.76 & 3.26 & 4.05 & 0.83\tabularnewline
 & \textbf{SMCO-BR} & 57.10 & 110.34 & 0.65 & 37.52 & 53.17 & 0.75 & 5.51e-3 & 9.52e-3 & 0.77 & 3.15 & 3.95 & 0.84\tabularnewline
\cline{3-14} \cline{4-14} \cline{5-14} \cline{6-14} \cline{7-14} \cline{8-14} \cline{9-14} \cline{10-14} \cline{11-14} \cline{12-14} \cline{13-14} \cline{14-14} 
 & GenSA & 33.81 & 79.10 & 0.91 & 4.55e-13 & 4.55e-13 & 0.84 & 2.39e-3 & 5.42e-3 & 0.93 & 2.16 & 3.44 & 1.65\tabularnewline
 & SA & 70.63 & 204.01 & 0.016 & 79.76 & 385.30 & 0.02 & 0.150 & 0.205 & 0.02 & 4.97 & 5.58 & 0.02\tabularnewline
 & DEoptim & 18.34 & 29.17 & 0.11 & 1.48 & 2.18 & 0.13 & 0.085 & 0.118 & 0.13 & 3.54 & 4.07 & 0.14\tabularnewline
 & CMAES & 117.68 & 417.88 & 1.93 & 593.59 & 1434.37 & 0.16 & 0.274 & 0.551 & 1.28 & 3.67 & 6.08 & 1.93\tabularnewline
 & STOGO & 320.71 & 320.71 & 1.12 & Error & Error & 1.23 & 0 & 0 & 1.30 & 2.88 & 2.88 & 1.40\tabularnewline
 & GA & 301.36 & 420.60 & 0.21 & 677.92 & 890.64 & 0.21 & 0.122 & 0.189 & 0.22 & 4.28 & 5.25 & 0.22\tabularnewline
 & PSO & 78.35 & 215.34 & 0.37 & 297.61 & 748.59 & 0.37 & 0.070 & 0.111 & 0.38 & 2.37 & 3.86 & 0.39\tabularnewline
\hline 
\multirow{17}{*}{Min} & \textbf{Best Value} & \multicolumn{3}{c|}{\textbf{1.989918}} & \multicolumn{3}{c|}{\textbf{0}} & \multicolumn{3}{c|}{\textbf{4.440892e-16}} & \multicolumn{3}{c|}{\textbf{-8.003921}}\tabularnewline
\cline{3-14} \cline{4-14} \cline{5-14} \cline{6-14} \cline{7-14} \cline{8-14} \cline{9-14} \cline{10-14} \cline{11-14} \cline{12-14} \cline{13-14} \cline{14-14} 
 & GD & 55.16 & 95.33 & 2.06 & 142.89 & 240.02 & 3.8 & 19.58 & 19.92 & 2.2 & 4.42 & 5.45 & 2.78\tabularnewline
 & SignGD & 109.19 & 184.16 & 3.35 & 354.16 & 611.64 & 3.8 & 19.58 & 19.92 & 3.4 & 4.49 & 5.50 & 0.47\tabularnewline
 & ADAM & 90.30 & 149.75 & 0.31 & 2.84e-13 & 4.73e-13 & 0.95 & 19.58 & 19.92 & 0.27 & 3.30 & 4.39 & 1.2\tabularnewline
 & SPSA & 82.54 & 103.80 & 0.49 & 382.48 & 645.54 & 0.18 & 19.61 & 19.94 & 0.08 & 5.56 & 6.56 & 0.11\tabularnewline
 & L-BFGS & 73.12 & 119.39 & 0.056 & 1.92e-03 & 7.40e-3 & 0.15 & 19.58 & 19.92 & 0.05 & 3.06 & 4.02 & 0.26\tabularnewline
 & BOBYQA & 69.45 & 130.30 & 0.16 & 7.16e-10 & 1.97e-9 & 0.22 & 18.79 & 19.91 & 0.18 & 2.89 & 3.94 & 0.42\tabularnewline
\cline{3-14} \cline{4-14} \cline{5-14} \cline{6-14} \cline{7-14} \cline{8-14} \cline{9-14} \cline{10-14} \cline{11-14} \cline{12-14} \cline{13-14} \cline{14-14} 
 & \textbf{SMCO} & 37.54 & 46.59 & 0.68 & 0.968 & 1.03 & 0.78 & 3.08 & 3.53 & 0.80 & 3.88 & 4.43 & 0.87\tabularnewline
 & \textbf{SMCO-R} & 15.28 & 29.86 & 0.69 & 0.175 & 0.235 & 0.79 & 0.0955 & 0.0850 & 0.81 & 3.05 & 3.95 & 0.89\tabularnewline
 & \textbf{SMCO-BR} & 16.86 & 29.64 & 0.70 & 0.203 & 0.267 & 0.80 & 0.919 & 1.67 & 0.82 & 3.02 & 3.81 & 0.90\tabularnewline
\cline{3-14} \cline{4-14} \cline{5-14} \cline{6-14} \cline{7-14} \cline{8-14} \cline{9-14} \cline{10-14} \cline{11-14} \cline{12-14} \cline{13-14} \cline{14-14} 
 & GenSA & 23.48 & 45.28 & 0.83 & 5.44e-3 & 0.0246 & 1.3 & 1.66 & 3.22 & 1.0 & 2.14 & 3.35 & 1.6\tabularnewline
 & SA & 106.34 & 149.25 & 0.009 & 55.55 & 118.28 & 0.01 & 18.07 & 20.68 & 0.01 & 4.88 & 5.58 & 0.02\tabularnewline
 & DEoptim & 25.32 & 33.57 & 0.11 & 0.367 & 0.548 & 0.12 & 3.60 & 6.31 & 0.13 & 3.48 & 4.15 & 0.14\tabularnewline
 & CMAES & 45.91 & 116.10 & 1.73 & 23.41 & 91.16 & 1.7 & 20.33 & 21.60 & 1.30 & 3.87 & 6.32 & 2.0\tabularnewline
 & STOGO & 68.65 & 68.65 & 1.07 & 6.29e-4 & 6.29e-4 & 1.2 & 18.54 & 18.54 & 1.2 & 3.13 & 3.13 & 1.4\tabularnewline
 & GA & 49.94 & 81.75 & 0.22 & 12.41 & 24.85 & 0.22 & 14.34 & 20.54 & 0.22 & 4.29 & 5.50 & 0.23\tabularnewline
 & PSO & 13.68 & 33.15 & 0.36 & 0.132 & 0.487 & 0.37 & 2.20 & 13.02 & 0.37 & 2.31 & 3.65 & 0.40\tabularnewline
\hline 
\end{tabular}

Based on 250 replications.  32 starting points. SMCO: max iter = 200, tol = 1e-8. 
\end{sidewaystable}

\begin{sidewaystable}
\centering{}\caption{\label{tab:Comp_20d}Dimension $d=20$, Uniform Multi-Start }
\begin{tabular}{|c|ccccccccccccc|}
\hline 
\multicolumn{2}{|c|}{Test Function} & \multicolumn{3}{c}{\textbf{Rastrigin}} & \multicolumn{3}{c|}{\textbf{Griewank}} & \multicolumn{3}{c|}{\textbf{Ackley}} & \multicolumn{3}{c|}{\textbf{Michalewicz}}\tabularnewline
\hline 
\multicolumn{2}{|c|}{Performance} & \multicolumn{1}{c|}{RMSE} & \multicolumn{1}{c|}{AE99} & \multicolumn{1}{c|}{Time} & \multicolumn{1}{c|}{RMSE} & \multicolumn{1}{c|}{AE99} & \multicolumn{1}{c|}{Time} & \multicolumn{1}{c|}{RMSE} & \multicolumn{1}{c|}{AE99} & \multicolumn{1}{c|}{Time} & \multicolumn{1}{c|}{RMSE} & \multicolumn{1}{c|}{AE99} & Time\tabularnewline
\hline 
\multirow{16}{*}{Max} & \textbf{Best} & \multicolumn{3}{c|}{\textbf{2417.617}} & \multicolumn{3}{c|}{\textbf{7240.209 }} & \multicolumn{3}{c|}{\textbf{22.34938 }} & \multicolumn{3}{c|}{\textbf{14.37489 }}\tabularnewline
\cline{3-14} \cline{4-14} \cline{5-14} \cline{6-14} \cline{7-14} \cline{8-14} \cline{9-14} \cline{10-14} \cline{11-14} \cline{12-14} \cline{13-14} \cline{14-14} 
 & GD & 253.35  & 399.29  & 6.7  & 2276.85  & 2,858.55  & 12  & 3.54e-3  & 6.57e-3  & 7.3  & 8.67  & 10.08  & 10.7 \tabularnewline
 & SignGD & 1,113.38  & 1,239.69  & 10 & 4,056.59  & 4,472.41  & 12  & 3.51e-3  & 6.44e-3 & 11.5  & 8.48  & 9.67  & 1.6 \tabularnewline
 & SPSA & 638.42  & 767.70  & 0.86  & 4,184.55  & 4,596.13  & 0.32  & 0.0908  & 0.130  & 0.15  & 11.70  & 12.71  & 1.0 \tabularnewline
 & L-BFGS & 445.62  & 614.87  & 0.23  & 279.33  & 647.21  & 1.1  & 3.47e-3  & 6.44e-3  & 0.14  & 6.57  & 7.88  & 1.6 \tabularnewline
 & BOBYQA & 747.83  & 987.66  & 1.1  & 277.67  & 627.07  & 0.086 & 3.46e-3  & 6.58e-3  & 1.2  & 6.25  & 7.66  & 1.3 \tabularnewline
\cline{3-14} \cline{4-14} \cline{5-14} \cline{6-14} \cline{7-14} \cline{8-14} \cline{9-14} \cline{10-14} \cline{11-14} \cline{12-14} \cline{13-14} \cline{14-14} 
 & \textbf{SMCO} & 632.59  & 731.91  & 1.9  & 30.34  & 37.95  & 2.3  & 0.174  & 0.196  & 2.4  & 8.45  & 9.06  & 2.6 \tabularnewline
 & \textbf{SMCO-R} & 353.79  & 433.02  & 1.91  & 48.23  & 58.73  & 2.4  & 8.89e-3  & 0.0159  & 2.4  & 6.89  & 7.68  & 2.7 \tabularnewline
 & \textbf{SMCO-BR} & 367.18  & 446.75  & 1.9  & 101.57  & 121.98  & 2.4  & 8.68e-3  & 0.0137  & 2.4  & 6.72  & 7.58  & 2.7 \tabularnewline
\cline{3-14} \cline{4-14} \cline{5-14} \cline{6-14} \cline{7-14} \cline{8-14} \cline{9-14} \cline{10-14} \cline{11-14} \cline{12-14} \cline{13-14} \cline{14-14} 
 & GenSA & 219.27  & 410.31  & 2.1  & 9.09e-13  & 9.09e-13  & 2.2  & 4.77e-3  & 9.91e-3  & 2.3  & 4.96  & 6.83  & 4.8 \tabularnewline
 & SA & 341.55  & 540.56  & 0.02  & 394.09  & 935.70  & 0.022  & 0.252  & 0.307  & 0.023  & 10.19  & 11.11  & 0.028 \tabularnewline
 & DEoptim & 252.20  & 299.05  & 0.26  & 138.09  & 168.08  & 0.31  & 0.174  & 0.208  & 0.31  & 8.30  & 9.02  & 0.35 \tabularnewline
 & CMAES & 230.09  & 636.56  & 9.7  & 758.94  & 1,863.78  & 0.64  & 0.360  & 0.595  & 4.2  & 6.38  & 10.80  & 8.9 \tabularnewline
 & STOGO & 1,385.39  & 1,385.39  & 2.3  & ERROR & ERROR  & 2.8  & 6.19e-3  & 6.19e-3  & 2.9  & 7.65  & 7.65  & 3.1 \tabularnewline
 & GA & 988.53  & 1,218.89  & 0.24  & 2,134.93  & 2,567.87  & 0.25  & 0.205  & 0.291  & 0.26  & 8.70  & 10.23  & 0.25 \tabularnewline
 & PSO & 266.43  & 550.05 & 0.46  & 829.96 & 1,877.95 & 0.49  & 0.162  & 0.217  & 0.52  & 6.92  & 8.57  & 0.53 \tabularnewline
\hline 
\multirow{17}{*}{Min} & \textbf{Best} & \multicolumn{3}{c|}{\textbf{10.94455}} & \multicolumn{3}{c|}{\textbf{2.220446e-16 }} & \multicolumn{3}{c|}{\textbf{2.586065e-11}} & \multicolumn{3}{c|}{\textbf{-14.48265 }}\tabularnewline
\cline{3-14} \cline{4-14} \cline{5-14} \cline{6-14} \cline{7-14} \cline{8-14} \cline{9-14} \cline{10-14} \cline{11-14} \cline{12-14} \cline{13-14} \cline{14-14} 
 & GD & 134.67  & 206.36  & 6.0  & 377.08  & 492.81  & 11  & 19.82  & 19.93  & 5.4  & 8.36  & 9.77  & 9.9 \tabularnewline
 & SignGD & 274.76  & 370.22  & 9.3  & 951.60  & 1,253.15  & 11  & 19.82  & 19.93  & 11  & 8.55  & 9.93  & 1.4 \tabularnewline
 & ADAM & 244.08  & 331.84  & 0.53  & 6.09e-13  & 8.84e-13  & 2.6  & 19.82  & 19.93  & 0.51  & 7.14  & 8.62  & 6.2 \tabularnewline
 & SPSA & 248.98  & 304.39  & 0.79  & 1,017.81  & 1330.50  & 0.30  & 19.87  & 19.93  & 0.20  & 11.75  & 12.85  & 1.0 \tabularnewline
 & L-BFGS & 190.98  & 260.00  & 0.18  & 1.27e-10  & 2.68e-10  & 0.36  & 19.82  & 19.93  & 0.15  & 6.31  & 8.00  & 1.4 \tabularnewline
 & BOBYQA & 236.82  & 324.35  & 1.1  & 7.67e-9  & 2.20e-8 & 1.2  & 19.81  & 19.93  & 1.1  & 6.27  & 7.59  & 1.2 \tabularnewline
\cline{3-14} \cline{4-14} \cline{5-14} \cline{6-14} \cline{7-14} \cline{8-14} \cline{9-14} \cline{10-14} \cline{11-14} \cline{12-14} \cline{13-14} \cline{14-14} 
 & \textbf{SMCO} & 85.32  & 101.94  & 1.9 & 0.992  & 1.02  & 3.5  & 5.82  & 19.92  & 3.5  & 7.16  & 7.61  & 3.9 \tabularnewline
 & \textbf{SMCO-R} & 63.52  & 98.04  & 1.9 & 0.340  & 0.423  & 3.5  & 2.24  & 2.34  & 3.5  & 6.53  & 7.04  & 4.0 \tabularnewline
 & \textbf{SMCO-BR} & 60.01  & 92.56  & 2.0 & 0.414  & 0.510  & 3.5  & 1.83  & 2.34  & 3.5  & 6.30  & 7.03  & 4.0 \tabularnewline
\cline{3-14} \cline{4-14} \cline{5-14} \cline{6-14} \cline{7-14} \cline{8-14} \cline{9-14} \cline{10-14} \cline{11-14} \cline{12-14} \cline{13-14} \cline{14-14} 
 & GenSA & 74.49  & 130.24  & 1.9 & 9.31e-13  & 2.02e-12  & 2.8  & 11.53  & 19.93  & 2.1  & 4.88  & 7.00  & 4.1 \tabularnewline
 & SA & 314.88  & 399.62  & 0.01 & 274.79  & 498.07  & 0.009  & 20.88  & 21.09  & 0.021  & 10.32  & 11.13  & 0.026 \tabularnewline
 & DEoptim & 108.35  & 127.56  & 0.24 & 1.05  & 1.07  & 0.28  & 14.93  & 18.40  & 0.28  & 8.38  & 9.14  & 0.32 \tabularnewline
 & CMAES & 90.65  & 233.08  & 5.4 & 48.11  & 147.13 & 4.6  & 20.48  & 21.64  & 4.2  & 6.52  & 11.20  & 8.5 \tabularnewline
 & STOGO & 174.12  & 174.12  & 2.1 & 2.94e-6  & 2.94e-6  & 2.6  & 19.34  & 19.34  & 2.6  & 8.46  & 8.46  & 2.9 \tabularnewline
 & GA & 171.04  & 236.88  & 0.24 & 90.23  & 139.16  & 0.25  & 20.05  & 20.90  & 0.25  & 8.81  & 10.32  & 0.25 \tabularnewline
 & PSO & 55.44  & 104.09  & 0.43 & 0.018  & 0.055  & 0.45  & 3.63  & 20.77  & 0.46  & 7.27  & 8.79  & 0.50 \tabularnewline
\hline 
\end{tabular}

Based on 250 replications. SMCO \& Group I: 45 starting points. SMCO: max iter = 200, tol = 1e-8.
\end{sidewaystable}

\begin{sidewaystable}
\centering{}\caption{\label{tab:Comp50d}Dimension $d=50$, Uniform Multi-Start}
\begin{tabular}{|c|ccccccccccccc|}
\hline 
\multicolumn{2}{|c|}{Test Function} & \multicolumn{3}{c|}{\textbf{Rastrigin}} & \multicolumn{3}{c|}{\textbf{Griewank}} & \multicolumn{3}{c|}{\textbf{Ackley}} & \multicolumn{3}{c|}{\textbf{Michalewicz}}\tabularnewline
\hline 
\multicolumn{2}{|c|}{Performance} & \multicolumn{1}{c|}{RMSE} & \multicolumn{1}{c|}{AE99} & \multicolumn{1}{c|}{Time} & \multicolumn{1}{c|}{RMSE} & \multicolumn{1}{c|}{AE99} & \multicolumn{1}{c|}{Time} & \multicolumn{1}{c|}{RMSE} & \multicolumn{1}{c|}{AE99} & \multicolumn{1}{c|}{Time} & \multicolumn{1}{c|}{RMSE} & \multicolumn{1}{c|}{AE99} & Time\tabularnewline
\hline 
\multirow{16}{*}{Max} & \textbf{Best} & \multicolumn{3}{c}{\textbf{5860.922}} & \multicolumn{3}{c|}{\textbf{18061.87 }} & \multicolumn{3}{c|}{\textbf{22.34741 }} & \multicolumn{3}{c|}{\textbf{25.11029}}\tabularnewline
\cline{3-14} \cline{4-14} \cline{5-14} \cline{6-14} \cline{7-14} \cline{8-14} \cline{9-14} \cline{10-14} \cline{11-14} \cline{12-14} \cline{13-14} \cline{14-14} 
 & GD & 615.25  & 868.45  & 37  & 6,867.83  & 7634.46  & 59 & 4.13e-3  & 7.06e-3  & 22  & 12.83  & 15.05  & 56 \tabularnewline
 & SignGD & 2,879.52  & 3059.61  & 49  & 10,992.09  & 11,651.92  & 59 & 4.18e-3  & 7.10e-3 & 56  & 17.50  & 18.76  & 7.7 \tabularnewline
 & SPSA & 2,533.16  & 2758.01  & 1.9  & 11,286.52  & 11,936.42  & 0.95 & 0.236  & 0.282  & 0.28  & 20.91  & 22.53  & 2.7 \tabularnewline
 & L-BFGS & 1,288.07  & 1547.96  & 1.7  & 1,280.51  & 1,816.09  & 12 & 4.13e-3  & 7.06e-3  & 0.70  & 10.43  & 11.98  & 9.2 \tabularnewline
 & BOBYQA & 2,718.94 & 2914.27  & 8.1 & 1,280.51 & 1,816.09  & 0.38 & 7.35e-3  & 0.0111 & 8.2  & 9.91  & 11.52  & 8.0 \tabularnewline
\cline{3-14} \cline{4-14} \cline{5-14} \cline{6-14} \cline{7-14} \cline{8-14} \cline{9-14} \cline{10-14} \cline{11-14} \cline{12-14} \cline{13-14} \cline{14-14} 
 & \textbf{SMCO} & 2,307.95  & 2,422.66  & 9.3 & 89.09  & 100.55  & 12  & 0.191  & 0.204  & 11  & 15.09  & 15.92  & 14 \tabularnewline
 & \textbf{SMCO-R} & 1,424.85  & 1,581.45  & 9.5 & 152.51  & 168.25  & 12  & 0.0102  & 0.0150  & 11  & 10.28  & 12.25  & 14 \tabularnewline
 & \textbf{SMCO-BR} & 1,482.92  & 1,646.48 & 9.6 & 314.89  & 342.48  & 12  & 0.0118  & 0.0163  & 11  & 9.65  & 11.39  & 14 \tabularnewline
\cline{3-14} \cline{4-14} \cline{5-14} \cline{6-14} \cline{7-14} \cline{8-14} \cline{9-14} \cline{10-14} \cline{11-14} \cline{12-14} \cline{13-14} \cline{14-14} 
 & GenSA & 1,167.86  & 1,772.72  & 7.9 & 0 & 0 & 8.8  & 6.61e-3  & 0.0109  & 8.2  & 6.23  & 8.80  & 18 \tabularnewline
 & SA & 1,556.56  & 1,917.50  & 0.045 & 3,539.77  & 4688.62  & 0.035  & 0.363  & 0.413  & 0.039  & 19.14  & 20.34  & 0.052 \tabularnewline
 & DEoptim & 1,805.04  & 1,937.33  & 0.87 & 3,292.97  & 3484.97  & 1.0  & 0.284  & 0.312  & 1.0  & 16.23  & 17.20  & 1.3 \tabularnewline
 & CMAES & 418.84  & 868.91  & 136 & 672.77  & 1,614.30  & 8.4  & 0.428  & 0.580  & 28  & 11.84  & 19.75  & 74 \tabularnewline
 & STOGO & 3,837.20  & 3,837.20  & 7.3 & ERROR & ERROR & 8.5  & 0.113  & 0.113  & 8.4  & 15.58  & 15.58  & 11 \tabularnewline
 & GA & 3,113.40  & 3,361.89  & 0.28 & 8,637.59  & 9,412.69  & 0.29  & 0.308  & 0.389  & 0.28  & 16.60  & 19.08  & 0.29 \tabularnewline
 & PSO & 827.38 & 1,475.95 & 0.74 & 3,097.68  & 5,063.57  & 0.76  & 0.290  & 0.328  & 0.8  & 14.39  & 16.16  & 0.90 \tabularnewline
\hline 
\multirow{17}{*}{Min} & \textbf{Best} & \multicolumn{3}{c|}{\textbf{69.80496}} & \multicolumn{3}{c|}{\textbf{7.320811e-13 }} & \multicolumn{3}{c|}{\textbf{1.179612 }} & \multicolumn{3}{c|}{\textbf{-25.24325 }}\tabularnewline
\cline{3-14} \cline{4-14} \cline{5-14} \cline{6-14} \cline{7-14} \cline{8-14} \cline{9-14} \cline{10-14} \cline{11-14} \cline{12-14} \cline{13-14} \cline{14-14} 
 & GD & 350.64  & 462.75  & 32 & 1,186.21  & 1439.17  & 56 & 18.74  & 18.77  & 34  & 13.19  & 15.06  & 56 \tabularnewline
 & SignGD & 853.17  & 1,011.03  & 49 & 3,015.98  & 3685.10  & 56 & 18.74  & 18.77  & 53  & 17.60  & 18.93  & 7.3 \tabularnewline
 & ADAM & 792.90  & 937.19  & 1.4 & 1.97e-13  & 3.00e-13  & 9.5 & 18.74  & 18.77  & 1.4  & 10.85  & 12.59  & 28 \tabularnewline
 & SPSA & 966.96  & 1,077.18  & 1.9 & 3204.76  & 3,890.44  & 0.93 & 18.91  & 19.10  & 0.44  & 21.09  & 22.54  & 2.6 \tabularnewline
 & L-BFGS & 599.33  & 711.08  & 1.0 & 4.15e-10  & 5.99e-10  & 1.0 & 18.74  & 18.77  & 0.86  & 10.12  & 11.84  & 9.0 \tabularnewline
 & BOBYQA & 813.91  & 943.09 & 8.3 & 1.70e-3 & 5.12e-3 & 8.6 & 18.76  & 18.79  & 7.9  & 9.10  & 10.87  & 7.9 \tabularnewline
\cline{3-14} \cline{4-14} \cline{5-14} \cline{6-14} \cline{7-14} \cline{8-14} \cline{9-14} \cline{10-14} \cline{11-14} \cline{12-14} \cline{13-14} \cline{14-14} 
 & \textbf{SMCO} & 264.32  & 296.12  & 10  & 1.08  & 1.09  & 12  & 19.06  & 19.18  & 12  & 16.51  & 17.38  & 15 \tabularnewline
 & \textbf{SMCO-R} & 261.10  & 289.54  & 10  & 0.680  & 0.755  & 12  & 18.37  & 18.63  & 12  & 8.33  & 10.24  & 15 \tabularnewline
 & \textbf{SMCO-BR} & 258.94  & 294.48  & 10  & 0.970  & 1.01  & 12  & 18.27  & 18.63  & 12  & 9.13  & 10.68  & 15 \tabularnewline
\cline{3-14} \cline{4-14} \cline{5-14} \cline{6-14} \cline{7-14} \cline{8-14} \cline{9-14} \cline{10-14} \cline{11-14} \cline{12-14} \cline{13-14} \cline{14-14} 
 & GenSA & 348.10  & 483.46  & 7.7 & 7.51e-12  & 1.11e-11  & 9.7  & 18.74  & 18.79  & 8.0  & 5.54  & 8.34  & 17 \tabularnewline
 & SA & 1,099.59  & 1,261.13  & 0.03 & 1,403.02  & 2,217.09  & 0.01  & 20.10  & 20.18  & 0.039  & 19.09  & 20.36  & 0.055 \tabularnewline
 & DEoptim & 539.37  & 591.06  & 0.87 & 53.19  & 62.06  & 0.98  & 19.93  & 20.01  & 0.96  & 16.12  & 17.18  & 1.2 \tabularnewline
 & CMAES & 195.48  & 569.49  & 42 & 140.47  & 410.28  & 25  & 19.94  & 20.46  & 26  & 12.74  & 20.27  & 84 \tabularnewline
 & STOGO & 579.89  & 579.89  & 7.2 & 3.85e-9  & 3.85e-9  & 8.2  & 18.66  & 18.66  & 8.1  & 15.99  & 15.99  & 11 \tabularnewline
 & GA & 649.13  & 757.55  & 0.29 & 756.14  & 946.54  & 0.28  & 19.81  & 20.01  & 0.29  & 16.29  & 19.23  & 0.30\tabularnewline
 & PSO & 229.28  & 338.40 & 0.77 & 0.0175  & 0.05 & 0.75  & 10.88  & 20.03  & 0.79  & 14.47  & 16.20  & 0.91 \tabularnewline
\hline 
\end{tabular}

Based on 100 replications. SMCO \& Group I: 71 starting points. SMCO: max iter = 200, tol = 1e-8.
\end{sidewaystable}



\begin{sidewaystable}
\begin{centering}
\caption{\label{tab:Comp_200d_diag} Dimension $d=200$, Diagonal Multi-Start}
\begin{tabular}{|c|cccc|ccc|ccc|ccc|}
\hline 
\multicolumn{2}{|c|}{Test Function} & \multicolumn{3}{c|}{\textbf{Rastrigin}} & \multicolumn{3}{c|}{\textbf{Griewank}} & \multicolumn{3}{c|}{\textbf{Ackley}} & \multicolumn{3}{c|}{\textbf{Michalewicz}}\tabularnewline
\hline 
\multicolumn{2}{|c|}{Performance} & \multicolumn{1}{c|}{RMSE} & AE99 & Time & \multicolumn{1}{c|}{RMSE} & \multicolumn{1}{c|}{AE99} & Time & \multicolumn{1}{c|}{RMSE} & \multicolumn{1}{c|}{AE99} & Time & \multicolumn{1}{c|}{RMSE} & AE99 & Time\tabularnewline
\hline 
\multirow{10}{*}{Min} & \textbf{Best} & \multicolumn{3}{c|}{\textbf{1596.819 }} & \multicolumn{3}{c|}{\textbf{2.35e-13 }} & \multicolumn{3}{c|}{\textbf{19.85761 }} & \multicolumn{3}{c|}{\textbf{-62.77543 }}\tabularnewline
\cline{3-14} \cline{4-14} \cline{5-14} \cline{6-14} \cline{7-14} \cline{8-14} \cline{9-14} \cline{10-14} \cline{11-14} \cline{12-14} \cline{13-14} \cline{14-14} 
 & \textbf{SMCO} & \textbf{285 } & \textbf{365 } & 23.2  & \textbf{1.08 } & \textbf{1.09 } & 23.6  & \textbf{0.103 } & \textbf{0.112 } & 23.9  & \textbf{42.51 } & \textbf{45.05 } & 31.5 \tabularnewline
 & \textbf{SMCO-R} & \textbf{318 } & \textbf{376 } & 23.3  & \textbf{0.95 } & \textbf{0.98 } & 24.0  & \textbf{0.011 } & \textbf{0.015 } & 24.1  & \textbf{28.00 } & \textbf{29.09 } & 31.6 \tabularnewline
 & \textbf{SMCO-BR} & \textbf{287 } & \textbf{340 } & 23.4  & \textbf{0.97 } & \textbf{0.99 } & 24.7  & \textbf{0.012 } & \textbf{0.018 } & 24.2  & \textbf{30.84 } & \textbf{32.49 } & 31.6 \tabularnewline
\cline{3-14} \cline{4-14} \cline{5-14} \cline{6-14} \cline{7-14} \cline{8-14} \cline{9-14} \cline{10-14} \cline{11-14} \cline{12-14} \cline{13-14} \cline{14-14} 
 & GenSA & 1718  & 2011  & 268.3  & 3.05e-11  & 6.10e-11  & 302.2  & 0.109  & 0.114  & 271.7  & 6.27  & 12.03  & 813.4 \tabularnewline
 & SA & 3947  & 4338  & 0.2  & 11626  & 14002  & 0.1  & 1.653  & 1.678  & 0.3  & 52.76  & 54.53  & 0.36 \tabularnewline
 & DEoptim & 3813  & 3968  & 31.5  & 5856  & 6086  & 33.4  & 1.548  & 1.564  & 32.2  & 47.41  & 48.91  & 43.5 \tabularnewline
 & STOGO & 1293  & 1293  & 285.6  & 4.38e-09  & 4.38e-09  & 307.1  & 0.002  & 0.002  & 295.4  & 43.48  & 43.48  & 404.4 \tabularnewline
 & GA & 2496  & 2623  & 0.5  & 7372  & 7920  & 0.5  & 1.452  & 1.502  & 0.5  & 49.28  & 52.30  & 0.61 \tabularnewline
 & PSO & 590 & 812  & 4.9  & 81  & 143 & 5.1  & 1.548 & 1.587 & 5.1 & 43.13  & 45.36  & 6.19 \tabularnewline
\hline 
\multirow{10}{*}{Max} & \textbf{Best} & \multicolumn{3}{c|}{\textbf{18239.57 }} & \multicolumn{3}{c|}{\textbf{72331.4}} & \multicolumn{3}{c|}{\textbf{22.33939 }} & \multicolumn{3}{c|}{\textbf{69.53472 }}\tabularnewline
\cline{3-14} \cline{4-14} \cline{5-14} \cline{6-14} \cline{7-14} \cline{8-14} \cline{9-14} \cline{10-14} \cline{11-14} \cline{12-14} \cline{13-14} \cline{14-14} 
 & \textbf{SMCO} & \textbf{7219 } & \textbf{7596 } & 22.9  & \textbf{189 } & \textbf{200 } & 24.6  & \textbf{0.121 } & \textbf{0.124 } & 23.7  & \textbf{36.91 } & \textbf{38.65 } & 31.3 \tabularnewline
 & \textbf{SMCO-R} & \textbf{3918 } & \textbf{4061 } & 23.3  & \textbf{272 } & \textbf{283 } & 24.9  & \textbf{0.106 } & \textbf{0.110 } & 23.5  & \textbf{20.64 } & \textbf{21.82 } & 31.4 \tabularnewline
 & \textbf{SMCO-BR} & \textbf{4050 } & \textbf{4095 } & 23.3  & \textbf{587 } & \textbf{605 } & 25.0  & \textbf{0.012 } & \textbf{0.012 } & 23.9  & \textbf{25.22 } & \textbf{26.92 } & 31.6 \tabularnewline
 & GenSA & 2605  & 2975  & 275.9  & 0  & 0  & 312.7  & 0.004  & 0.007  & 273.5  & 4.65  & 8.98  & 814.1 \tabularnewline
 & SA & 5653  & 7180  & 0.3  & 32221  & 36635  & 0.3  & 0.488  & 0.505  & 0.3  & 55.48  & 57.14  & 0.37 \tabularnewline
 & DEoptim & 7978  & 8191  & 31.8  & 36012  & 36465  & 33.7  & 0.421  & 0.432 & 32.6  & 49.83  & 50.98  & 43.8 \tabularnewline
 & STOGO & 10987  & 10987  & 290.0  & 9908  & 9908  & 310.0  & 0.120  & 0.120  & 299.9  & 45.26  & 45.26  & 407.8 \tabularnewline
 & GA & 8977  & 9310  & 0.5  & 45978  & 46953  & 0.5  & 0.450  & 0.489  & 0.5  & 49.86  & 51.20  & 0.61 \tabularnewline
 & PSO & 998 & 1719  & 5.0  & 20792  & 23940  & 5.2  & 0.429  & 0.448  & 5.2  & 46.50 & 48.88  & 6.23 \tabularnewline
\hline 
\end{tabular} 
\par\end{centering}
\centering{}Based on: 10 reps. SMCO/-R/-BR run from 14 diagonal start points in  parallel.
\end{sidewaystable}

\begin{sidewaystable}
\centering{}
\caption{\label{tab:Comp_d500_3glob} Dimension $d=500$, Diagonal Multi-Start}
\begin{tabular}{|c|ccccccccccccc|}
\hline 
\multicolumn{2}{|c|}{} & \multicolumn{3}{c|}{Rastrigin} & \multicolumn{3}{c|}{Michalewicz}  & \multicolumn{3}{c|}{Griewank} & \multicolumn{3}{c|}{Ackley} \tabularnewline
\hline 
\multicolumn{2}{|c|}{} & \multicolumn{1}{c|}{RMSE} & \multicolumn{1}{c|}{AE99} & \multicolumn{1}{c|}{Time} & \multicolumn{1}{c|}{RMSE} & \multicolumn{1}{c|}{AE99} & \multicolumn{1}{c|}{Time} & \multicolumn{1}{c|}{RMSE} & \multicolumn{1}{c|}{AE99} & \multicolumn{1}{c|}{Time} & \multicolumn{1}{c|}{RMSE} & \multicolumn{1}{c|}{AE99} & Time\tabularnewline
\hline 

\multirow[t]{9}{*}{Max} & Best 
& \multicolumn{3}{c|}{\textbf{41635.04}} 
& \multicolumn{3}{c|}{\textbf{123.50}} 
& \multicolumn{3}{c|}{\textbf{181936.18}} 
& \multicolumn{3}{c|}{\textbf{22.35}} \tabularnewline
\cline{3-14}

& \textbf{SMCO}
& \textbf{16141.12} & \textbf{16305.43} & \textbf{207.61}
& \textbf{109.90}     & \textbf{113.51}     & \textbf{223.86}
& \textbf{778.66}    & \textbf{805.66}    & \textbf{212.92}
& \textbf{0.16}      & \textbf{0.17}      & \textbf{206.97} \tabularnewline

& \textbf{SMCO-R}
& \textbf{5939.57}  & \textbf{5945.41}  & \textbf{696.47}
& \textbf{72.93}     & \textbf{75.79}     & \textbf{228.85}   & \textbf{\textbf{201.05}} 
& \textbf{214.89}           & \textbf{724.20}
& \textbf{0.12}      & \textbf{0.12}      & \textbf{211.83} \tabularnewline

& \textbf{SMCO-BR}
& \textbf{5683.68}  & \textbf{5768.69}  & \textbf{206.89}
& \textbf{78.84}     & \textbf{80.96}     & \textbf{227.13}
& \textbf{1958.97}  & \textbf{2008.52}  & \textbf{212.96}
& \textbf{0.01}      & \textbf{0.01}      & \textbf{213.54} \tabularnewline
\cline{3-14}

& GenSA
& 2007.79  & 3382.04  & 2088.85
& 6.81      & 10.37     & 3941.08
& 0.00      & 0.00      & 2734.51
& 0.01      & 0.01      & 2085.72  \tabularnewline

& DEoptim
& 18596.52 & 18955.17 & 338.83
& 98.52     & 99.83     & 364.38
& 111051.66& 111240.77& 350.08
& 0.49      & 0.50      & 329.34 \tabularnewline

& PSO
& 3370.14  & 6057.01  & 22.12
& 91.60     & 95.21     & 24.76
& 67941.28 & 72225.00 & 23.92
& 0.49      & 0.50      & 23.35 \tabularnewline

\hline 
\multirow[t]{9}{*}{Min} & Best 
& \multicolumn{3}{c|}{\textbf{4953.99}} 
& \multicolumn{3}{c|}{\textbf{-140.11}} 
& \multicolumn{3}{c|}{\textbf{0.00}} 
& \multicolumn{3}{c|}{\textbf{19.86}} \tabularnewline
\cline{3-14}

& \textbf{SMCO}
& \textbf{188.23}    & \textbf{302.66}    & \textbf{206.59}
& \textbf{113.11}    & \textbf{115.00}    & \textbf{236.19}
& \textbf{1.49}      & \textbf{1.51}      & \textbf{215.68}
& \textbf{0.29}      & \textbf{0.29}      & \textbf{209.97} \tabularnewline

& \textbf{SMCO-R}
& \textbf{163.41}    & \textbf{248.16}    & \textbf{206.33}
& \textbf{78.85}     & \textbf{80.81}     & \textbf{226.41}
& \textbf{0.94}      & \textbf{0.97}      & \textbf{217.47}
& \textbf{0.02}      & \textbf{0.02}      & \textbf{211.02}\tabularnewline

& \textbf{SMCO-BR}
& \textbf{187.19}    & \textbf{281.60}    & \textbf{214.01}
& \textbf{84.48}     & \textbf{86.89}     & \textbf{222.84}
& \textbf{1.02}      & \textbf{1.04}      & \textbf{211.94}
& \textbf{0.02}      & \textbf{0.02}      & \textbf{209.58} \tabularnewline

\cline{3-14}

& GenSA
& 4994.21  & 5541.87  & 2057.16
& 5.88      & 6.56      & 3941.37
& 0.00      & 0.00      & 2286.15
& 0.12      & 0.13      & 2076.07 \tabularnewline

& DEoptim
& 11294.15 & 11294.15 & 354.66
& 109.41    & 109.41    & 390.09
& 29094.97 & 29094.97 & 367.80
& 1.66      & 1.66      & 355.15 \tabularnewline

& PSO
& 2495.80  & 3101.78  & 23.59
& 101.60    & 105.62    & 25.41
& 2277.66  & 2582.27  & 23.57
& 1.63      & 1.65      & 23.15 \tabularnewline

\hline 
 
\end{tabular}

Based on $10$ replications, tolerance =1e-06, maximum iteration=300 and the number of multi-start is 30 for SMCO method.
\end{sidewaystable}

\begin{table}
\caption{Sensitivity: Asymmetrization/Enlargement Transformation}
\label{tab:shift_sensitivity}
\centering
\setlength{\tabcolsep}{4pt} 
\renewcommand{\arraystretch}{1.15} 
\small
\begin{threeparttable}

\begin{tabular}{|l|c|c|c|c|c|c|}
\hline
\multicolumn{7}{|c|}{\textbf{Minimization of Rastrigin Function ($d=500$)}} \\ \hline
\hline
\multirow{2}{*}{Algorithm} & \multirow{2}{*}{\begin{tabular}[c]{@{}c@{}}Time\\(sec)\end{tabular}} & \multirow{2}{*}{RMSE} & \multicolumn{3}{c|}{Performance Metrics} & \multirow{2}{*}{Best Value} \\
\cline{4-6}
& & & AE 50\% & AE 95\% & AE 99\% & \\ \hline
\multicolumn{7}{|l|}{\textbf{Main-Text Configuration (25\% ±5\%, 50\% ±10\%)}} \\ \hline
\textbf{SMCO}       & 188.889  & 116.510    & 71.163     & 177.069    & 186.483    & 4978.24  \\ \hline
\textbf{SMCO-BR}    & 198.900  & 134.120    & 151.389    & 171.082    & 172.832    & 5010.15  \\ \hline
\textbf{SMCO-R}     & 200.381  & 154.574    & 178.063    & 197.746    & 199.496    & 5017.27  \\ \hline
\textbf{GenSA}      & 1371.450 & 5019.651   & 4893.066   & 5252.740   & 5284.711   & 9907.08  \\ \hline
\textbf{DEoptim}    & 166.636  & 10954.476  & 10917.423  & 11095.872  & 11111.734  & 15845.09 \\ \hline
\textbf{PSO}        & 10.535   & 2842.682   & 2855.954   & 2954.504   & 2963.264   & 7731.87  \\ \hline
\multicolumn{7}{|l|}{\textbf{Variant1 Configuration (25\% ±5\%, 25\% ±5\%)}} \\ \hline
\textbf{SMCO}       & 188.687  & 148.780    & 152.367    & 158.668    & 159.229    & 2705.21  \\ \hline
\textbf{SMCO-BR}    & 202.054  & 71.965     & 51.707     & 107.246    & 112.183    & 2611.88  \\ \hline
\textbf{SMCO-R}     & 200.051  & 116.283    & 127.235    & 136.212    & 137.010    & 2669.82  \\ \hline
\textbf{GenSA}      & 1371.725 & 3001.373   & 3095.112   & 3136.194   & 3139.845   & 5553.23  \\ \hline
\textbf{DEoptim}    & 166.118  & 7767.578   & 7736.598   & 7868.279   & 7879.984   & 10323.95 \\ \hline
\textbf{PSO}        & 10.397   & 4261.991   & 4238.920   & 4307.672   & 4313.784   & 6818.66  \\ \hline
\multicolumn{7}{|l|}{\textbf{Variant2 Configuration (25\% ±5\% , 75\% ±10\%)}} \\ \hline
\textbf{SMCO}       & 189.525  & 490.762    & 519.815    & 547.085    & 549.509    & 8016.38  \\ \hline
\textbf{SMCO-BR}    & 193.758  & 457.677    & 447.492    & 475.280    & 477.750    & 7988.20  \\ \hline
\textbf{SMCO-R}     & 193.144  & 483.404    & 521.899    & 537.952    & 539.378    & 8008.17  \\ \hline
\textbf{GenSA}      & 1373.213 & 6903.926   & 6854.474   & 7107.022   & 7129.471   & 14432.48 \\ \hline
\textbf{DEoptim}    & 166.342  & 15421.691  & 15338.121  & 15611.723  & 15636.043  & 22951.64 \\ \hline
\textbf{PSO}        & 10.510   & 414.609    & 36.133     & 649.106    & 703.593    & 7781.88  \\ \hline
\end{tabular}
\begin{tablenotes}
\item Minimization of Rastrigin function with $d$=500 based on 10 replications. SMCO/-R/-BR run from 30 diagonal starting points in parallel, with 300 max iterations. 
\end{tablenotes}

\end{threeparttable}
\end{table}

\begin{table}
\caption{Sensitivity: Tolerance Level}\label{tab:tolerance}
\centering
\begin{tabular}{|c|c|ccc|ccc|}
\hline
\multicolumn{2}{|c|}{Test Function} & \multicolumn{6}{c|}{Rastrigin} \\ \hline
\multicolumn{2}{|c|}{Setting} 
& \multicolumn{3}{c|}{$d=500$, Max Iter = 300} 
& \multicolumn{3}{c|}{$d=500$, Max Iter = 500} \\ \hline
\multicolumn{2}{|c|}{Tolerance} 
& \multicolumn{1}{c|}{1e-04} & \multicolumn{1}{c|}{1e-06} & \multicolumn{1}{c|}{1e-08} 
& \multicolumn{1}{c|}{1e-04} & \multicolumn{1}{c|}{1e-06} & \multicolumn{1}{c|}{1e-08} \\ \hline

\multicolumn{1}{|c|}{\multirow{8}{*}{Max}} 
& \textbf{SMCO}    
& \textbf{233.20} & \textbf{222.78} & \textbf{206.89} 
& \textbf{202.98} & \textbf{226.73} & \textbf{195.98} \\

& \textbf{SMCO-R}  
& \textbf{183.33} & \textbf{197.71} & \textbf{188.02} 
& \textbf{192.89} & \textbf{192.06} & \textbf{206.55} \\

& \textbf{SMCO-BR} 
& \textbf{251.21} & \textbf{252.67} & \textbf{253.95} 
& \textbf{185.63} & \textbf{202.86} & \textbf{195.82} \\
\cline{3-8}

& GD       
& 5700.08 & 5700.91 & 5700.91 
& 5700.08 & 5700.91 & 5700.91 \\

& SignGD   
& 5720.35 & 5741.50 & 5741.50 
& 5720.35 & 5741.50 & 5741.50 \\

& SPSA    
& 5714.29 & 5714.29 & 5714.29 
& 5714.29 & 5714.29 & 5714.29 \\

& L\mbox{-}BFGS 
& 6790.24 & 4037.30 & 4037.30 
& 6984.23 & 4042.53 & 4042.53 \\

& BOBYQA    
& 10471.71 & 10497.07 & 10497.07 
& 10471.71 & 10497.07 & 10497.07 \\
\hline

\multicolumn{1}{|c|}{\multirow{9}{*}{Min}} 
& \textbf{SMCO}    
& \textbf{10.26} & \textbf{10.50} & \textbf{10.35} 
& \textbf{4.52} & \textbf{4.63} & \textbf{4.20} \\

& \textbf{SMCO-R}  
& \textbf{1.76} & \textbf{1.75} & \textbf{1.36} 
& \textbf{1.42} & \textbf{1.46} & \textbf{1.70} \\

& \textbf{SMCO-BR} 
& \textbf{1.60} & \textbf{1.60} & \textbf{1.76} 
& \textbf{1.57} & \textbf{1.20} & \textbf{1.76} \\
\cline{3-8}

& GD        
& 6328.91 & 6328.91 & 6328.91 
& 6328.91 & 6317.07 & 6328.91 \\

& SignGD    
& 3993.21 & 3993.21 & 3993.21 
& 3993.21 & 4022.50 & 3993.21 \\

& ADAM      
& 2077.66 & 987.57 & 987.57 
& 1961.95 & 981.28 & 981.27 \\

& SPSA      
& 10987.23 & 10987.23 & 10987.23 
& 10987.23 & 12860.98 & 10987.23 \\

& L\mbox{-}BFGS  
& 2164.18 & 1074.12 & 1074.12 
& 2048.47 & 1067.83 & 1067.83 \\

& BOBYQA    
& 8775.96 & 8775.96 & 8775.96 
& 8775.96 & 8852.37 & 8775.96 \\
\hline
\end{tabular}

Based on 10 reps. On original (untransformed) domain. All algorithms use 30 diagonal starting points. 
\end{table}

\begin{sidewaystable}
\centering
\caption{Sensitivity: starting values and varying iteration numbers, $d=400$, 1e-06 tolerance level. Results based on 10 replications.}
\label{tab:rastrigin400d_combined}
\begin{tabular}{|c|c|cc|cc|cc|cc|cc|cc|}
\hline
\multicolumn{14}{|c|}{Test Function: Original Rastrigin, Domain: $[-5.12, 5.12]^{400}$, Dimension = 400} \\ \hline
\multicolumn{2}{|c|}{Start Point} & \multicolumn{4}{c|}{Single Uniform Start Points} & \multicolumn{4}{c|}{20 Uniform Start Points} & \multicolumn{4}{c|}{20 Diagonal Start Points} \\ \hline
\multicolumn{2}{|c|}{Max Iteration} & \multicolumn{2}{c|}{300} & \multicolumn{2}{c|}{500} & \multicolumn{2}{c|}{300} & \multicolumn{2}{c|}{500} & \multicolumn{2}{c|}{300} & \multicolumn{2}{c|}{500} \\ \hline
\multicolumn{2}{|c|}{Performance} & Time & RMSE & Time & RMSE & Time & RMSE & Time & RMSE & Time & RMSE & Time & RMSE \\ \hline
\multirow{8}{*}{\rotatebox{90}{Max}} & \textbf{SMCO} & \textbf{3.15} & \textbf{268.37} & \textbf{5.34} & \textbf{240.11} & \textbf{5.74} & \textbf{183.15} & \textbf{9.29} & \textbf{156.48} & \textbf{5.71} & \textbf{23.41} & \textbf{9.28} & \textbf{8.43} \\
& \textbf{SMCO-R} & \textbf{3.16} & \textbf{229.00} & \textbf{5.28} & \textbf{225.45} & \textbf{5.87} & \textbf{165.02} & \textbf{9.44} & \textbf{155.62} & \textbf{5.80} & \textbf{1.51} & \textbf{9.48} & \textbf{1.30} \\
& \textbf{SMCO-BR} & \textbf{3.18} & \textbf{299.93} & \textbf{5.31} & \textbf{229.76} & \textbf{5.85} & \textbf{227.54} & \textbf{9.47} & \textbf{153.77} & \textbf{5.83} & \textbf{0.42} & \textbf{9.51} & \textbf{1.30} \\ \cline{2-14}
& GD & 0.36 & 4568.38 & 0.47 & 4570.12 & 7.28 & 4556.40 & 7.27 & 4556.40 & 0.25 & 863.47 & 0.31 & 863.47 \\
& SignGD & 11.18 & 4867.92 & 11.07 & 4786.97 & 18.36 & 4591.01 & 18.43 & 4591.01 & 14.82 & 0.00 & 14.86 & 0.00 \\
& SPSA & 0.01 & 5077.67 & 0.00 & 4571.43 & 0.29 & 4571.43 & 0.17 & 4393.22 & 0.20 & 4310.67 & 0.22 & 4478.01 \\
& L-BFGS & 0.05 & 4571.09 & 0.05 & 4571.09 & 0.35 & 4561.61 & 0.26 & 4561.61 & 0.21 & 0.00 & 0.26 & 0.00 \\
& BOBYQA & 0.75 & 7123.97 & 0.76 & 7115.03 & 3.70 & 6407.59 & 3.65 & 6407.59 & 3.63 & 0.27 & 3.69 & 0.27 \\ \hline
\multirow{9}{*}{\rotatebox{90}{Min}} & \textbf{SMCO} & \textbf{3.18} & \textbf{12.14} & \textbf{5.31} & \textbf{5.92} & \textbf{5.82} & \textbf{8.33} & \textbf{9.56} & \textbf{3.40} & \textbf{5.93} & \textbf{7.69} & \textbf{9.47} & \textbf{2.98} \\
& \textbf{SMCO-R} & \textbf{3.18} & \textbf{5.53} & \textbf{5.32} & \textbf{3.89} & \textbf{5.98} & \textbf{1.01} & \textbf{9.74} & \textbf{1.21} & \textbf{6.06} & \textbf{1.10} & \textbf{9.82} & \textbf{0.61} \\
& \textbf{SMCO-BR} & \textbf{3.21} & \textbf{5.11} & \textbf{5.38} & \textbf{4.04} & \textbf{5.98} & \textbf{1.27} & \textbf{9.79} & \textbf{1.29} & \textbf{6.03} & \textbf{0.47} & \textbf{9.67} & \textbf{0.44} \\ \cline{2-14}
& GD & 2.50 & 5663.71 & 3.85 & 5715.96 & 13.31 & 5089.62 & 13.37 & 5089.62 & 12.77 & 43.44 & 12.71 & 43.44 \\
& SignGD & 10.89 & 3492.43 & 10.87 & 3493.89 & 18.11 & 3227.54 & 18.07 & 3227.54 & 10.59 & 0.00 & 10.55 & 0.00 \\
& ADAM & 0.07 & 3690.40 & 0.07 & 3660.24 & 0.24 & 3369.64 & 0.30 & 3369.64 & 0.20 & 0.00 & 0.20 & 0.00 \\
& SPSA & 0.00 & 11635.37 & 0.01 & 10904.31 & 0.13 & 8422.62 & 0.17 & 8071.49 & 0.20 & 8375.80 & 0.21 & 7430.25 \\
& L-BFGS & 0.24 & 3051.31 & 0.23 & 3013.33 & 0.62 & 2588.30 & 0.61 & 2588.30 & 0.41 & 0.00 & 0.44 & 0.00 \\
& BOBYQA & 0.85 & 5757.97 & 0.85 & 5610.33 & 3.97 & 5205.29 & 4.00 & 5205.29 & 4.06 & 413.25 & 4.07 & 413.25 \\ \hline
\end{tabular}
\end{sidewaystable}

\begin{sidewaystable}
\caption{Sensitivity: bounds buffer values, $d=400$, 1e-06 tolerance level. Results based on 10 replications.}
\label{tab:rastrigin400d_combined1}
\small
\begin{tabular}{|c|c|c|cc|cc|cc|cc|cc|cc|}
\hline
\multicolumn{15}{|c|}{Test Function: Original Rastrigin, Domain: $[-5.12, 5.12]^{400}$, Dimension = 400} \\ \hline
\multicolumn{3}{|c|}{Start Point} & \multicolumn{4}{c|}{Single Uniform Start Points} & \multicolumn{4}{c|}{20 Uniform Start Points} & \multicolumn{4}{c|}{20 Diagonal Start Points} \\ \hline
\multicolumn{3}{|c|}{Max Iteration} & \multicolumn{2}{c|}{300} & \multicolumn{2}{c|}{500} & \multicolumn{2}{c|}{300} & \multicolumn{2}{c|}{500} & \multicolumn{2}{c|}{300} & \multicolumn{2}{c|}{500} \\ \hline
\multicolumn{3}{|c|}{Performance} & Time & RMSE & Time & RMSE & Time & RMSE & Time & RMSE & Time & RMSE & Time & RMSE \\ \hline
\multirow{9}{*}{\rotatebox{90}{Max}} & \multirow{3}{*}{$\d^*=0.05$}&  {SMCO}  &  {3.15} &  {268.37} &  {5.34} &  {240.11} &  {5.74} &  {183.15} &  {9.29} &  {156.48} &  {5.71} &  {23.41} &  {9.28} &  {8.43} \\
& &  {SMCO-R} &  {3.16} &  {229.00} &  {5.28} &  {225.45} &  {5.87} &  {165.02} &  {9.44} &  {155.62} &  {5.80} &  {1.51} &  {9.48} &  {1.30} \\
& & {SMCO-BR:} &  {3.18} &  {299.93} &  {5.31} &  {229.76} &  {5.85} &  {227.54} &  {9.47} &  {153.77} &  {5.83} &  {0.42} &  {9.51} &  {1.30} \\
\cline{2-15}
& \multirow{3}{*}{$\d^*=0.02$} & \textbf{SMCO} & \textbf{3.23} & \textbf{157.33} & \textbf{5.24} & \textbf{142.82} & \textbf{5.91} & \textbf{102.98} & \textbf{9.75} & \textbf{86.71} & \textbf{5.99} & \textbf{18.02} & \textbf{9.79} & \textbf{5.94} \\
& &\textbf{SMCO-R} & \textbf{3.22} & \textbf{135.19} & \textbf{5.24} & \textbf{130.87} & \textbf{6.15} & \textbf{69.81} & \textbf{9.74} & \textbf{74.62} & \textbf{6.07} & \textbf{1.47} & \textbf{9.85} & \textbf{1.28} \\
& &\textbf{SMCO-BR} & \textbf{3.23} & \textbf{161.84} & \textbf{5.27} & \textbf{125.82} & \textbf{5.98} & \textbf{101.87} & \textbf{9.81} & \textbf{73.22} & \textbf{6.10} & \textbf{1.49} & \textbf{9.98} & \textbf{1.35} \\ \hline
\multirow{9}{*}{\rotatebox{90}{Min}} &\multirow{3}{*}{$\d^*=0.05$} & {SMCO} &  {3.18} &  {12.14} &  {5.31} &  {5.92} &  {5.82} &  {8.33} &  {9.56} &  {3.40} &  {5.93} &  {7.69} &  {9.47} &  {2.98} \\
& & {SMCO-R} &  {3.18} &  {5.53} &  {5.32} &  {3.89} &  {5.98} &  {1.01} &  {9.74} &  {1.21} &  {6.06} &  {1.10} &  {9.82} &  {0.61} \\
& & {SMCO-BR} &  {3.21} &  {5.11} &  {5.38} &  {4.04} &  {5.98} &  {1.27} &  {9.79} &  {1.29} &  {6.03} &  {0.47} &  {9.67} &  {0.44} \\\cline{2-15}
& \multirow{3}{*}{$\d^*=0.02$}& \textbf{SMCO} & \textbf{3.23} & \textbf{7.59} & \textbf{5.30} & \textbf{2.74} & \textbf{6.06} & \textbf{7.03} & \textbf{9.89} & \textbf{2.53} & \textbf{6.15} & \textbf{3.54} & \textbf{10.19} & \textbf{1.81} \\
& &\textbf{SMCO-R} & \textbf{3.25} & \textbf{0.48} & \textbf{5.30} & \textbf{0.43} & \textbf{6.21} & \textbf{0.46} & \textbf{10.11} & \textbf{0.40} & \textbf{6.26} & \textbf{0.47} & \textbf{10.07} & \textbf{0.40} \\
& &\textbf{SMCO-BR} & \textbf{3.27} & \textbf{0.58} & \textbf{5.34} & \textbf{0.52} & \textbf{6.28} & \textbf{0.54} & \textbf{10.04} & \textbf{0.49} & \textbf{6.26} & \textbf{0.47} & \textbf{10.23} & \textbf{0.43} \\ \cline{2-14} \hline

\hline
\end{tabular}
\end{sidewaystable}

\begin{table}
\caption{Sensitivity: Running Time of Selected Competing Algorithms, Original Rastrigin with $d=1000$}
\label{tab:rastrigin1000_hd_combined}
\centering
\begin{threeparttable}
\begin{tabular}{|l|c|c|c|c|c|c|c|c|}
\hline
\multirow{3}{*}{Algorithm} & \multicolumn{4}{c|}{\textbf{Minimization (Target: 0)}} & \multicolumn{4}{c|}{\textbf{Maximization  (Target: 40353.29) }} \\
\cline{2-9}
& \multirow{2}{*}{\begin{tabular}[c]{@{}c@{}}Time\\(sec)\end{tabular}} & \multirow{2}{*}{RMSE} & \multirow{2}{*}{AE50} & \multirow{2}{*}{AE95} & \multirow{2}{*}{\begin{tabular}[c]{@{}c@{}}Time\\(sec)\end{tabular}} & \multirow{2}{*}{RMSE} & \multirow{2}{*}{AE50} & \multirow{2}{*}{AE95} \\
& & & & & & & & \\
\hline
\multicolumn{9}{|l|}{ Original Rastrigin Domain: $[-5.12,5.12]^{1000}$.} \\ \hline
\textbf{SMCO}    & \textbf{56.14} & \textbf{19.62} & \textbf{19.47} & \textbf{20.75} & \textbf{54.94} & \textbf{59.65} & \textbf{59.63} & \textbf{61.11} \\ \hline
\textbf{SMCO-R}  & \textbf{56.19} & \textbf{0.00}  & \textbf{0.00}  & \textbf{0.00}  & \textbf{54.98} & \textbf{0.00}  & \textbf{0.00}  & \textbf{0.00}  \\ \hline
\textbf{SMCO-BR} & \textbf{56.70} & \textbf{0.00}  & \textbf{0.00}  & \textbf{0.00}  & \textbf{55.44} & \textbf{0.00}  & \textbf{0.00}  & \textbf{0.00}  \\ \hline
GenSA            & 1867.35 & 0.00 & 0.00 & 0.00 & 2232.68 & 0.00 & 0.00 & 0.00 \\ \hline
GenSA*           & 57.78  & 387.64 & 379.08 & 432.56 & 61.14 & 50.22 & 27.14 & 85.64 \\ \hline
DEoptim          & 113.85 & 15847.66 & 15859.08 & 15915.79 & 115.16 & 19385.83 & 19387.06 & 19491.50 \\ \hline
PSO              & 5.64 & 5417.89 & 6048.48 & 6549.47 & 6.03 & 10959.24 & 11064.37 & 11418.06 \\ \hline
    PSO*              & 58.36 & 5409.06 & 5503.365 & 6537.48 & 61.72 & 9459.40 & 10552.61 & 11444.92 \\ \hline
\end{tabular}
\begin{tablenotes}
\item Results based on 10 replications with diagonal draw starting points. SMCO configuration: 35 starting points, 300 max iterations, parallel computing enabled with tolerance level 1e-06. 
RMSE: root mean square error; AE: absolute error percentiles. SMCO-R, SMCO-BR and
GenSA achieve perfect convergence (0.00) for both optimization directions, but GenSA needs a great deal of time.  GenSA and PSO use the defalt settings, relatively faster GenSA* uses the same defalt hyperparameters of GenSA except for
max.call= 1e6 (defalt max.call=1e7) and PSO* choose the same hyperparameters of PSO  except for $s$=700 (defalt $s$=10+2$\sqrt{d}$). 
\end{tablenotes}
\end{threeparttable}
\end{table}

\begin{table}
\begin{centering}
\caption{\label{tab:CrossLeg} Minimization of Cross-Leg-Table Function $\left(d=2\right)$}
\begin{tabular}{|ccccccc|}
\hline 
Starting Points & \multicolumn{3}{c|}{3} & \multicolumn{3}{c|}{25}\tabularnewline
\hline 
Performance & \multicolumn{1}{c|}{RMSE} & \multicolumn{1}{c|}{AE99} & \multicolumn{1}{c|}{Time} & \multicolumn{1}{c|}{RMSE} & \multicolumn{1}{c|}{AE99} & Time\tabularnewline
\cline{2-7} \cline{3-7} \cline{4-7} \cline{5-7} \cline{6-7} \cline{7-7} 
\textbf{Truth} & \multicolumn{3}{c}{\textbf{-1}} & \multicolumn{3}{c|}{\textbf{-1}}\tabularnewline
\cline{2-7} \cline{3-7} \cline{4-7} \cline{5-7} \cline{6-7} \cline{7-7} 
GD & 0.9999  & 0.9999  & 0.079  & 0.9998  & 0.9999  & 0.322 \tabularnewline
SignGD & 0.9989  & 0.9994  & 0.096  & 0.9987  & 0.9989  & 0.431 \tabularnewline
ADAM & 0.9563  & 0.9979  & 0.006  & 0.4455  & 0.9964  & 0.042 \tabularnewline
SPSA & 0.9999  & 0.9999  & 0.018  & 0.9999  & 0.9999  & 0.014 \tabularnewline
L-BFGS & 0.9997  & 0.9999  & 0.001  & 0.9995  & 0.9996  & 0.006 \tabularnewline
BOBYQA & 0.9995  & 0.9997  & 0.004  & 0.9993  & 0.9994  & 0.022 \tabularnewline
\cline{2-7} \cline{3-7} \cline{4-7} \cline{5-7} \cline{6-7} \cline{7-7} 
\textbf{SMCO} & 0.9998  & 0.9998  & 0.027  & 0.9997  & 0.9998  & 0.158 \tabularnewline
\textbf{SMCO-R} & 0.9995  & 0.9996  & 0.029  & 0.9993  & 0.9995  & 0.163 \tabularnewline
\textbf{SMCO-BR} & 0.9995  & 0.9996  & 0.028  & 0.9994  & 0.9995  & 0.166 \tabularnewline
\cline{2-7} \cline{3-7} \cline{4-7} \cline{5-7} \cline{6-7} \cline{7-7} 
GenSA & 0.3303  & 0.9965  & 2.372  & \multicolumn{3}{c|}{$/$}\tabularnewline
SA & 0.9998  & 0.9998  & 0.007  & \multicolumn{3}{c|}{$/$}\tabularnewline
DEoptim & 0.9996  & 0.9998  & 0.017  & \multicolumn{3}{c|}{$/$}\tabularnewline
CMAES & - & - & - & \multicolumn{3}{c|}{$/$}\tabularnewline
STOGO & 0.9997  & 0.9997  & 0.321  & \multicolumn{3}{c|}{$/$}\tabularnewline
GA & 0.9990  & 0.9998  & 0.093  & \multicolumn{3}{c|}{$/$}\tabularnewline
PSO & 0.9951 & 0.9992  & 0.188 & \multicolumn{3}{c|}{$/$}\tabularnewline
\hline 
\end{tabular}
\par\end{centering}
\centering{}{\small Based on 100 replications. Tolerance=1e-08; Max iteration = 400. Uniform start. ``-'' indicates
algorithm termination with error.}
\end{table}

\section{Supplemental Numerical Results on Maximum Score Estimation and Empirical Welfare Maximization} \label{append-E}

~

\textbf{Table \ref{tab:MSEW_rep10}} replicates Table \ref{tab:MSEW} on the two random test functions (maximum score estimation and empirical welfare maximization) \textbf{using 10 replications only} (instead of $500/d = 100,50,25,10$ replications as in Table \ref{tab:MSEW} for $d=5,10,20,50$). This replication serves to show that the results in Table \ref{tab:MSEW_rep10} obtained using only 10 replications are very similar to those in Table \ref{tab:MSEW} based on more replications. Due to constraints on computation time and resources, we report most results on higher dimensions in the main text and the appendix based on 10 replications. While more replications generally yield more accurate performance measures given the randomness in the optimization algorithms, we hope that the similarity of Table \ref{tab:MSEW}  and Table \ref{tab:MSEW_rep10} suggests that our numerical results are fairly stable with respect to the number of replications used.

\textbf{Table \ref{tab:MSEW_diag}} replicates Table \ref{tab:MSEW} using diagonal starting points instead of uniform ones for SMCO/-R/-BR and Group I (local) algorithms. While Group I algorithms display sensitivity with respect to the type of starting points, the performance of our SMCO algorithms are very stable, still outperforming both Group I and II algorithms. 

\textbf{Table \ref{tab:MS_Cauchy}} reports results on maximum score estimation for $d=50,100$, in which the error term $\e_i$ is drawn from the Cauchy distribution instead of the normal distribution in the main text. \textbf{Table \ref{tab:MS_Cauchy_buffer}} further reports the sensitivity checks of SMCO performance in maximum score estimation in $d=100$ with respect to the SMCO bound buffer parameter $\delta^*$. Overall our SMCO algorithms perform well across these test configurations.

\textbf{Table \ref{tab:EW-100d}} reports results on empirical welfare maximization in $d=100$, as well as the minimization of empirical welfare function as an additional test function (without economic meanings). The results further demonstrate the best performance of our SMCO algorithsm among all Group I and II algorithms in higher dimensions, as well as the robustness of our SMCO algorithms wrt the SMCO bound buffer parameter $\delta^*$.

\begin{sidewaystable}
\begin{centering}
\caption{\label{tab:MSEW_rep10} Maximization of $f_{MS}\left(\protect\t\right)$
and $f_{EW}\left(\protect\t\right)$}
\begin{tabular}{|ccccc|ccc|ccc|ccc|}
\hline 
\multicolumn{2}{|c}{Dimension} & \multicolumn{3}{c|}{\textbf{$d=5$}} & \multicolumn{3}{c|}{\textbf{$d=10$}} & \multicolumn{3}{c|}{\textbf{$d=20$}} & \multicolumn{3}{c|}{\textbf{$d=50$}}\tabularnewline
\hline 
\multicolumn{2}{|c}{Performance} & \multicolumn{1}{c|}{RMSE} & \multicolumn{1}{c|}{AE99} & Time & \multicolumn{1}{c|}{RMSE} & \multicolumn{1}{c|}{AE99} & Time & \multicolumn{1}{c|}{RMSE} & \multicolumn{1}{c|}{AE99} & Time & \multicolumn{1}{c|}{RMSE} & \multicolumn{1}{c|}{AE99} & Time\tabularnewline
\hline 
\multicolumn{14}{|c|}{\textbf{Maximum Score Estimation: $\max_{\b}f_{MS}\left(\b\right)\phantom{\frac{\frac{1}{1}}{\frac{1}{2}}}$}}\tabularnewline
\hline 
 & \textbf{Best} & \multicolumn{3}{c|}{\textbf{0.42 }} & \multicolumn{3}{c|}{\textbf{0.448 }} & \multicolumn{3}{c|}{\textbf{0.462 }} & \multicolumn{3}{c|}{\textbf{0.424 }}\tabularnewline
\hline 
\multirow{6}{*}{(I)} & GD & 0.137  & 0.169  & 0.06  & 0.154  & 0.183  & 0.03  & 0.172  & 0.197  & 0.08  & 0.157  & 0.174  & 0.13 \tabularnewline
 & SignGD & 0.137  & 0.169  & 0.08  & 0.154  & 0.183  & 0.03  & 0.172  & 0.197  & 0.09  & 0.157  & 0.174  & 0.11 \tabularnewline
 & ADAM & 0.137  & 0.169  & 0.07  & 0.154  & 0.183  & 0.07  & 0.172  & 0.197  & 0.08  & 0.157  & 0.174  & 0.11 \tabularnewline
 & SPSA & 0.137  & 0.169  & 0.11  & 0.154  & 0.183  & 0.08  & 0.172  & 0.197  & 0.12  & 0.157  & 0.174  & 0.12 \tabularnewline
 & L-BFGS & 0.137  & 0.169  & 0.07  & 0.154  & 0.183  & 0.08  & 0.172  & 0.197  & 0.10  & 0.157  & 0.174  & 0.13 \tabularnewline
 & BOBYQA & 0.090  & 0.105  & 0.07  & 0.110  & 0.136  & 0.09  & 0.131  & 0.146  & 0.11  & 0.134  & 0.155  & 0.22 \tabularnewline
\hline 
\multirow{3}{*}{} & \textbf{SMCO} & 0.008  & 0.010  & 0.26  & 0.009  & 0.012  & 0.40  & 0.013  & 0.016  & 0.86  & 0.056  & 0.068  & 5.09 \tabularnewline
 & \textbf{SMCO-R} & 0.002  & 0.004  & 0.22  & 0.004  & 0.006  & 0.37  & 0.004  & 0.006  & 0.86  & 0.010  & 0.016  & 5.19 \tabularnewline
 & \textbf{SMCO-BR} & 0.003  & 0.004  & 0.24  & 0.007  & 0.008  & 0.40  & 0.010  & 0.012  & 0.89  & 0.010  & 0.017  & 5.41 \tabularnewline
\hline 
\multirow{6}{*}{(II)} & GenSA & 0.017  & 0.054  & 1.96  & 0.032  & 0.034  & 5.01  & 0.032  & 0.036  & 16.74  & 0.010  & 0.018  & 89.93 \tabularnewline
 & SA & 0.071  & 0.076  & 0.07  & 0.073  & 0.087  & 0.08  & 0.106  & 0.116  & 0.14  & 0.114  & 0.124  & 0.34 \tabularnewline
 & DEoptim & 0.004  & 0.006  & 0.34  & 0.036  & 0.038  & 0.82  & 0.045  & 0.048  & 2.92  & 0.056  & 0.064  & 13.18 \tabularnewline
 & STOGO & 0.052  & 0.052  & 3.88  & 0.068  & 0.068  & 8.84  & 0.080  & 0.080  & 29.02  & 0.082  & 0.082  & 182.28 \tabularnewline
 & GA & 0.053  & 0.058  & 0.27  & 0.035  & 0.040  & 0.29  & 0.053  & 0.062  & 0.43  & 0.059  & 0.070  & 0.65 \tabularnewline
 & PSO & 0.048  & 0.054  & 0.84  & 0.034  & 0.038  & 1.11  & 0.039  & 0.044  & 1.80  & 0.017  & 0.026  & 4.32 \tabularnewline
\hline 
\multicolumn{14}{|c|}{\textbf{Empirical Welfare Maximization: $\max_b f_{EW}\left(\b\right)\phantom{\frac{\frac{1}{1}}{\frac{1}{2}}}$}}\tabularnewline
\hline 
 & \textbf{Best} & \multicolumn{3}{c|}{\textbf{1541.032 }} & \multicolumn{3}{c|}{\textbf{2061.409 }} & \multicolumn{3}{c|}{\textbf{2634.523 }} & \multicolumn{3}{c|}{\textbf{3825.49 }}\tabularnewline
\hline 
\multirow{6}{*}{(I)} & GD & 1058  & 1442  & 0.11  & 1351  & 1646  & 0.14  & 1767  & 1942  & 0.22  & 2934  & 3067  & 0.60 \tabularnewline
 & SignGD & 1058  & 1442  & 0.14  & 1351  & 1646  & 0.14  & 1767  & 1942  & 0.21  & 2934  & 3067  & 0.60 \tabularnewline
 & ADAM & 1058  & 1442  & 0.13  & 1351  & 1646  & 0.15  & 1767  & 1942  & 0.21  & 2934  & 3067  & 0.62 \tabularnewline
 & SPSA & 991  & 1442  & 0.12  & 1344  & 1721  & 0.13  & 1788  & 1942  & 0.17  & 2920  & 3063  & 0.69 \tabularnewline
 & L-BFGS & 886  & 1437  & 0.12  & 1259  & 1448  & 0.16  & 1766  & 1942  & 0.23  & 2927  & 3067  & 0.99 \tabularnewline
 & BOBYQA & 465  & 921  & 0.17  & 930  & 1235  & 0.20  & 1272  & 1503  & 0.35  & 2338  & 2599  & 1.51 \tabularnewline
\hline 
\multirow{3}{*}{} & \textbf{SMCO} & 214  & 270  & 1.81  & 191  & 309  & 3.88  & 249  & 312  & 10.41  & 428  & 498  & 78.37 \tabularnewline
 & \textbf{SMCO-R} & 139  & 218  & 1.73  & 69  & 131 & 3.85  & 82  & 159  & 10.46  & 158  & 233  & 78.59 \tabularnewline
 & \textbf{SMCO-BR} & 134  & 208  & 1.75  & 52  & 87  & 3.88  & 124  & 193  & 10.52  & 274  & 361  & 78.71 \tabularnewline
\hline 
\multirow{6}{*}{(II)} & GenSA & 105  & 145  & 15.15  & 191  & 245  & 35.40  & 636  & 903  & 102.96  & 1586  & 1762  & 455.24 \tabularnewline
 & SA & 199  & 247  & 0.51  & 416  & 527  & 0.63  & 860  & 1227  & 0.94  & 1939  & 2398  & 1.57 \tabularnewline
 & DEoptim & 136  & 155  & 2.69  & 232  & 272  & 6.61 & 611  & 681  & 19.74  & 1722  & 1805  & 87.03 \tabularnewline
 & STOGO & 134  & 134  & 31.85  & 470  & 470  & 71.97  & 835  & 835  & 204.88  & 2096  & 2096  & 891.62 \tabularnewline
 & GA & 155  & 206  & 1.21  & 277  & 441  & 1.48  & 657  & 779  & 2.17  & 1685  & 1831  & 3.69 \tabularnewline
 & PSO & 102  & 134  & 4.13  & 165  & 274  & 5.72  & 285  & 384  & 9.46  & 848  & 1051  & 21.52 \tabularnewline
\hline 
\end{tabular}
\par\end{centering}
\begin{centering}
SMCO/-R/-BR and Group (I) algos use $\sqrt{d}\approx2,3,4,7$ uniform
starting points.
\par\end{centering}
\centering{}All based on 10 replications. 
\end{sidewaystable}

\begin{sidewaystable}
\begin{centering}
\caption{\label{tab:MSEW_diag} Maximum Score Estimation and Empirical Welfare Maximization: Diagonal Multi-Start}
\begin{tabular}{|ccccc|ccc|ccc|ccc|}
\hline 
\multicolumn{2}{|c}{Dimension} & \multicolumn{3}{c|}{\textbf{$d=5$}} & \multicolumn{3}{c|}{\textbf{$d=10$}} & \multicolumn{3}{c|}{\textbf{$d=20$}} & \multicolumn{3}{c|}{\textbf{$d=50$}}\tabularnewline
\hline 
\multicolumn{2}{|c}{Performance} & \multicolumn{1}{c|}{RMSE} & \multicolumn{1}{c|}{AE99} & Time & \multicolumn{1}{c|}{RMSE} & \multicolumn{1}{c|}{AE99} & Time & \multicolumn{1}{c|}{RMSE} & \multicolumn{1}{c|}{AE99} & Time & \multicolumn{1}{c|}{RMSE} & \multicolumn{1}{c|}{AE99} & Time\tabularnewline
\hline 
\multicolumn{14}{|c|}{\textbf{Maximum Score Estimation: $\max_{\b}f_{MS}\left(\b\right)\phantom{\frac{\frac{1}{1}}{\frac{1}{2}}}$}}\tabularnewline
\hline 
 & \textbf{Best} & \multicolumn{3}{c|}{\textbf{0.42}} & \multicolumn{3}{c|}{\textbf{0.45}} & \multicolumn{3}{c|}{\textbf{0.46}} & \multicolumn{3}{c|}{\textbf{0.424}}\tabularnewline
\hline 
\multirow{6}{*}{(I)} & GD & 0.060  & 0.060  & 0.08  & 0.090  & 0.090  & 0.04  & 0.152  & 0.152  & 0.09  & 0.098  & 0.098  & 0.13 \tabularnewline
 & SignGD & 0.060  & 0.060  & 0.07  & 0.090  & 0.090  & 0.04  & 0.152  & 0.152  & 0.08  & 0.098  & 0.098  & 0.11 \tabularnewline
 & ADAM & 0.060  & 0.060  & 0.07  & 0.090  & 0.090  & 0.07  & 0.152  & 0.152  & 0.08  & 0.098  & 0.098  & 0.11 \tabularnewline
 & SPSA & 0.060  & 0.060  & 0.11  & 0.090  & 0.090  & 0.08  & 0.152  & 0.152  & 0.12  & 0.098  & 0.098  & 0.12 \tabularnewline
 & L-BFGS & 0.060  & 0.060  & 0.07 & 0.090  & 0.090  & 0.07  & 0.152  & 0.152  & 0.09  & 0.098  & 0.098  & 0.11 \tabularnewline
 & BOBYQA & 0.022  & 0.022  & 0.07 & 0.020  & 0.020  & 0.08  & 0.066  & 0.066  & 0.10  & 0.022  & 0.022  & 0.24 \tabularnewline
\hline 
\multirow{3}{*}{} & \textbf{SMCO} & 0.008 & 0.010 & 0.29 & 0.010 & 0.014 & 0.41 & 0.012 & 0.018 & 0.94 & 0.052 & 0.058 & 4.92 \tabularnewline
 & \textbf{SMCO-R} & 0.003 & 0.004 & 0.25 & 0.007 & 0.010 & 0.38 & 0.003 & 0.006 & 0.92 & 0.008 & 0.012 & 4.99 \tabularnewline
 & \textbf{SMCO-BR} & 0.002 & 0.004 & 0.28 & 0.007 & 0.010 & 0.41 & 0.008 & 0.012 & 0.97 & 0.009 & 0.014 & 5.16\tabularnewline
\hline 
\multirow{6}{*}{(II)} & GenSA & 0.017 & 0.054 & 1.96 & 0.034 & 0.036 & 5.01 & 0.030 & 0.034 & 16.74 & 0.010 & 0.018 & 89.93\tabularnewline
 & SA & 0.071 & 0.076 & 0.07 & 0.075 & 0.089 & 0.08 & 0.104 & 0.114 & 0.14 & 0.114 & 0.124 & 0.34\tabularnewline
 & DEoptim & 0.004 & 0.006 & 0.34 & 0.038 & 0.040 & 0.82 & 0.043 & 0.046 & 2.92 & 0.056 & 0.064 & 13.18\tabularnewline
 & STOGO & 0.052 & 0.052 & 3.88 & 0.070 & 0.070 & 8.84 & 0.078 & 0.078 & 29.02 & 0.082 & 0.082 & 182.28\tabularnewline
 & GA & 0.053 & 0.058 & 0.27 & 0.037 & 0.042 & 0.29 & 0.051 & 0.060 & 0.43 & 0.059 & 0.070 & 0.65\tabularnewline
 & PSO & 0.048 & 0.054 & 0.84 & 0.036 & 0.040 & 1.11 & 0.037 & 0.042 & 1.80 & 0.017 & 0.026 & 4.32\tabularnewline
\hline 
\multicolumn{14}{|c|}{\textbf{Empirical Welfare Maximization: $\max_b  f_{EW}\left(\b\right)\phantom{\frac{\frac{1}{1}}{\frac{1}{2}}}$}}\tabularnewline
\hline 
 & \textbf{Best} & \multicolumn{3}{c|}{\textbf{1541.032}} & \multicolumn{3}{c|}{\textbf{2050.312 }} & \multicolumn{3}{c|}{\textbf{2659.651 }} & \multicolumn{3}{c|}{\textbf{3827.822 }}\tabularnewline
\hline 
\multirow{6}{*}{(I)} & GD & 356  & 356  & 0.15  & 898  & 898  & 0.19  & 2122  & 2122  & 0.24  & 2477  & 2477  & 0.70 \tabularnewline
 & SignGD & 356  & 356  & 0.13  & 898  & 898  & 0.16  & 2122  & 2122  & 0.22  & 2477  & 2477  & 0.67 \tabularnewline
 & ADAM & 484  & 484  & 0.45  & 993  & 993  & 0.56  & 2122  & 2122  & 0.21  & 2771  & 2771  & 1.41 \tabularnewline
 & SPSA & 622  & 857  & 0.14  & 1329  & 1489 & 0.17  & 2094  & 2122  & 0.17  & 2974  & 3153  & 0.73 \tabularnewline
 & L-BFGS & 236  & 236  & 0.43  & 797 & 797 & 0.56  & 2122  & 2122  & 0.20  & 1775  & 1775  & 8.83 \tabularnewline
 & BOBYQA & 386  & 386  & 0.19  & 656  & 656  & 0.25  & 1085   & 1085   & 0.36  & 1163  & 1163  & 1.50 \tabularnewline
\hline 
\multirow{3}{*}{} & \textbf{SMCO} & 200 & 268  & 1.89  & 164  & 200  & 4.27  & 318  & 386  & 10.32 & 409 & 547  & 76.46 \tabularnewline
 & \textbf{SMCO-R} & 136 & 153  & 1.84  & 63 & 109  & 4.24  & 106  & 173  & 10.33  & 172  & 251  & 76.56 \tabularnewline
 & \textbf{SMCO-BR} & 146 & 177  & 1.91  & 66  & 120  & 4.27  & 118  & 157  & 10.42  & 217  & 316  & 76.68 \tabularnewline
\hline 
\multirow{6}{*}{(II)} & GenSA & 105 & 145 & 15.15 & 181 & 234 & 35.40 & 660 & 928  & 102.96 & 1589 & 1765  & 455.24\tabularnewline
 & SA & 199 & 247 & 0.51 & 405 & 516 & 0.63 & 885 & 1252 & 0.94 & 1942 & 2401  & 1.57\tabularnewline
 & DEoptim & 136 & 155 & 2.69 & 221 & 261 & 6.61 & 637 & 706  & 19.74 & 1724 & 1807  & 87.03\tabularnewline
 & STOGO & 134 & 134 & 31.85 & 459 & 459 & 71.97 & 860 & 860  & 204.88 & 2098 & 2098  & 891.62\tabularnewline
 & GA & 155 & 206 & 1.21 & 267 & 430 & 1.48 & 682 & 804  & 2.17 & 1687 & 1833  & 3.69\tabularnewline
 & PSO & 102 & 134 & 4.13 & 155 & 262 & 5.72 & 309  & 409  & 9.46 & 851 & 1054  & 21.52\tabularnewline
\hline 
\end{tabular}
\par\end{centering}
\centering{}Based on 10 replications. SMCO/-R/-BR and Group (I) algos
use $\sqrt{d}\approx2,3,4,7$ diagonal starting points.
\end{sidewaystable}

\begin{table}
\caption{\label{tab:MS_Cauchy}Maximum Score Estimation with Cauchy Errors}
\centering
\scalebox{1}{
\begin{tabular}{|c|ccc|ccc|}
\hline 
\multicolumn{1}{|c|}{} & \multicolumn{3}{c|}{$d = 50$} & \multicolumn{3}{c|}{$d = 100$} \\ 
\hline 
\multicolumn{1}{|c|}{Algorithm} & RMSE & AE99 & Time & RMSE & AE99 & Time \\ 
\hline 

\textbf{Best} 
& \multicolumn{3}{c|}{\textbf{0.4900}} 
& \multicolumn{3}{c|}{\textbf{0.4540}} \\ 
\cline{2-7}

\textbf{SMCO}
& \textbf{0.1118} & \textbf{0.1176} & 19.8548
& \textbf{0.1410} & \textbf{0.1458} & 70.1807 \\

\textbf{SMCO-R}
& \textbf{0.0210} & \textbf{0.0276} & 19.2028
& \textbf{0.0341} & \textbf{0.0398} & 74.7200 \\

\textbf{SMCO-BR}
& \textbf{0.0046} & \textbf{0.0060} & 19.6607
& \textbf{0.0143} & \textbf{0.0180} & 72.0697 \\

\cline{2-7}

GD
& 0.1420 & 0.1420 & 0.3721
& 0.1020 & 0.1020 & 1.9571 \\

SignGD
& 0.1420 & 0.1420 & 0.3554
& 0.1020 & 0.1020 & 1.8590 \\

SPSA
& 0.1420 & 0.1420 & 0.0506
& 0.1020 & 0.1020 & 0.0853 \\

optimLBFGS
& 0.1420 & 0.1420 & 0.3520
& 0.1020 & 0.1020 & 2.8355 \\

BOBYQA
& 0.0420 & 0.0420 & 2.2335
& 0.0240 & 0.0240 & 16.2511 \\

GenSA
& 0.0344 & 0.0380 & 63.6950
& 0.0152 & 0.0218 & 286.4138 \\

SA
& 0.1520 & 0.1600 & 0.2773
& 0.1429 & 0.1553 & 0.5721 \\

DEoptim
& 0.0841 & 0.0898 & 10.7332
& 0.0922 & 0.0998 & 44.2958 \\

CMAES
& 0.0581 & 0.0789 & 3.4109
& 0.0315 & 0.0478 & 13.1697 \\

stogo
& 0.1060 & 0.1060 & 136.1108
& 0.0880 & 0.0880 & 541.6599 \\

GA
& 0.0874 & 0.1206 & 0.6455
& 0.0872 & 0.1095 & 1.3692 \\

PSO
& 0.0452 & 0.0589 & 3.7318
& 0.0357 & 0.0526 & 9.0746 \\

\hline
\end{tabular}}

Based on 10 replications. SMCO/-R/-BR and Group I algorithms run from $4\sqrt{d}$ diagonal starting points. SMCO: iter max = 1000, bounds buffer = 0.02, parallelization across start points.
\end{table}

\begin{table}
\caption{\label{tab:MS_Cauchy_buffer} Maximum Score Estimation with Cauchy Errors and $d=100$: Sensitivity wrt Bound Buffer $\delta^*$}
\centering
\scalebox{1}{
\begin{tabular}{|c|ccc|ccc|}
\hline 
\multicolumn{1}{|c|}{Bound Buffer} & \multicolumn{3}{c|}{$\delta^* = 0.02$} & \multicolumn{3}{c|}{$\delta^* = 0.05$} \\
\hline 
\hline 
\multicolumn{1}{|c|}{Algorithm} & RMSE & AE99 & Time (s) & RMSE & AE99 & Time (s)\\
\hline 

\textbf{Best} 
& \multicolumn{3}{c|}{\textbf{0.4640}} & \multicolumn{3}{c|}{\textbf{0.458}} \\  
\cline{2-7}

\textbf{SMCO}    
& \textbf{0.1216} & \textbf{0.1280} & 46.6113 & \textbf{0.1154} & \textbf{0.1240} & 45.4532 \\

\textbf{SMCO-R}  
& \textbf{0.0457} & \textbf{0.0514} & 45.1869 & \textbf{0.0388} & \textbf{0.0468} & 44.4921 \\

\textbf{SMCO-BR} 
& \textbf{0.0302} & \textbf{0.0354} & 48.1126 & \textbf{0.0236} & \textbf{0.0314} & 46.9663 \\

\cline{2-7}

GD         
& 0.1120 & 0.1120 & 1.8369 & 0.1060 & 0.1060 & 1.8720 \\

SignGD     
& 0.1120 & 0.1120 & 1.8657 & 0.1060 & 0.1060 & 1.8777 \\

SPSA       
& 0.1120 & 0.1120 & 0.0865 & 0.1060 & 0.1060 & 0.0931 \\

optimLBFGS 
& 0.1120 & 0.1120 & 3.0625 & 0.1060 & 0.1060 & 3.3938 \\

BOBYQA    
& 0.0340 & 0.0340 & 16.3390 & 0.0280 & 0.0280 & 16.6256 \\

GenSA        
& 0.0255 & 0.0314 & 281.1997 & 0.0185 & 0.0308 & 308.4736 \\

SA           
& 0.1521 & 0.1660 & 0.5504 & 0.1476 & 0.1580 & 0.6533 \\

DEoptim      
& 0.1042 & 0.1100 & 45.6777 & 0.0976 & 0.1054 & 47.3026 \\

CMAES        
& 0.0373 & 0.0588 & 14.2327 & 0.0374 & 0.0771 & 13.0083 \\

stogo        
& 0.0980 & 0.0980 & 565.0344    
& 0.0920 & 0.0920 & 616.7450 \\

GA           
& 0.0903 & 0.1063 & 1.2363 & 0.0882 & 0.1100 & 1.3138 \\

PSO          
& 0.0420 & 0.0574 & 8.9114 & 0.0382 & 0.0588 & 9.6422 \\

\hline
\end{tabular}}

Based on $30$ replications. SMCO and Group I algos run from diagonal start $4\sqrt{d}$. 

SMCO run with max iter = 500, parallelization across start points.
\end{table}

\begin{table}
\begin{centering}
\caption{\label{tab:EW-100d} Empirical Welfare Maximization \& Minimization:\textbf{
$d=100$}}
\begin{tabular}{|cc|cccccc|}
\hline 
\multicolumn{2}{|c|}{\# Starting Points} & \multicolumn{3}{c|}{Max} & \multicolumn{3}{c|}{Min}\tabularnewline
\hline 
\multicolumn{2}{|c|}{Performance} & \multicolumn{1}{c|}{RMSE} & \multicolumn{1}{c|}{AE99} & \multicolumn{1}{c|}{Time} & \multicolumn{1}{c|}{RMSE} & \multicolumn{1}{c|}{AE99} & Time\tabularnewline
\hline 
 & \textbf{Best} & \multicolumn{3}{c|}{\textbf{4536.174}} & \multicolumn{3}{c|}{\textbf{-3513.749}}\tabularnewline
\hline 
\multirow{6}{*}{(I)} & GD & 3589  & 3815 & 1.5 & 3732  & 3887  & 1.6 \tabularnewline
 & SignGD & 3589  & 3815 & 1.3 & 3732  & 3887  & 1.5 \tabularnewline
 & ADAM & 3589  & 3815 & 1.5 & 3732  & 3887  & 1.8 \tabularnewline
 & SPSA & 3616  & 3823 & 0.6 & 3709  & 3891  & 0.6 \tabularnewline
 & L-BFGS & 3589  & 3815 & 3.6 & 3701  & 3887  & 4.6 \tabularnewline
 & BOBYQA & 3062  & 3279 & 4.8 & 3177  & 3361  & 6.1 \tabularnewline
\hline 
\multirow{3}{*}{$\d^{*}=0.05$} & \textbf{SMCO} & \textbf{328 } & \textbf{406 } & \textbf{362} & \textbf{328 } & \textbf{455 } & \textbf{461 }\tabularnewline
 & \textbf{SMCO-R} & \textbf{92 } & \textbf{174 } & \textbf{352} & \textbf{99 } & \textbf{137 } & \textbf{475 }\tabularnewline
 & \textbf{SMCO-BR} & \textbf{267 } & \textbf{324 } & \textbf{351} & \textbf{311 } & \textbf{382 } & \textbf{477 }\tabularnewline\hline
\multirow{3}{*}{$\d^{*}=0.02$} & \textbf{SMCO} & \textbf{349 } & \textbf{404 } & \textbf{413 } & \textbf{386 } & \textbf{485 } & \textbf{421 }\tabularnewline
 & \textbf{SMCO-R} & \textbf{70 } & \textbf{142 } & \textbf{421 } & \textbf{132 } & \textbf{207 } & \textbf{420 }\tabularnewline
 & \textbf{SMCO-BR} & \textbf{294 } & \textbf{343 } & \textbf{411 } & \textbf{311 } & \textbf{421 } & \textbf{424 }\tabularnewline
\hline 
\multirow{5}{*}{(II)} & GenSA & 2079  & 2463  & 1478 & 2137  & 2330  & 2001 \tabularnewline
 & SA & 2517  & 2028  & 2.7 & 3195 2 & 3403  & 1.4 \tabularnewline
 & DEoptim & 2485  & 2052  & 282 & 487  & 2588  & 381 \tabularnewline
 & GA & 2400  & 2142  & 5.9 & 2364  & 2617  & 7.9 \tabularnewline
 & PSO & 1043  & 3509  & 43 & 1123  & 1340  & 59 \tabularnewline
\hline 
\end{tabular}
\par\end{centering}
\begin{centering}
Based on 10 replications. 

SMCO/-R/-BR and Group (I) algos use $\sqrt{d}=10$ uniform starting points (in parallel).

$\delta^*$ refers to bound buffer parameter.

STOGO excluded due to overly long computation time.
\par\end{centering}
\end{table}

\section{Additional Random Test Functions}\label{append-F}

\subsection{Conditional Moment Inequalities in Dynamic Binary Choice Models}

We now consider another random test function from economics, which summarizes identifying restrictions from a class of conditional moment inequalities in a dynamic (endogenous) binary choice model as studied in a recent paper by \cite{gao2025identification}. Specifically, consider the following binary choice model:
\[
Y_{it}=\ind\left\{ Z_{it}+\g_{0}X_{it}+\a_{i}+\e_{it}\geq0\right\} ,\quad t=1,...,T=3
\]
where $Z_{it}$ and $X_{it}$ are both taken to be scalar valued observed
covariates, with $Z_{it}$ being the exogenous covariate while $X_{it}$
being the endogenous covariates. We normalize the coefficient on $Z_{it}$
to $1$ and focus on the identification of the coefficient $\g_{0}$.
Here, $\a_{i}$ is an unobserved fixed effect, while $\e_{it}$ is
an unobserved time-varying error term whose distribution is left nonparametric
but assumed to satisfy a partial stationarity assumption across time
given the exogenous covariates $Z_{i}$. \cite{gao2025identification}  shows that under appropriate conditions, the identified
set for $\g_{0}$ is given by
\[
\G_{I}:=\left\{ \g:\ Q\left(\g\right)\leq0\right\} .
\]
where the criterion function is defined as a parametrized maximum
function 
$Q\left(\g\right):=\max_{c\in{\cal C},z\in{\cal Z}}Q\left(\g,c,z\right)$
with
\begin{align*}
Q\left(\g,c,z\right) := \max_{t}\P\left(\rest{Y_{it}=1,\ z_{t}+\g_{0}X_{it}\leq c}z\right)
+\max_{t}\P\left(\rest{Y_{it}=0,\ z_{t}+\g_{0}X_{it}\geq c}z\right)-1
\end{align*}
\cite{gao2025identification} numerically computes $Q\left(\g\right)$ under
a specific data generating process, and we replicate this exercise
here. Specifically, conditional on \emph{$Z_{i}=z\in{\cal Z}:=\left[-10,10\right]^{T}$},
we set the error term $\e_{it}\sim_{i.i.d.}N\left(0,1\right)$,
and the fixed effect 
$
\a_{i}=\rho_{\a}\cd\frac{1}{T}\sum_{t=1}^{T}z_{t,1}+\sqrt{1-\rho_{\a}^{2}}\xi_{i}
$
with $\rho_{\a}=0.1$ and $\xi_{i}\sim_{i.i.d.}\cN\left(0,1\right)$.
We consider the continuous specification where
$
X_{it}:=5\left(2\left(Y_{i,t-1}-0.5\right)+\eta_{it}\right)
$
with $\eta_{it}\sim_{i.i.d.}U\left(-1,1\right)$, and $\g_{0}=1$,
so that $X_{it}^{'}\g_{0}$ is again of the similar order of magnitude
as $z$. The initial condition $Y_{i,0}\sim_{i.i.d}Bernoulli\left(0.5\right).$ We then numerically compute $\hat{Q}\left(\g,c,z\right)$ as numerical
approximation of $Q\left(\g,c,z\right)$ via simulations as follows
\begin{align*}
\hat{Q}\left(\g,c,z\right):= \max_{t=1,...,T}\frac{1}{B}\sum_{b=1}^{B}Y_{bt}\ind\left\{ z_{t}+\g X_{bt}\leq c\right\} 
 +\max_{t=1,...,T}\frac{1}{B}\sum_{b=1}^{B}\left(1-Y_{bt}\right)\ind\left\{ z_{t}+\g X_{bt}\geq c\right\} -1
\end{align*}
based on $B=500$ simulations of $\left(Y_{bt},X_{bt}\right)$ based
on the DGP described above, conditional on each value of $z$. We
then optimize $\hat{Q}\left(\g,c,z\right)$ over $\left(c,z\right)\in{\cal C}\times{\cal Z}$,
with ${\cal C}:=\left[-30,30\right]$, to obtain $\hat{Q}\left(\g\right)$
as a numerical approximation of $Q\left(\g\right)$, using our SMCO
algorithms along with other competing optimization algorithms. We
then plot $\hat{Q}\left(\g\right)$ as a function of $\g$ in the
figure below, on a grid of 40 points for $\g\in\left[-20,20\right]$. \medskip

Note that the optimization of $\hat{Q}\left(\g,c,z\right)$ w.r.t.
$\left(c,z\right)$ is challenging for several reasons. To start,
the presence of indicator function $\ind\left\{ z_{t}+\g X_{bt}\leq c\right\} $
and $\ind\left\{ z_{t}+\g X_{bt}\geq c\right\} $ in $\hat{Q}\left(\g,c,z\right)$
makes the function nonsmooth and discontinuous in the parameters $\left(c,z\right)$
to be optimized over. Moreover, the indicator function also implies
that the gradient of $\hat{Q}$ is almost everywhere zero. The discreteness
in $\hat{Q}$ is further compounded by the fact that $Y_{bt}$ is
also a binary variable. Finally, there is little a priori information
on the structure of $\hat{Q}$ w.r.t $\left(c,z\right)$: the function
$\hat{Q}$ needs not be convex or monotone.\medskip

\begin{figure}
\caption{Dynamic Binary Choice: All}
\label{fig:Q_DBC}
\centering{}\includegraphics[scale=0.75]{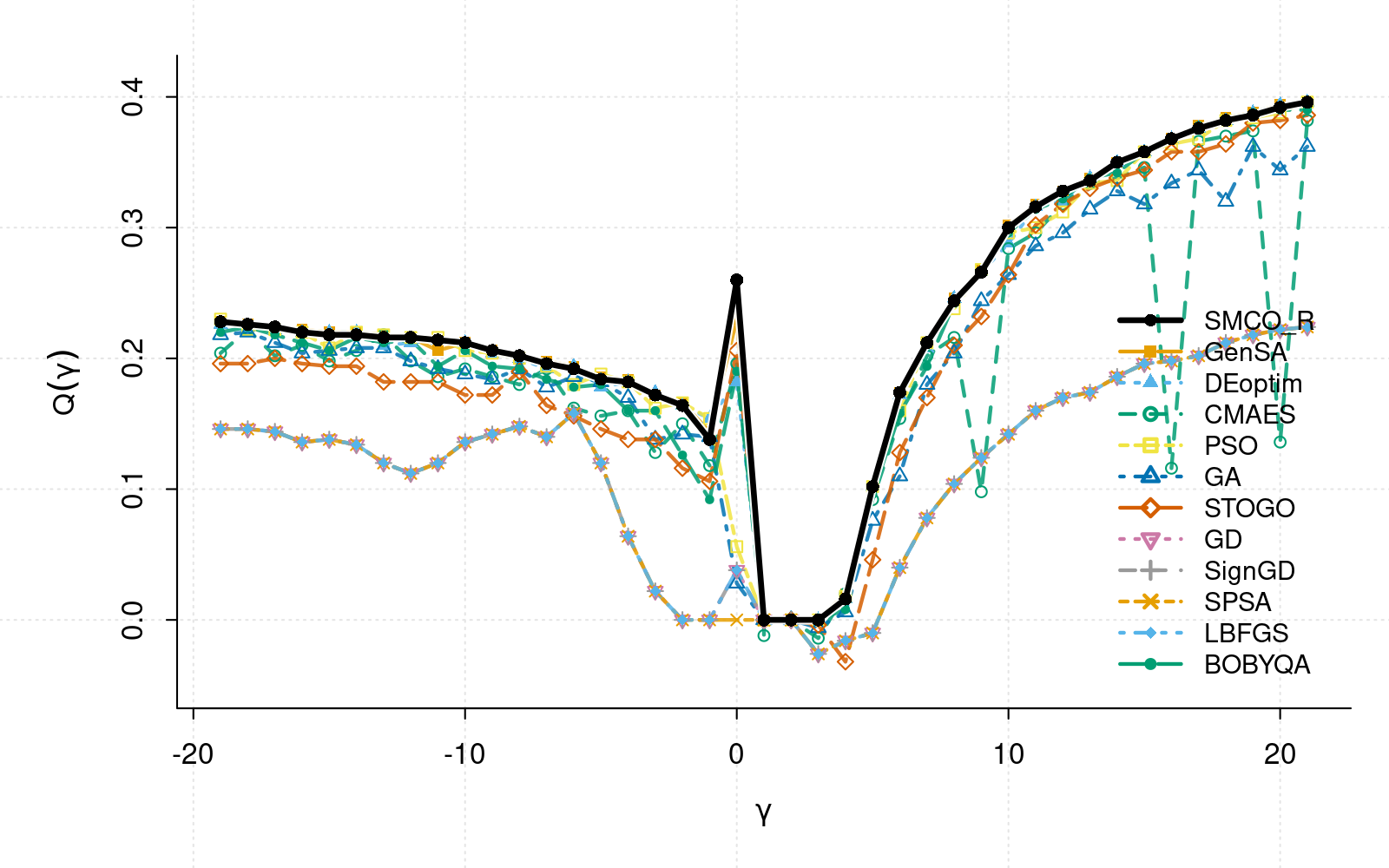}
\end{figure}

\begin{figure}
\caption{Dynamic Binary Choice: SMCO/-R/-BR, GenSA and DEoptim}
\label{fig:Q_DBC_3}
\centering{}\includegraphics[scale=0.75]{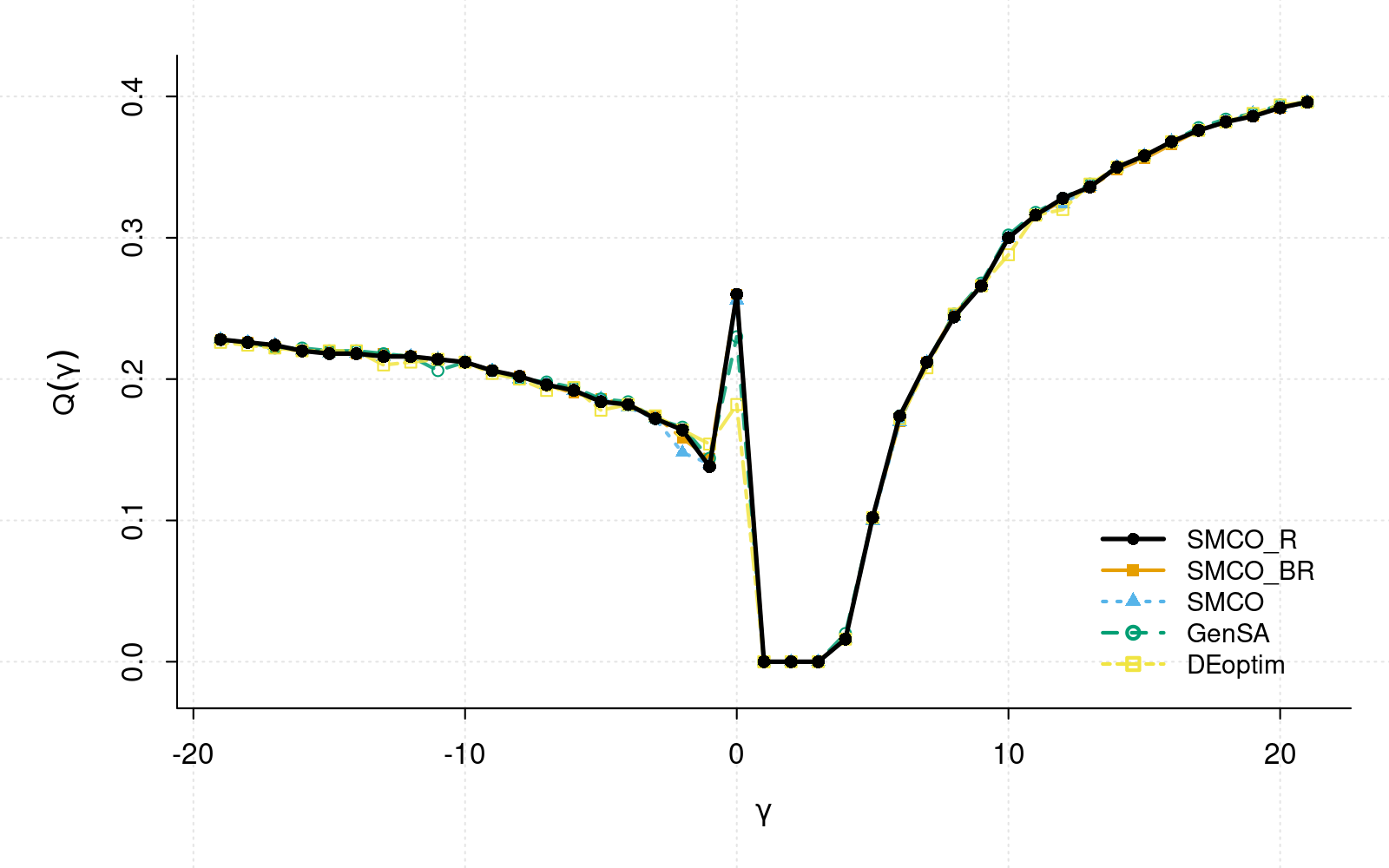}
\end{figure}

Figure \ref{fig:Q_DBC} contains plots of $\hat{Q}\left(\g\right)$
computed from SMCO\_BR along with other competing optimization algorithms,
while Figure \ref{fig:Q_DBC_3} contains results from the same exercise
but focuses on the comparison of SMCO\_BR with GenSA and DEoptim only.
Observe that SMCO\_BR performs robustly at the top across the whole
range of $\g$ under consideration, and performs on par with GenSA
and DEoptim while beating all other algorithms. It is also clear from
Figure \ref{fig:Q_DBC_3} that SMCO\_BR, GenSA and DEoptim yield very
similar curves, though our SMCO\_BR finds a noticeably larger maximum
than GenSA and DEoptim at $\g=0$.

\subsection{Square Loss Function for ReLU Neural Networks}

We also test the optimization algorithms on randomly generated functions
constructed from ReLU neural networks, which are well known to have
complicated and challenging optimization landscapes. Specifically,
we start with the following single-hidden-layer ReLU neural networks
with $d_{in}$ input nodes and $d_{hi}$ hidden nodes:
\begin{align*}
g\left(x,z\right)=\sum_{k=1}^{d_{hi}}w_{2k}\left[\sum_{j=1}^{d_{in}}w_{1kj}z_{j}+b_{1k}\right]_{+}+b_{2},
\end{align*}
where $z\in\mathbb{R}^{d_{in}}$ is a vector of input features and
$$x:=\left(\left(w_{1kj}\right)_{j=1,k=1}^{j=d_{in},k=d_{hi}},\left(w_{2k}\right)_{k=1}^{d_{hi}},\left(b_{1k}\right)_{k=1}^{d_{hi}},b_{2}\right)\in\mathbb{R}^{d_{x}}$$
collects all the neural-network parameters (weights and biases). Note
that the dimension of $x$ is $d_{x}:=d_{in}d_{hi}+2d_{hi}+1$. We construct two test functions based on the neural network described above.

We construct the test function based on the neural network using the square loss function, which widely used for regression purposes in statistics and machine learning. Essentially, we test how well the optimization algorithms can solve ``in-class'' regression problems under the neural network structure described above.

Specifically, we simulate $M$ copies of network parameters $x$ by sampling $w_{1kj}\sim_{i.i.d.}w_{2k}\sim_{i.i.d.}b_{2}\sim_{i.i.d.}U\left[-4,4\right]$
and $b_{1k}\sim_{i.i.d.}U\left[0,8\right]$. The positive support of
$b_{1k}$ is designed to make hidden nodes more likely to be activated,
creating relatively more challenging optimization landscapes. For
each simulated network parameter $x^{\left(m\right)}$, we construct an optimization test function as
\[
f_{m}\left(x\right):=\frac{1}{H}\sum_{h=1}^{H}\left(g\left(x,z^{\left(h\right)}\right)-g\left(x^{\left(m\right)},z^{\left(h\right)}\right)\right)^{2}
\]
where $\left\{ z^{\left(h\right)}:h=1,...,H\right\} $ are $H=1000$
input points, each of which is randomly sampled according to $z_{k}^{\left(h\right)}\sim_{i.i.d.}U\left[-4,4\right]$.
By construction, $f_{m}$ has a global minimum of value $0$ at $x^{\left(m\right)}$,
but there might be multiple global minima (and almost certainly many
local optima).

We then run the optimization algorithms to minimize $f_{m}$ on $\left[-10,10\right]^{d_{x}}$,
which contains the support of the simulated $x^{\left(m\right)}$
as a strict subset. The minimization of $f_{m}$ can be interpreted
as computing an approximation for the true neural-network regression
function $g\left(x^{\left(m\right)},\cdot\right)$ through the least-square
method evaluated at the $H=1000$ input points. We choose this ``noise-free''
setup to make sure that the true global minimum is known by design,
so that the comparison of the optimization algorithms below can be
cleanly interpreted as a comparison of optimization performances on
complex regression functions induced by the neural networks.\footnote{An alternative exercise would be to simulate a noisy outcome $y^{\left(m,h\right)}:=g\left(x^{\left(m\right)},z^{\left(h\right)}\right)+\epsilon^{\left(h\right)}$
and minimize $\tilde{f}_{m}\left(x\right)=\frac{1}{H}\sum_{h=1}^{H}\left(g\left(x,z^{\left(h\right)}\right)-y^{\left(m,h\right)}\right)^{2},$with
$\lim_{H\to\infty}\tilde{f}_{m}\left(x\right)$ being the (asymptotic)
out-of-sample error, whose global minimum becomes $\text{Var}\left(\epsilon^{\left(h\right)}\right)$.
We could then assess the out-of-sample MSE using an additional sample
of $\left(z,y\right)$, obtained by the optimization algorithms, which
resembles a real-world use of neural networks more. We note, however,
such an exercise can be misleading as a measure of optimization performance,
since it contaminates the optimization problem with statistical complications.
For example, it could be the case that a better optimization algorithm
is able to ``overfit'' the model, generates smaller in-sample errors
but larger out-of-sample errors in the test sample. However, this
should not be interpreted as evidence against the said optimization
algorithm, whose sole purpose is to minimize functions. } 

Table \ref{tab:Comp_ANN} contains our comparison using the following configurations: (i) We set the input dimension of the neural networks to be $d_{in} = 3$ and the number of nodes in the hidden layer to be $d_{hi} = 5$, with a total of $d=26$ parameters. (ii) We set $d_{in}=3$, $d_{hi} = 10$ with a total of $d=51$ parameters. (iii) We set $d_{in}=6$, $d_{hi} = 15$ with a total of $d=121$ parameters. We run our SMCO/-R/-BR algorithms with the default hyperparameters as well as another pass with the maximum iterations set to 1000 (the default is 500). We report results in minimizing the square loss test function with parameters randomly generated over $M = 10$ replications, for each of which we run the optimization algorithms once. We then compute the average performance measures across the $M = 10$ replications. It is clear from Table \ref{tab:Comp_ANN} that our SMCO algorithms still perform reasonably well, with RMSE competitive with or better than most of the included algorithms (with GenSA being an notable exception). Again, it is worth noting the remarkable stability of our algorithms in terms of AE99 (say, in comparison with PSO). 

\begin{sidewaystable}
\begin{centering}
\caption{\label{tab:Comp_ANN} Square Loss Minimization for ReLU Neural Networks}
\begin{tabular}{|ccccc|ccc|ccc|}
\hline 
\multicolumn{2}{|c}{Dimension} & \multicolumn{3}{c|}{\textbf{$d=26$}} & \multicolumn{3}{c|}{\textbf{$d=51$}} & \multicolumn{3}{c|}{\textbf{$d=121$}}\tabularnewline
\hline 
\multicolumn{2}{|c}{Layer Dimensions} & \multicolumn{3}{c|}{$d_{in}=3,\,d_{hi}=5$} & \multicolumn{3}{c|}{$d_{in}=3,\,d_{hi}=10$} & \multicolumn{3}{c|}{$d_{in}=6,\,d_{hi}=15$}\tabularnewline
\hline 
\multicolumn{2}{|c}{Performance} & \multicolumn{1}{c|}{RMSE} & \multicolumn{1}{c|}{AE99} & Time & \multicolumn{1}{c|}{RMSE} & \multicolumn{1}{c|}{AE99} & Time & \multicolumn{1}{c|}{RMSE} & \multicolumn{1}{c|}{AE99} & Time\tabularnewline
\hline 
\multicolumn{2}{|c}{\textbf{Best}} & \multicolumn{3}{c|}{\textbf{1.13e-8}} & \multicolumn{3}{c|}{\textbf{0.01778 }} & \multicolumn{3}{c|}{\textbf{0.008116 }}\tabularnewline
\hline 
\multirow{6}{*}{(I)} & GD & 547.8  & 1273  & 10.0  & 420.7  & 570.5  & 45.7  & 1225  & 2141  & 154.9 \tabularnewline
 & SignGD & 14.4  & 25.6  & 10.1  & 16.7  & 21.2  & 46.8  & 130.2  & 156.4  & 160.4 \tabularnewline
 & ADAM & 3.7  & 8.4  & 7.6  & 6.7  & 10.1  & 18.3  & 43.6  & 53.2  & 32.3 \tabularnewline
 & SPSA & 102576  & 146316 & 0.8  & 139017 & 204410  & 1.6  & 809811  & 1125694  & 2.4 \tabularnewline
 & L-BFGS & 5.3  & 11.4  & 1.3  & 6.3  & 8.9  & 5.4  & 68.7  & 85.5  & 17.8 \tabularnewline
 & BOBYQA & 22.9  & 40.3  & 0.4  & 43.7  & 56.1  & 0.9  & 907.1  & 1131  & 2.3 \tabularnewline
\hline 
\multirow{3}{*}{maxit = 500} & \textbf{SMCO} & \textbf{2.7 } & \textbf{4.0 } & \textbf{5.2 } & \textbf{6.0 } & \textbf{6.6 } & \textbf{23.2 } & \textbf{50.7 } & \textbf{54.9 } & \textbf{79.9 }\tabularnewline
 & \textbf{SMCO-R} & \textbf{1.5 } & \textbf{2.5 } & \textbf{5.2 } & \textbf{4.4 } & \textbf{4.8 } & \textbf{23.3 } & \textbf{37.3 } & \textbf{42.5 } & \textbf{80.4 }\tabularnewline
 & \textbf{SMCO-BR} & \textbf{2.7 } & \textbf{3.3 } & \textbf{5.3 } & \textbf{6.1 } & \textbf{7.0 } & \textbf{23.6 } & \textbf{55.1 } & \textbf{60.2 } & \textbf{80.4 }\tabularnewline
\hline 
\multirow{3}{*}{maxit = 1000} & \textbf{SMCO} & \textbf{0.9 } & \textbf{1.9 } & \textbf{26.9 } & \textbf{4.0 } & \textbf{4.8 } & \textbf{60.7 } & \textbf{35.8 } & \textbf{38.4 } & \textbf{153.7 }\tabularnewline
 & \textbf{SMCO-R} & \textbf{0.6 } & \textbf{1.3 } & \textbf{26.9 } & \textbf{3.5 } & \textbf{4.2 } & \textbf{60.8 } & \textbf{31.1 } & \textbf{34.3 } & \textbf{153.9 }\tabularnewline
 & \textbf{SMCO-BR} & \textbf{0.8 } & \textbf{1.3 } & \textbf{27.2 } & \textbf{4.4 } & \textbf{5.0 } & \textbf{61.2 } & \textbf{34.4 } & \textbf{39.6 } & \textbf{153.4 }\tabularnewline
\hline 
\multirow{6}{*}{(II)} & GenSA & 0.2  & 0.3  & 167.8  & 0.4  & 1.1  & 749.0  & 3.1  & 6.4  & 2769 \tabularnewline
 & SA & 2530  & 6747  & 0.2  & 2639 & 3936  & 0.3  & 36899  & 71465  & 0.1 \tabularnewline
 & DEoptim & 68.8  & 93.8  & 8.9  & 293.4  & 383.3  & 27.9  & 2755  & 3194  & 96.3 \tabularnewline
 & STOGO & 2.9  & 2.9  & 94.9  & 2181  & 2181  & 246.0  & 2315  & 2315  & 813.1\tabularnewline
 & GA & 132.7  & 191.6  & 0.9  & 145.3  & 197.4  & 1.3  & 823.6  & 1182  & 1.8 \tabularnewline
 & PSO & 6.9  & 13.2  & 4.1  & 13.3 & 19.3  & 7.4  & 106.3  & 130.8  & 14.1 \tabularnewline
\hline 
\end{tabular}
\par\end{centering}
Based on 10 replications. SMCO/-R/-BR and Group (I) algos use $\sqrt{d}\approx5,7,11$ uniform starting points (in parallel). 
\end{sidewaystable}

\section{Additional Random Statistical Test Functions }\label{append-G}

All numerical comparisons reported in this Appendix are for statistical estimation and prediction problems, and are implemented in \texttt{R} 4.4.1 on a Windows server with an NVIDIA GeForce RTX 4080 GPU.  
We use the \texttt{sgd} for SGD, \texttt{GenSA} solver for GenSA, and the \texttt{psoptim} solver in package \texttt{pso} for PSO available in \texttt{R}. 

We have implemented the following versions of the SMCO algorithms in this Appendix: SMCO (i.e., step-adaptive SMCO  with a single starting point as described in Appendix \ref{append-B}); SMCOfs (i.e., the original SMCO Algorithm \ref{alg:smco_algo} with a fixed single starting point); SMCOfs\_XXX (i.e., the existing XXX algorithm in \texttt{R} package 
using the SMCOfs's 20th iterative point $\hat{x}_{20}=\frac{S_{20}^{\bmtheta^*}}{20}$ as its starting value); and SMCO-BR (i.e., the “boosted” version of SMCO as described in Appendix  \ref{append-B}).
%

Examples \ref{example-1} and \ref{example-2} below are about correctly specified statistical models with Cauchy density, which is a leading example of a density with fat-tails. Cauchy likelihood is known to be difficult to optimize \cite{Brooks1995}. Examples \ref{example-3} and \ref{example-4} are about using neural networks to approximate true unknown functions for regression and classification problems, in which the random objective functions are not globally concave/convex with respect to the neural network parameters.

\begin{example}[Cauchy log-likelihood function with a small sample size $N$]\label{example-1}
    Let $X$ be a Cauchy random variable.
The log-likelihood surface of a Cauchy density function $g(X)=c/\{\pi(c^2+(X-x)^2)\}$ is known to be difficult to optimize when estimating $x$ under some known value $c>0$. Let a random sample $\{X_i\}_{i=1}^N$ with $N=8$ be generated from the Cauchy density $g(X)$ with $c=0.1$, which is
$X_i\in\{$-4.20, -2.85, -2.30, -1.02, 0.70, 0.98, 2.72, 3.50$\}$. The joint log-likelihood is
$N \log c-\sum_{i=1}^N \log \left\{c^2+\left(X_i-x\right)^2\right\}-N \log \pi $ with $N=8$. The maximum likelihood estimate for $x$ is the global maximizer $x^*=0.73$ of
$$f(x)=-\sum_{i=1}^8 \log \left\{c^2+\left(X_i-x\right)^2\right\},~~~\text{with}~f(x^*)=-5.36.$$
\end{example}

\begin{figure}
  \centering  
  \captionsetup[subfloat]{labelfont=scriptsize,textfont=tiny}
  
  \subfloat[SGD]  
  {
      \label{One-dimensional example(sgd)}\includegraphics[width=0.18\textwidth]{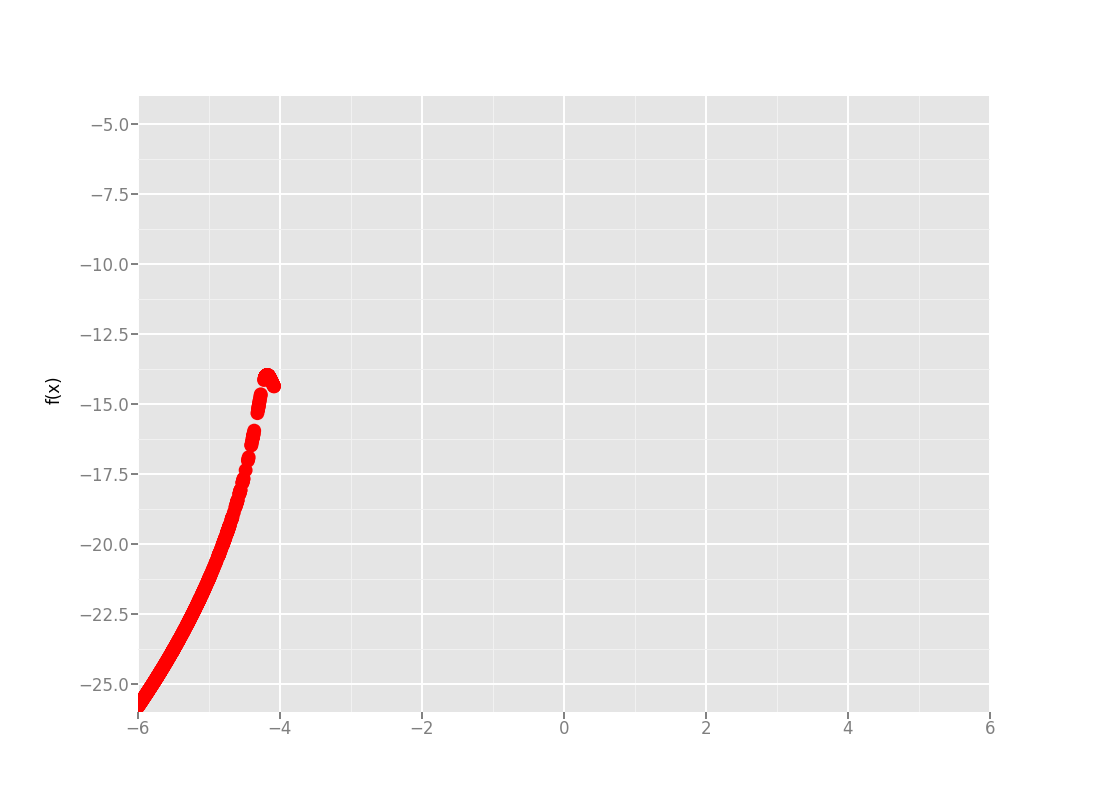}
  }
  \hfill
  \subfloat[GenSA]  
  {
      \label{One-dimensional example(simulated annealing)}\includegraphics[width=0.18\textwidth]{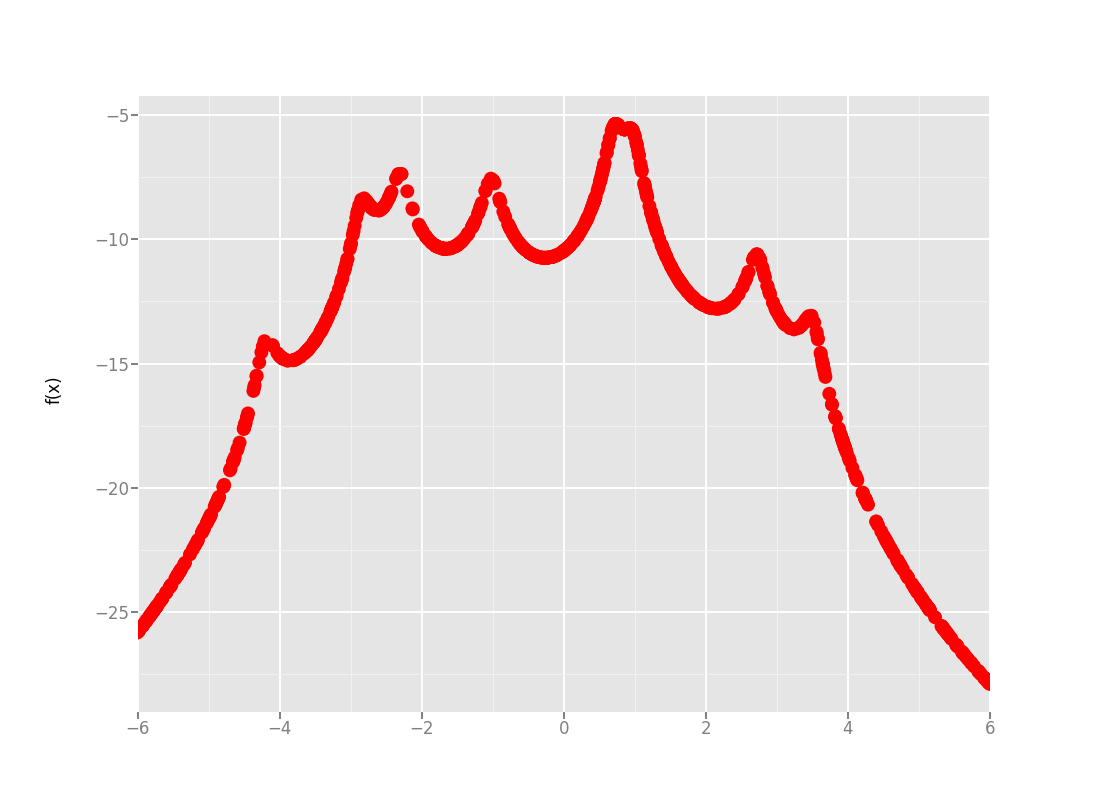}
  }
  \hfill
  \subfloat[PSO]  
  {
      \label{One-dimensional example(Particle Swarm Optimization)}\includegraphics[width=0.18\textwidth]{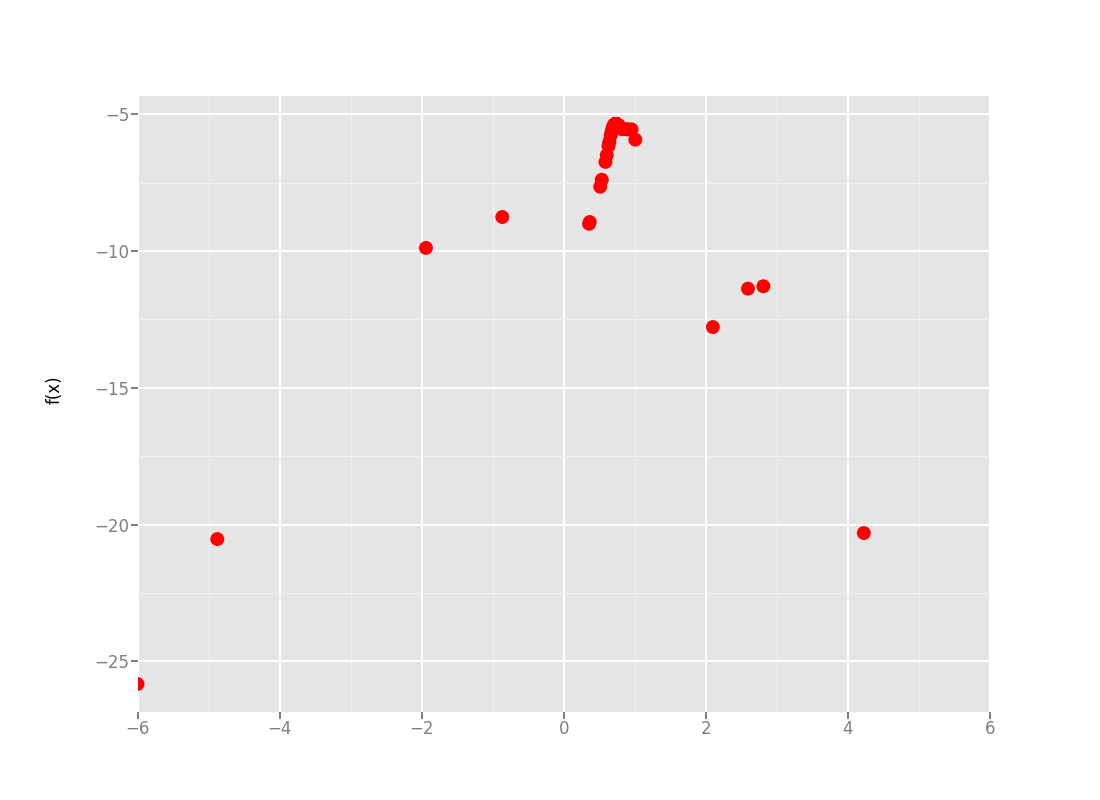}
  }
  \hfill
  \subfloat[SMCO]  
  {
      \label{One-dimensional example(Our approach)}\includegraphics[width=0.18\textwidth]{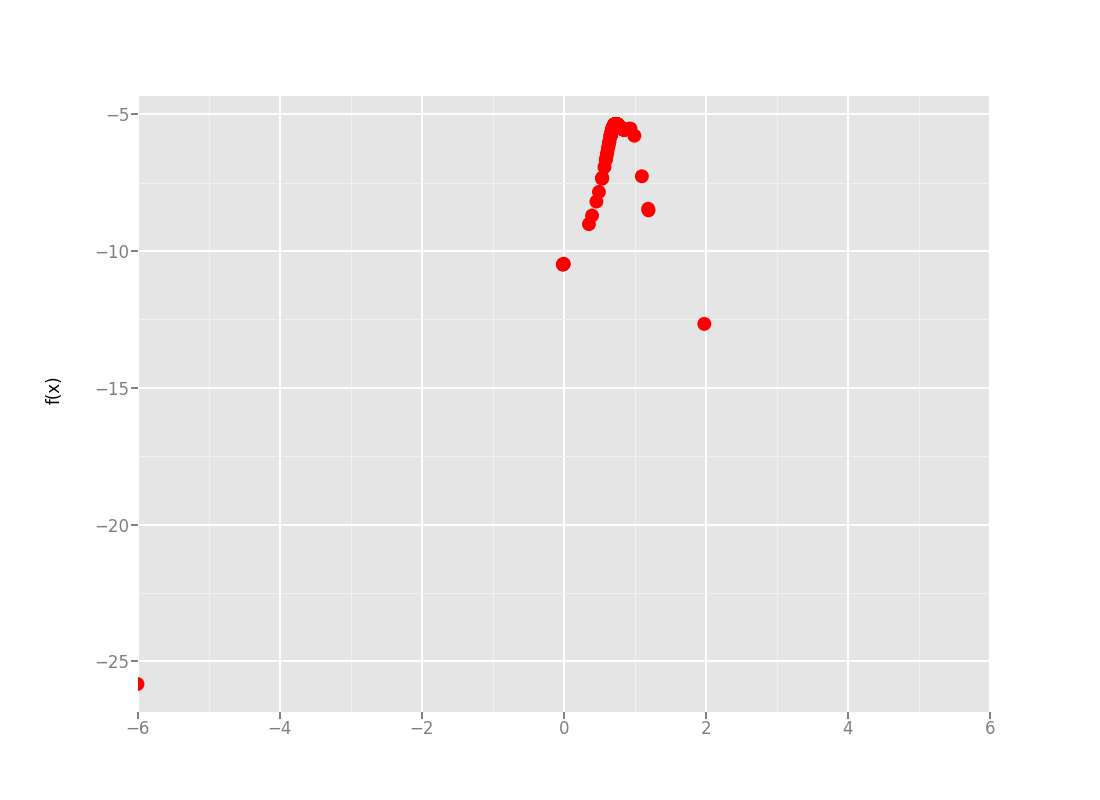}
  }
  \hfill
  \subfloat[SMCOfs]  
  {
      \label{One-dimensional example(Our approach)2}\includegraphics[width=0.18\textwidth]{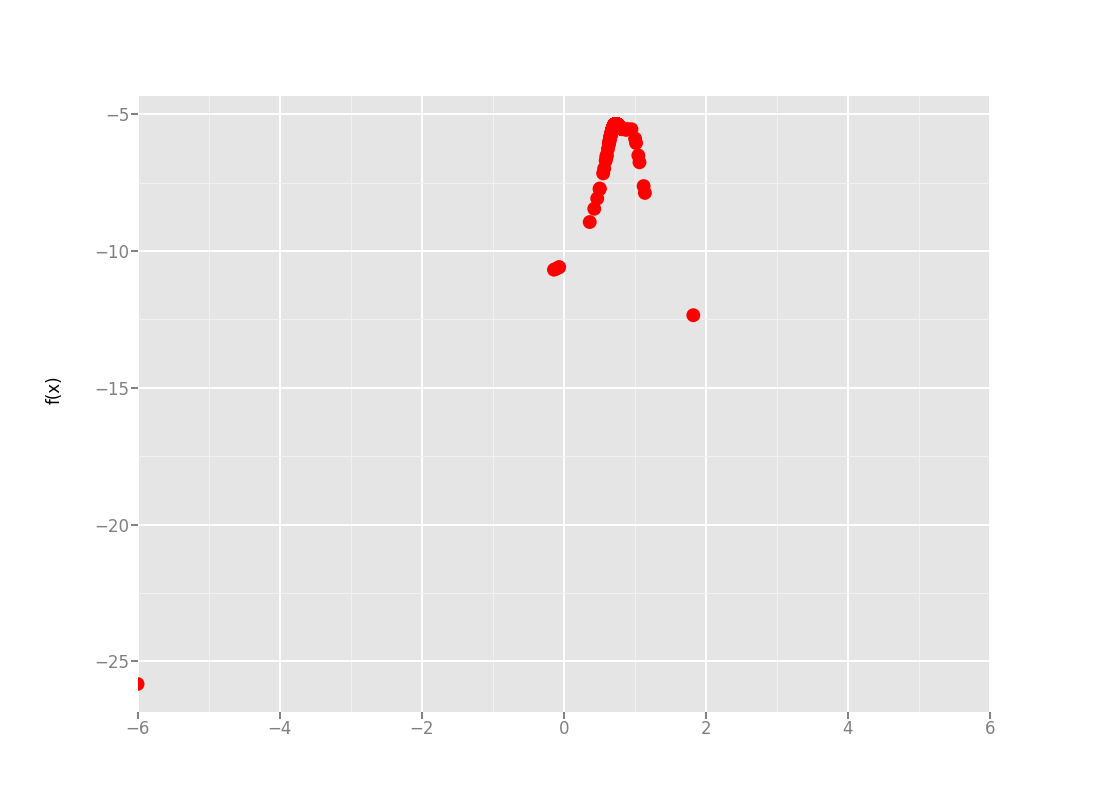}
  }
   
  \caption{The $1000$ iterative point-specific function values of Cauchy log-likelihood in Example \ref{example-1} under SGD, GenSA, PSO, SMCO and SMCOfs. SGD attains a local solution $\tilde{x}=-4.2$ with $f(\tilde{x})=-14.02$, while GenSA, PSO, SMCO and SMCOfs attain the global maximum $x^*=0.73$ with $f(x^*)=-5.36$.}
  \label{cauchycase}
\end{figure}

For Example \ref{example-1}, Figure \ref{cauchycase} depicts one thousand iterative points of Cauchy log-likelihood $f(x)$ under the five optimization algorithms: SGD, GenSA, PSO, SMCO and SMCOfs.  All algorithms have the same initial starting point value -6. The SGD has the learning rate 0.01 and batch size 1 observation. 
The SMCOfs has uniform distributions Unif (-6.5, -5.5) and Unif (5.5, 6.5) as its left and right arm distributions respectively. 
The depicted scatterplot of the Cauchy log-likelihood reveals the following: (1) the SGD attains the local maximum $(-4.2, -14.02)$; while GenSA, PSO, SMCO and SMCOfs  attain the global maximum $(0.73,-5.36)$. (2) Our SMCO delivers faster exploration to reach the global maximum  in terms of more accumulated points around the global maximum  under the accuracy level $\varepsilon = 10^{-7}$, while our SMCOfs performs very similar to SMCO for this one dimensional example. (3) In terms of running time, SGD takes 41 milliseconds but attains a local maximum point only; GenSA takes 205.2 milliseconds, PSO needs 3250.8 milliseconds, whereas our SMCO and  SMCOfs takes 110.5 and 128.2 milliseconds respectively.

\begin{example}[Regression model with Cauchy error]\label{example-2}
  Recently Feng et al. (2024) propose a new robust antitonic score matching (asm) method for linear regression models with heavy tailed error distributions such as Cauchy error.
We follow their numerical experiment design by generating data from a linear regression model with a Cauchy error:
$$Y_i= Z_i^\top x+\epsilon_i, i=1,\cdots, N,~~~\text{with}~N=600, $$
where $Z_i = (1, Z_{i1},\cdots, Z_{i5})^\top$ and $Z_{ij} $
is independently generated from normal distribution with $E(Z_{ij})={\rm var}(Z_{ij})=1$, and the true unknown values $x^*=(x_{1}^*,...,x_{6}^*)'$ with $x_{1}^*=1$ and $\|x_{-1}^*\|_2=3$. The error term $\epsilon_i$ is generated from the
Cauchy density function $g(\epsilon)=\frac{c}{\pi (c^2+\epsilon^2)}$ under two cases $c=1,~2$.  
\end{example}

For Example \ref{example-2}, we estimate/calculate $\hat{x}$ using three optimization algorithms: our SMCO-BR,  the R function \texttt{nlm}, and the antitonic score matching (asm) method of Feng et al. \cite{Feng2024} in the R package \texttt{asm}. 
We implement the \texttt{nlm} and \texttt{asm} packages using zero vector as the initial values  and set the initial value of \texttt{asm} minus and plus 2 to be the lower and upper bounds respectively for SMCO-BR in each replicate.
\texttt{nlm} and SMCO-BR use Cauchy density-based likelihood as the objective functions. Table \ref{example5} shows the slightly better estimates performance of our SMCO-BR in terms of both smaller median  and standard deviation of parameter root mean squared error (RMSE) $\sqrt{\sum_{j=1}^6(\hat{x}_j-x_j^*)^2/6}$ after 200 replicates.

\begin{table}
\renewcommand{\thetable}{S18}
\caption{Median (over 200 replicates) of the parameter root mean squared error (RMSE) $\|\hat{x}-x^*\|_2/\sqrt{6}$ and standard deviation (SD) of the parameter RMSE, with $\hat x$ computed using SMCO\text{-}R, SMCO\text{-}BR, \texttt{nlm}, and \texttt{asm} in Example~2.}\label{example5}
\begin{center}   
\begin{tabular}{|cccccccccccccc|}

\cline{1-14}
\multicolumn{2}{|c|}{Method} && \multicolumn{2}{c}{ Cauchy $c$} && \multicolumn{2}{|c|}{RMSE} && \multicolumn{2}{c|}{SD} && \multicolumn{2}{c|}{Time}\\
\cline{1-14}
\multicolumn{2}{|c|}{SMCO-R}  && \multicolumn{2}{c}{1} && \multicolumn{2}{|c|}{0.0700} && \multicolumn{2}{c|}{0.0008} && \multicolumn{2}{c|}{1.6446}\\
\multicolumn{2}{|c|}{}        && \multicolumn{2}{c}{2} && \multicolumn{2}{|c|}{0.0688} && \multicolumn{2}{c|}{0.0011} && \multicolumn{2}{c|}{2.0276}\\
\cline{1-14}
\multicolumn{2}{|c|}{SMCO-BR} && \multicolumn{2}{c}{1} && \multicolumn{2}{|c|}{0.0703} && \multicolumn{2}{c|}{0.0008} && \multicolumn{2}{c|}{2.6876}\\
\multicolumn{2}{|c|}{}        && \multicolumn{2}{c}{2} && \multicolumn{2}{|c|}{0.0695} && \multicolumn{2}{c|}{0.0011} && \multicolumn{2}{c|}{2.1741}\\
\cline{1-14}
\multicolumn{2}{|c|}{\texttt{nlm}} && \multicolumn{2}{c}{1} && \multicolumn{2}{|c|}{0.0731} && \multicolumn{2}{c|}{0.0009} && \multicolumn{2}{c|}{0.0022}\\
\multicolumn{2}{|c|}{}            && \multicolumn{2}{c}{2} && \multicolumn{2}{|c|}{0.0719} && \multicolumn{2}{c|}{0.0012} && \multicolumn{2}{c|}{0.0037}\\
\cline{1-14}
\multicolumn{2}{|c|}{\texttt{asm}} && \multicolumn{2}{c}{1} && \multicolumn{2}{|c|}{0.0763} && \multicolumn{2}{c|}{0.0010} && \multicolumn{2}{c|}{0.1906}\\
\multicolumn{2}{|c|}{}            && \multicolumn{2}{c}{2} && \multicolumn{2}{|c|}{0.0738} && \multicolumn{2}{c|}{0.0013} && \multicolumn{2}{c|}{0.2789}\\

\cline{1-14}
\end{tabular}
\end{center}
\end{table}


\begin{example}[Regression using neural network approximation]\label{example-3}
  We observe five-dimensional covariates $Z=(Z_{1},\cdots, Z_{5})^{\top}$ following the five-dimensional normal distribution with $E(Z)=(0)_5^{\top}$ and the covariance element $0.01^{l-l'}$, but the response generation only uses three dimensions with
$Y=\sin(Z_1) + 0.5 \cos(Z_2) + 0.3 Z_3^2+\epsilon$, where $\epsilon$ follows a normal distribution with zero mean and standard deviation 0.01. We use a single-hidden layer ReLU neural network $g(Z, x)$ to approximate the true unknown conditional mean function $E[Y|Z]$:
\begin{align*}
g(Z_i, x)=\sum_{k=1}^3\omega_{2k}{\rm max} (0, \sum_{j=1}^5\omega_{1kj}Z_{ij}+b_{1k})+b_{2},
\end{align*}
where its unknown neural network parameter vector $x$ consists of 22 weights and biases $\{\omega_1,b_1, \omega_2,b_2 \}$ with $\omega_{1, 5\times 3},b_{1,3\times 1}, \omega_{2,3\times 1},b_{2,1\times 1} $ for $k=1, 2, 3$.  
We obtain an estimator $\hat{x}$ by nonlinear least squares using a training data set of size $N$, that is, $\hat{x}=\arg\max_{x} f(x)$ with $f(x)=-\frac{1}{N}\sum_{i=1}^N\left(Y_i - g( Z_i, x)\right)^2$.  
\end{example} 

For Example \ref{example-3}, we compute $\hat{x}$ using three optimization algorithms: our SMCO-BR, the L-BFGS \cite{LiuBFGS1989}, the SMCOfs\_L-BFGS and R package \texttt{nnet}. L-BFGS is an improved variant of SGD with batch updates based on the backpropagation algorithm, which is used to calculate the gradient of the loss function with respect to the weights and biases of the network. We implement L-BFGS using uniform distribution Unif(0,1) to generate the initial values of the elements of $\omega_1$ and  $\omega_2$ and zero as the the initial values of $b_1$ and  $b_2$ in each replicate, which is denoted as L-BFGS.  We set L-BFGS's initial value minus and plus 1 to be the lower and upper bounds for SMCO-BR and SMCOfs.
SMCO-BR considers the number of starting points
is set to be 50 and the maximum number of iterations to be 300 with boosting iteration number 100.

We consider three training sizes with $N=1000, 3000$ and $N_{\rm test}=500$ testing dataset.
Table \ref{example3} shows the averaged (over 500 replicates) mean square errors (MSE) using the test sample of size $N_{\rm test}=500$ with ${\rm MSE}=\frac{1}{N_{\rm test}}\sum_{i\in \mathcal{D}_{\rm test}}\{Y_i-g(Z_i;\hat{x})\}^2$
across different accuracy levels $\varepsilon$ to stop the algorithms.
Table \ref{example3} shows the better prediction performance of our SMCO-BR over L-BFGS in terms of a smaller MSE and standard deviation using training and testing datasets for different accuracy levels. 
Interestingly, the L-BFGS algorithm under the initial SMCO-based value after 20 samplings, denoted SMCOfs\_L-BFGS, generates better prediction performance in the test data than that with L-BFGS. Table \ref{example3} again shows the advantage of combining SMCOfs with existing advanced optimization algorithms to improve performance. 
The R package \texttt{ nnet} is also considered for fitting the single hidden layer ReLU neural network but generates larger averaged MSE values for training and testing data sets, respectively.

\begin{table}     
\caption{Averaged (over 500 replicates) MSE (Training | Testing), standard deviation (Training | Testing) and running time (second), with $\hat x$ computed using SMCO-BR (without parallel), SMCOfs\_L-BFGS, L-BFGS, and R package \texttt{nnet} in Example \ref{example-3}.}      
\label{table1}      
\begin{center}   
\fontsize{8}{11}\selectfont
\begin{tabular}{ccccccccccc} 
\hline
&\mbox{MSE$\times 10^2$} & \mbox{SD$\times 10^2$} & \mbox{Time} &&\mbox{MSE$\times 10^2$} & \mbox{SD$\times 10^2$} & \mbox{Time} \\ \cline{2-4} \cline{6-8}
& \multicolumn{3}{c}{$N=1000~|~N_{\rm test}=500$}&& \multicolumn{3}{c}{$N=3000~|~N_{\rm test}=500$} \\ \hline
\multicolumn{8}{c}{training accuracy level $\varepsilon=10^{-5}$}\\ 
\hline 
SMCO-BR & \bf 9.36 | 10.17 & \bf 1.47 | 1.86& 51.02&& \bf 9.62 | 9.81 & \bf 1.06 | 1.73&118.27 \\
SMCOfs\_L-BFGS & 14.78 | 15.61 & 4.86 | 5.44&0.07 && 15.15 | 15.43 & 4.85 | 5.36&0.14 \\
L-BFGS & 46.41 | 46.83 & 28.10 | 28.34&0.03 && 48.67 | 48.24 & 29.51 | 28.90& 0.06\\
\texttt{nnet} & 29.64 | 30.17 & 7.95 | 8.56&0.28 && 32.30 | 32.44 & 5.81 | 6.94&0.49 \\
\hline
\multicolumn{8}{c}{training accuracy level $\varepsilon=10^{-8}$}\\ 
\hline 
SMCO-BR& {\bf 9.08 | 9.62} &{\bf  1.20 | 1.76} &61.25   && {\bf 9.21| 9.35} &{\bf 0.93 | 1.84}& 133.84 \\
SMCOfs\_L-BFGS& {13.83 | 14.67} & { 4.90 | 5.22} & 0.20  && {14.38 | 14.66} & { 4.69 | 4.88} & 0.37   \\
L-BFGS& 47.75 | 48.34 &33.54 | 33.19 & 0.05&& 46.58 |46.59 & 28.99 |28.81 & 0.10\\
\texttt{nnet}& 30.15 | 30.69  & 9.20 | 9.84 & 0.28 && 32.72 | 32.97  &  8.25 | 9.16 & 0.49 \\
\hline
\end{tabular}\label{example3}     
\end{center}      
\end{table}

\begin{example}[Classification using neural network approximation]\label{example-4}
   We generate binary response $Y$ based on five-dimensional covariates $Z=\left(Z_1, \cdots, Z_5\right)^{\top}$. The true data generating process is: $Y=1$ when $\{1/[1+\exp(2.5-3Z_1 Z_2-1.5\left|Z_3\right|-\left|Z_4 Z_5\right|]\}>0.5$; and $Y=0$ otherwise, where $Z$ follows the five-dimensional normal distribution with $E(Z)=(1)_5^{\top}$ for half of the data sets and the other half of the data sets using $E(Z)=(0)_5^{\top}$ with common covariance matrix $0.5^{l-l^{\prime}}$. 

We generate a random sample of 1500 observations from the distribution of $(Y, Z)$, with $N=1000$ as the training sample and the remaining $N_{\text {test }}=500$ as the test sample. We consider a single-layer neural network $g(Z, x)$ with three hidden units and one output to approximate the true unknown $E[Y \mid Z]$. Then the parameter vector $x$ consists of 22 weights and biases $\left\{\omega_1, b_1, \omega_2, b_2\right\}$ with $\omega_{1,5 \times 3}, b_{1,3 \times 1}, \omega_{2,3 \times 1}, b_{2,1 \times 1}$  and $o_k=\max \left(0, \sum_{j=1}^5 \omega_{1 k j} Z_j+b_{1 k}\right)$ for $k=1,2,3$ and output is $g\left(Z_i ; x\right)=\operatorname{sigmoid}\left(\sum_{k=1}^3 \omega_{2 k} o_k+b_2\right)$. We obtain the estimator $\hat{x}$ by maximizing the following negative quadratic loss function:

$$
f(x)=-\frac{1}{N} \sum_{i=1}^N\left\{Y_i-g\left(Z_i, x\right)\right\}^2, g\left(Z_i, x\right)=\operatorname{sigmoid}\left(\sum_{k=1}^3 \omega_{2 k} \max \left(0, \sum_{j=1}^5 \omega_{1 k j} Z_{i j}+b_{1 k}\right)+b_2\right) .
$$ 
\end{example}

For Example \ref{example-4}, we calculate $\hat{x}$ using three optimization algorithms: our SMCO-BR, the L-BFGS in the R package \texttt{nnet}, and the SMCOfs\_\texttt{nnet}. 
We implement the L-BFGS in R package \texttt{nnet} using uniform distribution Unif(0,1) to generate the initial values of the elements of $\omega_1$ and $\omega_2$ and zero as the initial values of $b_1$ and $b_2$ and set the initial value of L-BFGS minus and plus 3 to be the lower and upper bounds respectively for SMCO-BR and SMCOfs in each replicate.
SMCO-BR also considers that the number of starting points
is set to 50 and the maximum number of iterations to be 300 with boosting iteration number 100.
The SMCOfs\_\texttt{nnet} is implemented the same way as the \texttt{nnet} except using our SMCOfs' 20th sampling point as its warm starting point. Table \ref{example4} shows the misclassification rate (MR) values on the test sample of size $N_{\text {test }}=500$ with $\mathrm{MR}=1-\frac{1}{N_{\text {test }}} \sum_{i \in \mathcal{D}_{\text {test }}} I\left\{Y_i=I\left[g\left(Z_i ; \hat{x}\right)>0.5\right]\right\}$ across different accuracy levels $\varepsilon$ to stop the algorithms and take an average of MR's after 200 replicates at each accuracy level. Table \ref{example4} shows the best classification prediction performance of our SMCO-BR. Since \texttt{nnet} adequately accounts for the scalability of memory efficiency with storage of a small number of vectors and also takes advantage of the Hessian approximation for SGD iterates, its running speed is faster than that of our SMCO-BR on average. Interestingly, the SMCOfs\_\texttt{nnet}, the \texttt{nnet} algorithm under SMCOfs-based initial value after 20 samplings, generates better classification performance on the test data than that of \texttt{nnet}. Table \ref{example4} again shows the advantage of combining SMCOfs (or even better SMCO) with existing
advanced optimization algorithms to improve performance.

\begin{table}
\caption{Averaged (over 200 replicates) out-of-sample mis-classification rate (MR), standard deviation (SD) and running time (seconds) with $\hat x$ computed using SMCO-BR (without parallel), \texttt{nnet} and SMCOfs\_\texttt{nnet} in Example \ref{example-4}.}\label{example4}
\begin{center}   
\begin{tabular}{cccccccccccccc}
\cline{1-14} 
\multicolumn{2}{c}{Method} && \multicolumn{2}{c}{Training Accuracy $\varepsilon$} && \multicolumn{2}{c}{MR} && \multicolumn{2}{c}{SD}  && \multicolumn{2}{c}{Time}\\

\cline{1-14} 
\multicolumn{2}{c}{} && \multicolumn{2}{c}{$ 10^{-3}$} && \multicolumn{2}{c}{\textbf{0.121}} && \multicolumn{2}{c}{\textbf{0.016}}&& \multicolumn{2}{c}{19.47}\\
\multicolumn{2}{c}{SMCO-BR} && \multicolumn{2}{c}{$ 10^{-5}$} && \multicolumn{2}{c}{\textbf{0.112}} && \multicolumn{2}{c}{\textbf{0.013}}&& \multicolumn{2}{c}{75.82}\\
\multicolumn{2}{c}{} && \multicolumn{2}{c}{$ 10^{-7}$} && \multicolumn{2}{c}{\textbf{0.112}} && \multicolumn{2}{c}{\textbf{0.012}}&& \multicolumn{2}{c}{146.33}\\

\cline{1-14} 

\multicolumn{2}{c}{} && \multicolumn{2}{c}{$ 10^{-3}$} && \multicolumn{2}{c}{0.129} && \multicolumn{2}{c}{0.029}&& \multicolumn{2}{c}{0.205}\\
\multicolumn{2}{c}{SMCOfs\_\texttt{nnet}} && \multicolumn{2}{c}{$ 10^{-5}$} && \multicolumn{2}{c}{0.128} && \multicolumn{2}{c}{0.037}&& \multicolumn{2}{c}{0.217}\\
\multicolumn{2}{c}{} && \multicolumn{2}{c}{$ 10^{-7}$} && \multicolumn{2}{c}{0.127} && \multicolumn{2}{c}{0.049}&& \multicolumn{2}{c}{0.373}\\

\cline{1-14} 
\multicolumn{2}{c}{} && \multicolumn{2}{c}{$ 10^{-3}$} && \multicolumn{2}{c}{0.142} && \multicolumn{2}{c}{0.078}&& \multicolumn{2}{c}{0.166}\\
\multicolumn{2}{c}{\texttt{nnet}} && \multicolumn{2}{c}{$ 10^{-5}$} && \multicolumn{2}{c}{0.134} && \multicolumn{2}{c}{0.049}&& \multicolumn{2}{c}{0.178}\\
\multicolumn{2}{c}{} && \multicolumn{2}{c}{$ 10^{-7}$} && \multicolumn{2}{c}{0.133} && \multicolumn{2}{c}{0.061}&& \multicolumn{2}{c}{0.333}\\
\cline{1-14} 
\end{tabular}
\end{center}
\end{table}


\begin{table}
\caption{Averaged (over 200 replicates) in-sample (training) and out-of-sample (testing) mis-classification rate (MR), standard deviation (SD) and running time (seconds) with $\hat x$ computed in-sample (with Training Accuracy $\varepsilon=10^{-6}$) using SMCO-R/-BR (with parallel), \texttt{nnet} and SMCOfs\_\texttt{nnet} in Example \ref{example-4}. SMCO-R/-BR run with 50 diagonal multi-starting points, maximum iterations $=300$.}
\centering
\scalebox{1}{
\begin{tabular}{lcccc}
\toprule
Method & Dataset & MR & SD & Time (s) \\
\hline
\midrule
\multirow{2}{*}{\textbf{SMCO-R}}   
  & Test  & 0.1132 & 0.0154 & \multirow{2}{*}{25.413} \\
 & Train & 0.0943 & 0.0085 &                       \\
\midrule
\hline
\multirow{2}{*}{\textbf{SMCO-BR}}  
  & Test  & 0.1130 & 0.0156 & \multirow{2}{*}{30.422} \\
 & Train & 0.0943 & 0.0084 &                        \\
\midrule
\hline
\multirow{2}{*}{nnet}
 & Test  & 0.1281 & 0.0450 &\multirow{2}{*}{0.172} \\
  & Train & 0.0985 & 0.0460 &                         \\
\midrule
\hline
\multirow{2}{*}{SMCOfs\_nnet}
  & Test  & 0.1266 & 0.0413 &    \multirow{2}{*}{0.168} \\
  & Train & 0.0990 & 0.0409 &                     \\
\bottomrule
\hline
\end{tabular}}
\end{table}

\end{document}